\documentclass{emsprocart}
\usepackage{stmaryrd}
\usepackage[mathscr]{eucal}
\usepackage{amscd}
\usepackage[all]{xy}
\usepackage{mathrsfs}
\usepackage{xypic}
\usepackage{latexsym}
\newcommand{\ff}{\footnote}

\contact[igordon@ed.ac.uk]{I. Gordon, School of Mathematics and Maxwell Institute for Mathematical Sciences,
 Edinburgh  University, Edinburgh EH9 3JZ, Scotland.}

\newtheorem*{theorem}{Theorem}
\newtheorem*{corollary}{Corollary}
\newtheorem*{lemma}{Lemma}
\newtheorem*{proposition}{Proposition}
\newtheorem*{conjecture}{Conjecture}

\newtheorem*{definition}{Definition}

\theoremstyle{definition}

\newtheorem*{remark}{Remark}

\newtheorem*{question}{Question}

\newcommand{\Lmod}[1]{#1\text{-}{\mathsf{mod}}}

\newcommand{\Lgmod}[1]{#1\text{-}{\mathsf{grmod}}}
\newcommand{\tors}[1]{#1\text{-}{\mathsf{tors}}}

\newcommand{\Ext}{\mathrm{Ext}}
\newcommand{\Hom}{\mathrm{Hom}}
\newcommand{\End}{\mathrm{End}}
\newcommand{\im}{\mathrm{Im}}

\newcommand{\gr}{\mathrm{gr}\,}
\newcommand{\Ann}{\mathrm{Ann}}

\newcommand{\Tr}{\mathrm{Tr}}

\newcommand{\proj}{\mathsf{Proj}}
\newcommand{\ad}{\mathrm{ad}}
\newcommand{\modu}{\mathrm{mod}\,}

\newcommand{\erem}{\hfill$\lozenge$\end{remark}}
\newcommand{\beq}{\begin{equation}\label}
\newcommand{\eeq}{\end{equation}}

\newcommand{\f}[1]{\mathfrak{#1}}
\newcommand{\scr}[1]{\mathscr{#1}}
\newcommand{\Spec}{\mathrm{Spec} \,}

\newcommand{\GL}{GL}

\newcommand{\vi}{${\sf {(i)}}\;$}
\newcommand{\vii}{${\sf {(ii)}}\;$}

\newcommand{\g}{{\mathfrak{g}}}

\newcommand{\irr}{{\mathsf{Irrep}}}

\renewcommand{\part}{{\f P}}

\def\C{{\mathbb{C}}}

\def\E{{\mathsf{E}}}

\def\gl{{\mathfrak{g}\mathfrak{l}}}

\newcommand{\KZ}{{\mathsf{KZ}}}

\newcommand{\h}{{\mathfrak{h}}}

\newcommand{\D}{{\scr D}}

\newcommand{\B}{{\mathsf{B}}}

\newcommand{\lf}{\lfloor}
\newcommand{\rf}{\rfloor}

\newcommand{\Z}{{\mathbb{Z}}}

\newcommand{\Id}{\operatorname{Id}}
\newcommand{\Hilb}{{\operatorname{Hilb}^n{\mathbb{\C}}^2}}
\newcommand{\hil}{\operatorname{Hilb}}

\newcommand{\N}{{\mathbb{N}}}

\newcommand{\la}{{\lambda}}
\newcommand{\Mat}{{\mathrm{Mat}}}
\newcommand{\QQ}{{\overline{Q}}}
\newcommand{\Rep}{{\mathrm{Rep}}}

\newcommand{\Coh}{{\operatorname{Coh}}}

\newcommand{\ve}{\varepsilon}

\newcommand{\sym}{{\sf Sym}}
\newcommand{\fix}{{\sf Fix}}
\newcommand{\co}[2]{\C [#1]^{\text{co} #2}}
\newcommand{\defn}{\stackrel{\mbox{\tiny{\sf def}}}=}
\newcommand{\vr}{V_{\sf reg}}
\newcommand{\hr}{\h_{\sf reg}}
\newcommand{\codim}{\operatorname{codim}}
\newcommand{\id}{\operatorname{id}}
\newcommand{\blambda}{\boldsymbol{\lambda}}

\newcommand{\otot}{\otimes \cdots \otimes}

\newcommand{\diag}{{\sf diag}}

\newcommand{\coinv}{{\sf co}}
\newcommand{\triv}{{\sf triv}}
\newcommand{\sign}{{\sf sign}}
\newcommand{\tor}{{\sf tor}}
\newcommand{\Deltar}{\Delta^{\sf reg}}
\newcommand{\rank}{{\sf rank}}
\newcommand{\vect}{{\sf Vect}}
\newcommand{\bsupp}{{\mathrm{\bf Supp}}}

\title{Symplectic Reflection Algebras}

 \author[Iain G. Gordon]{Iain G. Gordon \thanks{I would like to thank very much Gbenga Akinbo for his help in preparing the first version of these notes; I am grateful to the ICRA committee for encouraging me to submit this paper to the ICRA XII book.} 
} 

\begin{document}

\begin{abstract}
We survey recent results on the representation theory of symplectic reflection algebras, focusing particularly on connections with symplectic quotient singularities and their resolutions, with category $\mathcal{O}$, and with spaces of representations of quivers. 
\end{abstract}

\begin{classification}
Primary 16G.
\end{classification}

\begin{keywords}
Symplectic reflection algebra, rational Cherednik algebra, symplectic resolution, quiver variety, differential operators.
\end{keywords}

\maketitle
\maketitle
\section*{Introduction} 

Apart from confirming conjectures in algebraic combinatorics, integrable systems and real algebraic geometry, having an interesting representation theory, and connecting with noncommutative, quiver, Hall and Hopf algebras, all symplectic reflection algebras have given us is an algebraic approach to resolutions of symplectic quotient singularities. They are new structures in representation theory which access parts of  ``double" algebra and ``double" geometry, but build largely on simple Lie algebras, preprojective algebras and deformation theory. Although many basic questions remain unanswered, their surprisingly diverse applications to a number of topics encourages more detailed investigation by more people from more fields. 

This article is a discursive introduction to symplectic reflection algebras. It attempts to explain, from a mostly algebraic point of view, the why rather than the how. The how can be found in every paper about symplectic reflection algebras, beginning with the groundbreaking, deep, beautiful, long paper of Etingof and Ginzburg that introduced symplectic reflection algebras, \cite{EG}: if you really want to learn about symplectic reflection algebras in detail there is no more inspiring place to start. 

Here we will introduce symplectic reflection algebras as deformations of orbit space singularities, echoing Crawley-Boevey and Holland's study of deformed preprojective algebras as deformations of kleinian singularities, and we will use this point of view as motivation for many of the constructions and results we present. These include the classification of symplectic singularities admitting symplectic resolutions; category $\mathcal{O}$, the $\KZ$-functor, highest weight covers and finite Hecke algebras; derived equivalences for some quiver varieties; the construction of quantisations of Hilbert schemes of points on the plane (by differential operators related to quiver varieties, by noncommutative algebraic geometry, and by microlocalisation) and associated geometric interpretations of representations of symplectic reflection algebras. 

\medskip
Several surveys on aspects of symplectic reflection algebras already exist, but they have different perspectives: Brown's on ring-theoretic aspects, \cite{kenny}, Etingof and Strickland's on quasi-invariants, integrable systems and differential operators, \cite{etst}, Rouquier's on rational Cherednik algebra representation theory, \cite{rou}, Etingof's lectures on deformation theory, \cite{Edef}, and on Calogero--Moser systems, \cite{ECM}, as well as Cherednik's book on the dahaddy of symplectic reflection algebras, \cite{cherednikbook}. In particular, since we deal with symplectic reflection algebras rather than just rational Cherednik algebras, quivers will play a more prominent r\^ole here than in the other surveys.
 
 \medskip
There are many exciting parts of the theory that we will not touch. These include positive characteristic,  \cite{BFG}, \cite{latour}, \cite{BN}; continuous Hecke algebras, \cite{egg}, \cite{mon3}; generalisations of Schur--Weyl duality, \cite{guay2}, \cite{guay3}, \cite{guay4}; connections with affine Lie algebras, \cite{arakawasuzuki}; Cherednik algebras for varieties with group actions, \cite{et}; relations with Hopf algebras, \cite{bazber}; real algebraic geometry, \cite{gordshapiro}; symmetric function theory, \cite{berbur}; and most of all double affine Hecke algebras, \cite{cherednikbook} and the references in it and \cite{Suzuki}, (with a new relation to Hall algebras -- \cite{schiff0} and \cite{schiff}). To try compensate for these omissions the bibliography includes references to all the papers involving symplectic reflection algebras that I know about. As well as dealing with the omitted topics, these fill in the missing how. 

\medskip
The article is divided into the following sections.

\begin{enumerate}
\item
Symplectic group actions.
\item Symplectic reflection algebras.
\item Dimension 2.
\item Basic properties.

\item Specific case: rational Cherednik algebras I ($t\neq 0$).

\item Specific case: rational Cherednik algebras II ($t=0$). 

\item Specific case: quivers and hamiltonian reduction.
\item Very specific case: the symmetric group and Hilbert schemes. 

\item Problems. 
\end{enumerate}

The content of the article is as follows. The first section introduces symplectic singularities, a class of varieties which motivate the introduction of symplectic reflection algebras. Section 2 presents a little background from deformation theory which is then used to construct symplectic reflection algebras. (If you want, it is possible to skip all of this background and start reading from the definition of a symplectic reflection algebra at \ref{PBWthm}.)  Section 3 recalls work on deformed preprojective algebras associated to affine Dynkin quivers: these are the two-dimensional case of symplectic reflection algebras. Section 4 explains a few properties  which hold for all symplectic reflection algebras: these mostly concern the case $t=0$ where the centre of the algebra is very large. Following this, the theory splits into two cases depending on the group defining the symplectic reflection algebra. In the first case the group is a complex reflection group and symplectic reflection algebras are called rational Cherednik algebras: in Section 5 we take $t\neq 0$ and introduce category $\mathcal{O}$ and the $\KZ$-functor, and describe some finite dimensional representations; in Section 6 we let $t=0$ and make connections to symplectic resolutions and conjecturally to Rouquier families. In Section 7 we discuss the second case: the group is related to an affine Dynkin quiver and we present relatively conclusive results on representation theory for $t=0$ and then relate symplectic reflection algebras to differential operators on quivers for $t=1$. It turns out that the two cases above are not disjoint -- they both include $S_n$ -- and so we are able to prove more by combining both approaches. We explain this in Section 8 where we relate representations of the symplectic reflection algebra associated to $S_n$ to sheaves on the Hilbert scheme of points on the plane. We close with a list of problems.

\section{Symplectic group actions}

We introduce a particular type of orbit singularity associated to symplectic group actions on complex vector spaces. It will turn out that this type of singularity will have a wonderfully rich structure with connections to many parts of representation theory and beyond. We ask a number of questions about these singularities which mirror classical theorems on smooth orbit spaces; these will motivate in part the introduction of symplectic reflection algebras. 

\subsection{Orbit spaces} 
Let $(G,V)$ be the data of a finite group $G$ acting linearly on a finite dimensional complex vector space $V$ over $\C$. In what follows the kernel of the action will be irrelevant, so we will always assume that $$\text{$G$ acts on $V$ faithfully.}$$ Let $\C [V]$ denote the algebra of regular functions on $V$: it is isomorphic to the symmetric algebra $\sym (V^*)$. This algebra inherits an action of $G$ via $({}^gf)(v) = f({}^{g^{-1}}v)$ for all $f\in \C[V], g\in G, v\in V$. We define the {\it invariant ring} to be $$\C[V]^G = \{ f\in \C[V] : {}^g f = f \text{ for all }g\in G\}.$$ This is a finitely generated $\C$-algebra which is by definition the ring of regular functions on the orbit space $V/G$. Under the correspondence between commutative algebra and algebraic geometry, the inclusion $\C[V]^G \hookrightarrow \C[V]$ corresponds to the orbit map $\rho : V \longrightarrow V/G$ which sends an element $v\in V$ to its $G$-orbit.

\subsection{$\!\!\!\!\!$}
$V/G$ is an irreducible affine variety and we would like to know about its geometric properties. 
\subsection{$\!\!\!\!\!$}
We begin with smoothness. Recall that $(G,V)$ is {\it generated by (complex) reflections} if there exists a generating set $S \subset G$ such that for every $s\in S$ we have $\fix(s) = \{ v\in V: {}^s v = v \}$ has codimension $1$ in $V$. In other words if we were to diagonalise the action of $s$ on $V$ there would be exactly one non-trivial eigenvalue.
The following theorem is a jewel in the crown of invariant theory with a dazzling array of applications across many fields.
\begin{theorem}[Shephard-Todd, see {\cite[Theorem 7.2.1]{BenINV}}]\label{STC}
The following are equivalent:
\begin{enumerate}
\item the orbit space $V/G$ is smooth;
\item the orbit map $\rho : V\longrightarrow V/G$ is flat;
\item the invariant ring $\C[V]^G$ is a polynomial algebra (on $\dim V$ generators);
\item $(G,V)$ is generated by complex reflections.
\end{enumerate}
\end{theorem}  
It is helpful to provide a gloss on the flatness in (b). Since we are assuming that $G$ acts faithfully on $V$, the generic $G$-orbit has $|G|$ elements and so the generic fibre of $\rho$ consists of $|G|$ points. By semi-continuity $|G|$ then provides a lower bound for the vector space dimension of any (scheme-theoretic) fibre of $\rho$, and flatness simply means that all the fibres of $\rho$ actually have this dimension. The fibre of the zero orbit is particularly interesting. By definition its coordinate ring is the {\it coinvariant algebra} of $(G,V)$ \begin{equation}\label{coinvariant} \co{V}{G} \defn \frac{\C[V]}{\langle\C[V]_+^G\rangle}\end{equation} where $\langle C[V]_+^G\rangle$ denotes the ideal of $\C[V]$ generated by all $G$-invariant regular functions with zero constant term. This special fibre inherits structure from $V$: an action of the group $G$; a grading, or equivalently a $\C^*$-action. Flatness and the rigidity of $G$-modules implies that $\co{V}{G}$ carries the regular representation of $G$ since the generic fibre does: the combinatorics of this space, and also of $\C[V]^G$, is  remarkable; we mention that if $(G,V)$ is the Weyl group of a semisimple Lie algebra acting on the Cartan subalgebra then the coinvariant ring is isomorphic as a $G$-equivariant graded algebra to the cohomology ring of the corresponding flag manifold with complex coefficients. 

\subsection{$\!\!\!\!\!$} When studying singularities of $V/G$, Part (3) of the above theorem shows that we might as well as assume that $G$ contains no complex reflections, since we can replace $G$ by the quotient $G/H$ where $H$ is the (obviously normal) subgroup of $G$ generated by complex reflections: $\C[V]^G = (\C[V]^H)^{G/H}$ and $\C[V]^H$ is a polynomial ring. We will say that $G$ is {\it small} if it contains no complex reflections. 

\subsection{$\!\!\!\!\!$} In general there are only a few results about the singularities of $V/G$: they are normal, \cite[Proposition 1.1.1]{BenINV}, (so smooth in codimension 1) and by Hochster--Eagon they are Cohen-Macaulay, \cite[Theorem 4.3.6]{BenINV}, (so a vector bundle over an affine space). Finally,

\begin{theorem}[Watanabe, see {\cite[Theorem 4.6.2]{BenINV}}] \label{watanabe} Suppose that $G$ is small. Then $\C[V]^G$ has finite injective dimension (i.e. is Gorenstein) if and only if $G \leq SL(V)$.
\end{theorem}

\subsection{Symplectic singularities and resolutions} 
We now move attention to a special class of pairs $(G,V)$ in order to say something more explicit. Henceforth we will assume that $V$ is a complex symplectic vector space with symplectic form $\omega_V$ such that $G$ preserves $\omega_V$, meaning that $$\omega_V (v_1,v_2) = \omega_V ({}^gv_1, {}^gv_2) \quad \text{for all }g\in G, v_1, v_2\in V.$$ A good example to keep in mind is to begin with $(G, \h)$ and then produce $V = \h\oplus \h^*$ with $\omega_V ((x_1,f_1), (x_2,f_2)) = f_2(x_1) - f_1(x_2)$. Here $G$ acts on naturally on $V$ with $G\leq Sp (V,\omega)$. This is a special case of a geometric construction: start with a smooth $G$-variety $X$; then $G$ acts on the cotangent bundle $T^*X$ which has a canonical symplectic structure

\subsection{$\!\!\!\!\!$} We will refer to the data $(G,V, \omega_V)$ as a {\it symplectic triple}. In what follows we can assume without loss of generality that $(G,V,\omega_V)$ is an {\it indecomposable symplectic triple} which means that there is no $G$-equivariant splitting $V= V_1 \oplus V_2$ with $\omega (V_1, V_2) = 0$. Clearly any symplectic triple is a direct sum of indecomposable triples.

\subsection{$\!\!\!\!\!$} \label{sympres} Observe that if $(G,V,\omega_V)$ is a symplectic triple then $G$ is a small group, for $G\leq Sp(V) \leq SL(V)$ but by definition any reflection has non-trivial determinant. Furthermore for any $g\in G$ we can make a $\langle g \rangle$-equivariant splitting $V = \fix(g) \oplus \overline{\fix}(g)$ and then see that $\omega (\fix (g), \overline{\fix}(g))=0$. It follows that $\omega$ restricts to a non-degenerate form on both $\fix(g)$ and $\overline{\fix}(g)$ and so $\fix(g)$ must have codimension at least two in $V$ if $g\neq 1$. Now set $$\vr \defn V\setminus \bigcup_{g\neq 1} \fix(g),$$ a $G$-equivariant open subset of $V$ of codimension at least $2$ which inherits a $G$-equivariant complex symplectic form from $V$. Since $(G,V)$ is {\it not} a complex reflection group, the orbit space $V/G$ is singular: however $G$ acts freely on $\vr$ and we see that $\vr /G$ is $(V/G)_{\sf sm}$, the smooth locus of $V/G$. Since the symplectic form on $\vr$ was $G$-equivariant we can push it down to $\vr /G$ to induce a symplectic form on $(V/G)_{\sf sm}$ which we denote by $\omega_{(V/G)_{\sf sm}}$. This leads to a key definition.

\begin{definition} 
A {\it symplectic resolution} of $V/G$ is a resolution of singularities $\pi : X \longrightarrow V/G$ such that there exists a complex symplectic form $\omega_X$ on $X$ for which the isomorphism $$\pi|_{\pi^{-1}((V/G)_{\sf sm})} : \pi^{-1}((V/G)_{\sf sm}) \longrightarrow (V/G)_{\sf sm}$$ is a symplectic isomorphism, i.e. $\pi^{\ast} (\omega_{(V/G)_{\sf sm}})$ equals the restriction of $\omega_X$ to the open set $\pi^{-1}((V/G)_{\sf sm})\subset X$.
\end{definition}
There are several useful comments to make around this definition. 
\begin{itemize}
\item In \cite{beauville} Beauville introduced the weaker notion of a {\it symplectic singularity}: by definition $Y$ is a symplectic singularity if the smooth locus $Y_{\sf sm} \subseteq Y$ carries a symplectic form $\omega$ and that for some (thus any) resolution of singularities $\pi: X\longrightarrow Y$, the pull-back $\pi^{\ast} (\omega)$ defined on $\pi^{-1}(Y_{\sf sm})\subseteq X$ extends to a (possibly degenerate) $2$-form on all of $X$. \cite[Proposition 2.4]{beauville} shows that any $V/G$ as above is a symplectic singularity. There are other very interesting affine examples coming from closures of nilpotent orbits of simple Lie algebras which have been studied extensively by Fu, \cite{Funilp}, and have beautiful representation theory associated to them: an example is the Springer resolution of the nullcone of a simple Lie algebra $\pi : T^*\mathcal{B} \longrightarrow \mathcal{N}$. For a survey see \cite{kaledinsurvey}.
\item It is not true that any $V/G$ admits a symplectic resolution: that classification of such $(G,V)$ is (almost) complete thanks to the representation theory of symplectic reflection algebras.  
\item When they exist, symplectic resolutions of $V/G$ need not be unique. It is conjectured, however, that there should be only finitely many non-isomorphic symplectic resolutions of any given symplectic variety, \cite[Conjecture 1]{FN}.\item A theorem of Fu, Kaledin and Namikawa -- see the survey \cite[Proposition 1.6]{Fusurvey} -- shows that $\pi : X\longrightarrow V/G$ is a symplectic resolution if and only if it is a crepant resolution, i.e. if and only if $\pi^{\ast} K_{V/G} \cong K_X$, where $K_{V/G}$ and $K_X$ denote the canonical bundles of $V/G$ and $X$ respectively. Note that $K_{V/G}$ is trivial since $V/G$ is Gorenstein by Theorem \ref{watanabe}. So $\pi$ is crepant if and only if $K_X$ is trivial, in other words if and only if $X$ is a Calabi-Yau variety. One direction of this equivalence is clear since if $X$ is symplectic with form $\omega_X$ then $\wedge^{{\sf top}/2}\omega_X$ trivialises $K_X$. Crepant resolutions and Calabi-Yau varieties are important concepts in algebraic geometry and mathematical physics. However, much of the focus has been in three dimensions -- the examples which we discuss here are not three dimensional since the symplectic structure implies that $V/G$ is even dimensional.  
 \end{itemize}
 \subsection{$\!\!\!\!\!$} Motivated by Theorem \ref{STC} and the above discussion we ask the following questions for indecomposable symplectic triples $(G,V,\omega_V)$.
 \begin{question}
 When does $V/G$ admit a symplectic resolution? What obstructions are there?
 \end{question}
  \begin{question}
 How does the coinvariant ring $\co{V}{G}$ behave?
 \end{question}
  \begin{question}
 Is there interesting representation theory attached to $V/G$ and its resolutions?
 \end{question}
  \begin{question}
 Are there interesting combinatorics attached to $V/G$?
 \end{question}
 To some extent these questions have been part of the motivation for the research in symplectic reflection algebras and they will play that r\^ole here too. It turns out, however, that symplectic reflection algebras are not bound to the world of symplectic singularities: as we will see there are great many other uses.
 
 \subsection{Toy example} We will consider the example throughout the survey: $G = \mu_2$ acting on $V = \C \oplus \C^*$ by multiplication by $-1$ with $\omega_V$ be the standard symplectic form. Letting $\C[V] = \C[x,y]$ we see that $\C[V]^G = \C[x^2, xy, y^2] = \C[A,B,C] / (AC-B^2)$, the quadric cone. This has an isolated singularity at the origin, i.e. at the zero orbit, which can be resolved by blowing up there. The resulting resolution $\pi: T^*\mathbb{P}^1 \longrightarrow V/G$ collapses the zero fibre of $T^*\mathbb{P}^1$ to a point: it is a symplectic resolution where $T^\ast\mathbb{P}^1$ has its canonical symplectic structure. There is an almost endless amount of interesting representation theory attached to this.

\subsection{Algebra} 
 It is easy to explain why the restriction to symplectic triples $(G,V, \omega_V)$ is of interest from the algebraic point of view. The non-degeneracy of $\omega_{\vr/G}$ allows us to identify $T^*(\vr/G)$ and $T(\vr /G)$ and hence to identify $k$-forms with $k$-vectors. Thus $\omega_{\vr/G}$ corresponds to some $2$-vector $\Theta \in \wedge^2 T(\vr/G)$. Since $V/G$ is normal and $\vr/G$ has codimension $2$ in $V/G$ this $2$-vector can be extended to $V/G$ by Hartog's theorem; we continue to call it $\Theta$. This encodes the data of a Poisson bracket on $\C [V]^G$ via $\{ f_1, f_2 \} = \Theta (df_1, df_2)$: the Jacobi identity is equivalent to the fact that $d\omega_{\vr/G} = 0$.
(There is another more down-to-earth description of this $2$-vector. 
The symplectic form on $V$ induces a Poisson bracket on $\C[V]$ as follows. Let $v_1, \ldots ,v_n$ be a basis of $V$ and let $x_1, \ldots , x_n$ be the dual basis, so that $\C [V] = \C[x_1, \ldots , x_n]$. Then $$\{ f, g \} = \sum_{1\leq i, j \leq n} \frac{\partial f}{\partial x_i}\frac{\partial g}{\partial x_j} \omega_V (v_i, v_j).$$
Since $\omega$ is $G$-invariant the bracket restricts to $\C[V]^G$ and this is the same as the one above described by means of $\Theta$.)  

So whenever we have a symplectic triple $(V,G,\omega_V)$ we have an induced Poisson bracket on $\C[V]^G$. Poisson brackets are the residue of deformations of $\C[V]^G$, and so we are immediately lead to noncommutative algebras deforming $\C[V]^G$. We hope will help us understand the finer structure of the orbit space $V/G$. 

\section{Symplectic reflection algebras}

In this chapter, after discussing a little relevant deformation theory, we introduce symplectic reflection algebras.

\subsection{Generalities on deformations}
Throughout this section $k$ will denote a semisimple
artinian $\C$-algebra and $A$ a $k$--algebra, i.e. a $k$--bimodule with
a $k$--bimodule mapping $A\otimes_k A\longrightarrow A$.

\begin{remark} This is an {\it unusual} definition, but it's here for a serious
reason. We are going to study two types of deformations -- formal
and graded. While formal deformations only take place in an infinitesimal
neighbourhood of the algebra, graded deformations take place along the
$k$-affine line. So although $k$ is a semisimple algebra, and so has
trivial Hochschild cohomology and therefore no formal
deformations, it is not true that it has no global deformations.
We want to avoid this happening: when
we deform $A$ we won't want to
deform $k$. To ensure this means using the above definition
of a $k$--algebra.
\end{remark}

\subsection{Formal deformations}
Recall that a {\it formal deformation of $A$} is a $k[[\hbar]]$-bimodule map $\star : A[[\hbar]] \times A[[\hbar]] \longrightarrow A[[\hbar]]$ that makes $A[[\hbar]]$ a $k[[\hbar]]$--algebra with $\hbar$ central, and that deforms the multiplication on $A$ in the sense that $a\star b \equiv ab  \text{ mod } \hbar A[[\hbar]]$ for all $a,b\in A$. Similarly, an {\it $i$th level deformation of $A$} is the same as above, replacing $k[[\hbar]]$ with the truncated polynomial ring $k[\hbar]/(\hbar^{i+1})$. A first level deformation of $A$ is also called an {\it infinitesimal deformation of $A$}.

Let $J$ be the group of $k[[\hbar]]$--bimodule automorphisms $g$ of $A[[\hbar]]$ such that $g (u) \equiv u \text{ mod } \hbar A[[\hbar]]$ for all $u\in A[[\hbar]]$. Two formal deformations $\star$ and $\star'$ are said to be {\it equivalent} if there is an element $g\in J$ such that $g(u\star v) = g(u) \star' g(v)$ for all $u,v\in A[[\hbar]]$. 

\subsection{$\!\!\!\!\!$}
In a deformation we can write the product of two elements $a,b\in A$ as
$$ a\star b = ab + B_1(a,b)\hbar + B_2(a,b)\hbar^2 + \cdots + B_i(a,b)\hbar^i + \cdots $$
for bimodule mappings $B_i : A \otimes A \longrightarrow A$, and these mappings determine the multiplication $\star$. 
If we let $H^{\bullet}(A,A)$ denote the Hochschild cohomology of $A$ then it is an elemntary exercise to check the following.
\begin{itemize}
\item The set of isomorphism classes of infinitesimal deformations of $A$ canonically identifies with $H^2(A,A)$.
\item Given $A_i$, an $i$th level deformation of $A$, the obstruction for its continuation to the $(i+1)$--st level lies in $H^3(A,A)$.
\item Let $A_i$ be as above. Then the set of isomorphism classes of continuations of $A_i$ to the $(i+1)$--st level is an $H^2(A,A)$--homogeneous space.
\end{itemize}

In particular, if $H^3(A,A) = 0 $ then any infinitesimal deformation on extends to a formal deformation of $A$.

\subsection{Poisson brackets} \label{defpoisson} If we assume that $A$ is a commutative $k$-algebra then there is a Poisson bracket associated to any deformation: for $a,b \in A$ we define \begin{equation}\label{getPoisson} \{ a, b \} \defn \left( \hbar^{-1} (a\star b - b\star a)\right)_{|_{\hbar=0}}. \end{equation} Since $A$ is commutative $a\star b - b\star a \in \hbar A[[\hbar]]$, so this definition does make sense. Of course, this Poisson bracket may be trivial, for instance if $B_1(a,b) = 0 $ for all $a,b \in A$. In this case we can attempt to modify the definition by defining $m = \min \{ i: \text{ there exist $x,y\in A$ with $B_i(x,y) \neq B_i(y,x)$} \}$ (if this exists) and then setting $\{a,b\}   \defn \left( \hbar^{-m} (a\star b - b\star a)\right)_{|_{t=0}}$ for all $a,b\in A$. If the above minimum does not exist we set $m=\infty$ and we simply take the trivial bracket. In this guise, the Poisson bracket on $A$ is a residue of the noncommutativity of the deformation $(A[[\hbar]], \star)$; in particular $(A[[\hbar]], \star)$ is a commutative deformation if and only if the Poisson bracket is trivial.
\subsection{Graded deformations} \label{gradeddef}
Suppose that $A$ is an $\N$-graded algebra $A = \oplus_{j\geq 0} A_j$. A {\it graded deformation of $A$} is a $k[\hbar]$--bimodule map $$\star : A[\hbar] \times A[\hbar] \longrightarrow A[\hbar]$$ that makes $A[\hbar]$ an $\mathbb{N}$-graded $k[\hbar]$--algebra where $\deg \hbar = 1$. We can write the product of two elements $a,b\in A$ as
$$ a\star b = ab + B_1(a,b)\hbar + B_2(a,b)\hbar^2 + \cdots + B_i(a,b)\hbar^i + \cdots $$
for bimodule mappings $B_i : A \otimes A \longrightarrow A$ of degree $-i$ and these mappings determine the multiplication $\star$. Isomorphisms of graded $k[\hbar]$--algebras give rise to equivalent deformations.

\begin{remark}
It makes perfect sense to deal with $k[\hbar]$-deformations rather than $k[[\hbar]]$--deformations since any graded $k[[\hbar]]$-deformation of $A$ actually comes from a $k[\hbar]$--deformation. For if we take $u,v\in A[\hbar]\subset A[[\hbar]]$ then we can write $u = \sum_{i=0}^n u_i \hbar^i$ and $v = \sum_{i=0}^n v_i \hbar^i$ for some large enough $n$ and for some $u_i, v_i \in A$. Now we can find a positive integer $m$ such for all $i$ we have $u_i,v_i \in \oplus_{j=0}^m A_j$. Since the formal deformation of $A$ is graded and the degree of $B_l$ is $-l$ we see that $B_l(u_i,v_j) = 0$ for all $l > 2m$. Therefore $u\star v \in A[\hbar]$ and so $(A[\hbar], \star)$ is a subalgebra of $(A[[\hbar]], \star)$. The advantage of graded deformations is that we can specialise $\hbar$ to {\it any} value, i.e. the deformation is defined over $\mathbb{A}^1_k$ rather some formal neighbourhood of zero.
\end{remark}

\subsection{$\!\!\!\!\!$} If $A$ is a graded algebra then we can construct graded Hochschild cohomology groups in the category of graded $A$-bimodules and these groups are themselves graded $\hat{H}^{\bullet}(A,A) = \bigoplus_{j \geq 0} \hat{H}^{\bullet}_j(A,A).$ The same argument as in the ungraded case then shows the following.

\begin{itemize}
\item The set of isomorphism classes of infinitesimal graded deformations of $A$ canonically identifies with $\hat{H}^2_{-1}(A,A)$.
\item Given $A_i$, an $i$th level deformation of $A$, the obstruction for its continuation to the $(i+1)$--st level lies in $\hat{H}^3_{-i-1}(A,A)$.
\item Let $A_i$ be as in (2). Then the set of isomorphism classes of continuations of $A_i$ to the $(i+1)$--st level is an $\hat{H}^2_{-i-1}(A,A)$--homogeneous space.
\end{itemize}

\subsection{Koszul deformation principle} Koszul algebras are important for many reasons:
here we will exploit the fact that they have relatively easy
deformation theory, \cite{BrGa}.
\begin{definition}
A graded $k$--algebra $A$ is {\it Koszul} if $k$ considered as a left $A$--module by $k = A/A_{>0}$ has a graded projective resolution $\cdots \longrightarrow P^2 \longrightarrow P^1 \longrightarrow P^0 \longrightarrow k \longrightarrow 0$ such that $P^i$ is generated by its component of degree $i$ (i.e. $P^i = AP^i_i$).
\end{definition}

The {\bf key example} for us will be $A=\C [V]\rtimes G$, the smash product of $\C[V]$ and $G$, considered as a $k=\C G$-algebra. As a vector space $\C[V]\rtimes G$ is isomorphic to $\C[V] \otimes_{\C} \C G$ and both $\C[V]$ and $\C G$ are $\C$-subalgebras of $\C[V]\rtimes G$, but the multiplication between $\C[V]$ and $\C G$ is twisted to take account of the action of $G$ on $V$, namely $g\cdot p = {}^gp\otimes g$ for any $p\in \C[V]$ and $g\in G$. The grading on $\C[V]\rtimes G$ is given by putting $G$ in degree $0$ and $V^*\subset \C[V]$ in degree $1$. To see that $\C[V]\rtimes G$ is Koszul we use the obvious generalisation of the Koszul resolution for $\C[V]$  $$\cdots \longrightarrow \C[V]\otimes_{\C} \wedge^p V^*
\otimes_{\C} \C G \longrightarrow \cdots \longrightarrow
\C[V]\otimes_{\C} \wedge^1 V^* \otimes_{\C} \C G
\longrightarrow \C[V]\otimes_{\C} {\C}G \longrightarrow
 \C G \longrightarrow 0,$$ where $G$ acts on each term diagonally from the left.

\subsection{$\!\!\!\!\!$}
The critical lemma is the following.
\begin{lemma} Let $A$ be a Koszul $k$-algebra. Then for all $p < -q$ we have $\hat{H}^p_q (A,A) = 0$.
\end{lemma}
The proof of this is quite straightforward. Normally one uses the bar resolution to calculate Hochschild cohomology, but in this situation it is possible to cook a smaller complex from the given projective resolution of $k$ which can also be used. Now the restrictions on the gradings in the Koszul complex translate to the stated vanishing of cohomology.

Thanks to our discussion on deformation theory, this means that quadratic Koszul $k$-algebras should have a controllable deformation theory.

\subsection{Graded and filtered algebras}
The key example of a Koszul $k$-algebra $\C[V]\rtimes G$ is a particular case of the following general construction. Let $W$ be a $k$--bimodule and $T_k(W)$ the tensor algebra of $W$ over $k$. Let $P \subseteq T_k^{\leq 2}(W) = k \oplus W \oplus W\otimes W$ and let $J(P)$ denote the two--sided ideal of $T_k(W)$ generated by $P$. The algebra $Q(W,P) = T_k(W)/J(P)$ is called a {\it nonhomogeneous quadratic algebra}. If $R\subseteq W\otimes W$ then $Q(W,R)$ is called a {\it quadratic algebra}. Note that quadratic algebras are graded with $Q(W, R)_m = W^{\otimes m}/(J(R) \cap W^{\otimes m}).$ 

To get $\C[V]\rtimes G = Q(W, R)$ take $k = \C G$, $W = V^*\otimes \C G$ with $G$ acting diagonally on the left and by right multiplication on the right, and $R$ to be the $k$-span of the elements $(x\otimes_{\C} 1) \otimes (y \otimes_{\C}
g) - (y\otimes_{\C} 1) \otimes (x \otimes_{\C} g) \in W\otimes_k W$ with $x, y\in V^*$ and $g\in G$.

\subsection{$\!\!\!\!\!$} Suppose that $B = Q(W, P)$ is a nonhomogeneous quadratic
algebra. Then there is a quadratic algebra canonically associated to
$B$. To define it let $\pi : k\oplus W\oplus W\otimes W
\longrightarrow W\otimes W$ be the projection map and set $R = \pi
(P)$: we associate $A = Q(W,R)$. We will say
that $B$ is a {\it PBW deformation} of $A$ if there exists a graded
deformation $\tilde{A}$ of $A$ such that $\tilde{A}_1 \defn
\tilde{A}/(\hbar-1)\tilde{A} = B$.

There is another useful way to interpret a PBW deformation. While the algebra $B$ is not usually graded, it is filtered. This means that there is an ascending chain of $k$--bimodules $$0 = F^{-1}B  \subseteq F^{0}B \subseteq F^{1}B \subseteq F^{2}B \subseteq \cdots \subseteq F^iB \subseteq \cdots \subseteq B$$ where $F^iB$ consists of the image of elements of $T_k V$ of degree less than or equal to $i$. Note that $F^iB F^jB \subseteq F^{i+j}B$ and $k\subseteq F^0B$.  Now set $gr_i B \defn F^iB/F^{i-1}B$ and let $\gr B = \oplus_{i\geq 0} \gr_i B$. This is naturally a $k$--algebra if we define multiplication as follows: given $a = f + F^{i-1}B \in gr_i B$ and $b = g + F^{j-1}B \in \gr_jB$; then $ab = fg + F^{i+j-1}B\in \gr_{i+j}B$. This new algebra is called the {\it associated graded algebra of $B$}. Now by construction $\gr B$ is generated over $k$ by the image of $W$. Thus the universal property of $T_kW$ ensures a surjective homomorphism $\theta : T_kW \longrightarrow \gr B$. Furthermore, if $r\in R$ then there exists $p\in P$ such that $r +T_k^{\leq 1}W = p+T_k^{\leq 1}W$, from which it follows that $\theta (r) = p + F^1B = 0$. Hence there is a surjective homomorphism $$\theta : A=Q(W,R) \longrightarrow \gr B= \gr Q(W,P).$$ This is an isomorphism if and only if $B$ is a PBW deformation.

\begin{theorem}[{\cite{BrGa}}]
\label{bravermangaitsgory}
Let $A= Q(W,R)$ be a Koszul algebra with $W$ a free $k$--module (on the left and on the right) and assume we are given $\alpha : R \longrightarrow W$ and $\beta : R\longrightarrow k$. Set $P = \{ r + \alpha(r) + \beta (r) : r\in R \} \subseteq T_k^{\leq 2}W$ so that $\pi (P) = R$. Then $\theta : A \longrightarrow \gr B = \gr Q(W,P)$ is an isomorphism if and only if the following three conditions are satisfied (where the domain of each of the mappings is $(R\otimes W)\cap (W\otimes R)$):
\begin{enumerate}
\item $\alpha \otimes \text{id} - \text{id} \otimes \alpha$ has image in $R$;
\item $\alpha \circ (\alpha \otimes \text{id} - \text{id} \otimes \alpha) =   \text{id} \otimes \beta  - \beta \otimes \text{id} $;
\item $\beta \circ ( \alpha \otimes \text{id} - \text{id} \otimes \alpha) = 0.$
\end{enumerate}
\end{theorem}

\subsection{Symplectic reflection algebras}
At last we are in position to introduce symplectic reflection algebras. Our goal earlier was to construct deformations of $\C[V]^G$ where $(G,V,\omega_V)$ is an indecomposable symplectic triple. This is difficult to do, however, because we have little explicit understanding of the invariant ring $\C[V]^G$ except in some special cases. Instead our tactic will be to consider $\C[V]^G$ as the centre of the smash product $\C[V]\rtimes G$, then to deform $\C[V]\rtimes G$ -- which we can manage well because it is a Koszul $\C G$-algebra --, and finally to check that $\C[V]^G$ has deformed nicely too. Symplectic reflection algebras are actually then deformations of $\C[V]\rtimes G$; they will, however, have subalgebras deforming $\C[V]^G$.

\subsection{$\!\!\!\!\!$} Since $ \C[V]\rtimes G$ is a Koszul $k$-algebra of the form $Q(W, R)$ it makes sense to speak of PBW deformations of $\C[V]\rtimes G$. The key theorem--definition in the subject is this.

\begin{theorem}[{\cite[Theorem 1.3]{EG}}] \label{PBWthm} Let $(G,V,\omega_V)$ be an indecomposable symplectic triple. Then the PBW deformations of $\C[V]\rtimes G = \sym(V^*)\rtimes G$ are {\bf precisely} the algebras $H_{\kappa} \defn T_{\C}(V^*)\rtimes G/\langle x\otimes y - y\otimes x - \kappa(x,y): x,y\in V^*\rangle$ where  $\kappa : V^*\otimes V^* \longrightarrow \C G$ is an alternating form on $V^*$ of the form $$\kappa (x,y) = t\omega_{V^*} (x,y) -2 \sum_{s\in \mathcal{S}} c(s)
\omega_s (x,y) s .$$ Here $\omega_{V^*}$ is the symplectic form on $V^*$ corresponding to $\omega_V$ under the identification of $V$ and $V^*$ induced by $\omega_V$; $\mathcal{S}$ is the set of symplectic reflections of $(G,V, \omega_V)$, that is $\mathcal{S} \defn \{ s\in G: \codim_V \fix (s) =2 \}$; $t\in \C$ and $c\in \C[\mathcal{S}]^{\operatorname{ad} G}$ is a class function on $\mathcal{S}$; $\omega_s$ is  the restriction of $\omega_V$ to $(\id - s)(V)$ whose radical is $\ker (\id - s)$.
\end{theorem}
\begin{proof}
Thanks to the Koszul deformation principle, this is now a straightforward check of Conditions (1)-(3) of Theorem \ref{bravermangaitsgory}. The details are given in \cite{EG}, and also (in a more general context) in \cite{RS}. 
\end{proof}
\begin{remark}
For general pairs $(G,V)$ a description of PBW deformations of $\C[V]\rtimes G$ was originially given by Drinfel'd, \cite{drinfeld}. In the symplectic case this was rediscovered by Etingof and Ginzburg as above, and Drinfeld's general case was described in detail by Ram--Shepler, \cite{RS}. In this generality these algebras are sometimes called {\it
graded Hecke algebras} because, when $G$
is a Weyl group acting on its reflection representation, there are connections with
graded affine Hecke algebras which are important in the
representation theory of Lie groups and beyond, \cite{Lu}.
\end{remark}

\subsection{$\!\!\!\!\!$} Let's record the official definition of a symplectic reflection algebra.
\begin{definition} Given the data of an indecomposable symplectic triple $(G,V, \omega_V)$ and $t\in \C$, $c\in \C[\mathcal{S}]^{\operatorname{ad}G}$ as in Theorem \ref{PBWthm}, we write the corresponding PBW deformation of $\C[V]\rtimes G$ as $H_{t,c}$ and call it a {\it symplectic reflection algebra}. 
\end{definition}
\noindent
If it's necessary to specify the group, we will write $H_{t,c}(G)$.

So, by construction, $H_{t,c}$ is generated by $V^*$ and $G$ and we have a filtration on $H_{t,c}$ with $F^0 = \C G$, $F^1 = V\oplus \C G$ and $F^i = (F^1)^i$ such that there is a $\C G$-algebra isomorphism \begin{equation} \label{PBW} \C [V]\rtimes G \stackrel{\sim}\longrightarrow \gr H_{t,c}. \end{equation} In other words there is a left $\C G$-module isomorphism $H_{t,c} \stackrel{\sim} \longrightarrow \C [V]\otimes \C G$: this is called the {\it PBW isomorphism}. Note that $H_{0,0} = \C[V]\rtimes G$.

\subsection{Toy example} \label{TEdef}For $G = \mu_2$ acting on $V = \C^2$ there is a unique symplectic reflection, namely $s = -1 \in \mu_2$ and $\omega_s = \omega_{V*}$ since $(\id - s)(V) = V$. Thus the symplectic reflection algebra $H_{t,c}(\mu_2)$ depends on two parameters $t,c\in \C$ and is the quotient of $\C\langle x, y \rangle \rtimes G$ by the relation \begin{equation} \label{commmreln} yx - xy = t - 2cs.\end{equation} This relation allows us to put all elements of $H_{t,c}(\mu_2)$ into normal form and we find a basis $\{ x^iy^j , x^i y^j s: i, j \geq 0\}$.

\subsection{Spherical subalgebra} The algebra $\C[V]\rtimes G$ is noncommutative since if $f\in \C[V]$ and $g\in G$ we have $g \cdot f = {}^gf \cdot g$. Since $G$ acts faithfully on $V$ we have $Z(\C[V]\rtimes G) = \C[V]^G \subseteq \C[V]\rtimes G$. 

Now $\C[V]\rtimes G$ contains another subalgebra which is isomorphic to $\C[V]\rtimes G$. Let $e = |G|^{-1} \sum_{g\in G} g\in \C G$ be the trivial idempotent. The algebra $e(\C[V]\rtimes G)e$ is a subalgebra of $\C [V]\rtimes G$ (with identity element $e$) and there is an isomorphism $\C[V]^G \stackrel{\sim}\longrightarrow e(\C[V]\rtimes G)e$ given by $f \mapsto fe$ for any $f\in \C[V]^G$. Since $e\in \C G \subseteq H_{t,c}$ by the PBW theorem, this inspires the following definition.

\begin{definition}
The subalgebra $eH_{t,c}e$ of $H_{t,c}$ is called the {\it spherical subalgebra} of $H_{t,c}$. For the rest of the survey we will denote it by $U_{t,c}$.
\end{definition}
Since $U_{t,c}$ is a subspace of $H_{t,c}$ it inherits a filtration which is defined by $F^i(U_{t,c}) = U_{t,c}\cap F^i(H_{t,c})$. It's straightforward to see that \eqref{PBW} then implies 
\begin{equation} \label{PBWspher} \C[V]^G \cong e(\C[V]\rtimes G) e \stackrel{\sim}\longrightarrow \gr eH_{t,c}e. \end{equation} So the spherical subalgebras provide a good, i.e. flat, family of deformations of the coordinate ring of the symplectic singularity $V/G$, as required. 

It is the insight of Etingof and Ginzburg, which actually goes back at least to Crawley-Boevey and Holland, to study deformations of $\C[V]\rtimes G$ instead of $\C[V]^G$. Indeed we may hope to study the symplectic reflection algebras in order to understand better the $G$-equivariant geometry of $V$ since the category of $\C[V]\rtimes G$-modules is equivalent to the category of $G$-equivariant coherent sheaves on $V$.

\subsection{Toy example}\label{firsttoyex} In the example of \ref{TEdef} $e = \frac{1}{2}(1+s)$ and $U_{t,c}(\mu_2)$ is generated as a $\C$-algebra by ${\bf h} =  -\frac{1}{2}e(xy+yx)e, {\bf e}\defn \frac{1}{2}ex^2 e$ and ${\bf f}\defn \frac{1}{2}ey^2e$. There are relations \begin{align*} [{\bf e},{\bf f}] = t {\bf h},\, [{\bf h},{\bf e}]  = - 2t {\bf e}, \, [{\bf h},{\bf f}] =  2t{\bf f} \text{ and } {\bf e}{\bf f} = (2c- {\bf h}/2)({t}/{2}-c- {\bf h}/2). \end{align*} So if $t=0$ then $U_{t,c}(\mu_2)$ is commutative, whilst if $t = 1$ $U_{t,c}(\mu_2)$ is a central quotient of the enveloping algebra of $\mathfrak{sl}(2,\C)$.

\subsection{Symplectic reflection groups} We end this section with a brief reality check. The explicit form of the deformation $\kappa$ appearing in Theorem \ref{PBWthm} shows that the symplectic reflection algebras really only rely on the subgroup $H\leq G$ generated by the set $\mathcal{S}$ of symplectic reflections. Before symplectic reflection algebras were introduced, Verbitsky had seen this geometrically.
\begin{theorem}[{\cite[Theorem 3.2]{verbitsky}}]
Let $(G,V,\omega_V)$ be an indecomposable symplectic triple. if $V/G$ admits a symplectic resolution then $G$ is a symplectic reflection group.
\end{theorem}
This theorem should be considered as a {\it partial} analogue of the Theorem \ref{STC}: it gives necessary condition for the existence of a symplectic resolution, but, as we will see, this is not a necessary condition.

\subsection{$\!\!\!\!\!$} \label{groups}  We will assume for the rest of these notes that $$\text{$G$ is generated by symplectic reflections.}$$ We call such groups {\it symplectic reflection groups}.

There are two straightforward examples of symplectic reflection groups.
\begin{enumerate}
\item{\bf Complex reflections.} Let $G\leq GL(\h)$ be a complex reflection group. Set $V = \h \oplus
\h^*$ with its canonical symplectic form and with $G$ acting diagonally. Then $G\leq GL(V)$ is a symplectic reflection group. (We have doubled-up everything here, so since $G$ was generated by elements fixing hyperplanes in the action on $\h$ the same elements become symplectic reflections in the action on $\h\oplus \h^*$.) 
\item{\bf Wreath products.} Let $\Gamma \leq SL(2, \C) = Sp(2, \C)$ be finite: such groups are called kleinian subgroups and they preserve the canonical symplectic structure on $\C^2$. Set $$V = \underbrace{\C^2\oplus \C^2\oplus \cdots \oplus \C^2}_{n \text{ summands}}$$ with the symplectic form induced from that on $\C^2$ and let $G = \Gamma^n \rtimes \mathfrak{S}_n$ act in the obvious way on $V$. In this action $G$ is generated by symplectic reflections. \end{enumerate}

There is a little overlap in these two families. If we take $\Gamma$ to be the cyclic subgroup of $SL(2,\C)$ generated by $\diag (\exp(2\pi\sqrt{-1}/\ell), \exp(-2\pi \sqrt{-1}/\ell))$ where $\ell$ is some positive integer, then the action of $G$ on $\C^2$ restricts to the subspace $\C\times \{0\} \subset \C^2$ and so we can restrict the action of $G= \Gamma^n \rtimes \mathfrak{S}_n$ to the corresponding lagrangian subspace $\h \defn (\C\times \{0\})^n$. The action of $G$ on $\h$ is generated by complex reflections: in this guise the group $G$ is a complex reflection group, called $G(\ell, 1 , n)$, and we find the overlap between (1) and (2) above. 

The good news is that indecomposable symplectic reflection groups were classified by
Huffmann--Wales, \cite{HuWa}, see also Cohen, \cite{cohen}, and Guralnick--Saxl, \cite{guralsaxl}. Roughly speaking, the classification states that the above two classes are the only examples.

\section{Dimension 2}
In the special case that the symplectic triple $(G,V,\omega_V)$ is two dimensional, the family of symplectic reflection algebras was discovered by Crawley-Boevey and Holland, \cite{CBH}. They are basically deformed preprojective algebras for affine Dynkin quivers and the spherical subalgebras were called deformations of kleinian singularities in \cite{CBH}. In order to illustrate a number of the properties of symplectic reflection algebras we will recall some of the results of \cite{CBH}.

\subsection{Kleinian subgroups and singularities} \label{intersection} Let $G$ be a non-trivial finite subgroup of $SL(2,\C)$. The symplectic singularity $\C^2 /G$ is called a {\it kleinian} or {\it du Val singularity}: it is an extremely rich and well-studied meeting place for geometry, algebra and combinatorics. Each kleinian singularity has a unique singular point -- the zero orbit ${0}$ -- and has a unique symplectic resolution which can be constructed by a sequence of blow-ups $\pi: X = \widetilde{\C^2/G}\longrightarrow \C^2/G$. Each irreducible component of $\pi^{-1}({0})$ is a projective line. We form a graph whose vertices are labeled by the irreducible components of $\pi^{-1}({0})$ and which has an edge between two vertices if and only if the two components have non-trivial intersection. This graph turns out to determine the isomorphism type of the group $G$ and to be a Dynkin diagram of type $A$, $D$ or $E$. 

\subsection{$\!\!\!\!\!$} A remarkable observation of McKay, \cite{mckay}, gives a construction of the Dynkin diagram directly from the representation theory of $G$. We form a graph whose vertices are labeled by the isomorphism classes of irreducible representations $\{ S_i \}$ of $G$ and which has an edge between two vertices $i$ and $j$ if and only if $S_i$ is a summand of $\C^2 \otimes S_j$. The resulting graph determines the isomorphism type of the group and is an affine Dynkin diagram of type $A$, $D$ or $E$; removing the node corresponding to the trivial representation produces the intersection graph of \ref{intersection}. 

\begin{align*}
\widetilde{\mathrm{A}}_n && &
\xygraph{
[]
!{<0pt,0pt>;<20pt,0pt>:}
*\cir<2pt>{}
!{\save -<0pt,6pt>*\txt{$_1$}  \restore}
- [r]
*\cir<2pt>{}
!{\save -<0pt,6pt>*\txt{$_1$}  \restore}
- [r]
*{ \;  \dots \; }
- [r]
*\cir<2pt>{}
!{\save -<0pt,6pt>*\txt{$_1$}  \restore}
- [r]
*\cir<2pt>{}
!{\save -<0pt,6pt>*\txt{$_1$}  \restore}
- [llu]
\bullet{}
!{\save -<0pt,6pt>*\txt{$_1$}  \restore}
- [lld]
*\cir<2pt>{}
}\\\nonumber
\widetilde{\mathrm{D}}_n && &
\xygraph{
[]
!{<0pt,0pt>;<20pt,0pt>:}
[u]
*\cir<2pt>{}
!{\save -<0pt,6pt>*\txt{$_1$}  \restore}
- [rd]
*\cir<2pt>{}
!{\save -<0pt,6pt>*\txt{$_2$}  \restore}
- [ld]
*\cir<2pt>{}
!{\save -<0pt,6pt>*\txt{$_1$}  \restore}
[ru]
*\cir<2pt>{}
- [r]
*\cir<2pt>{}
!{\save -<0pt,6pt>*\txt{$_2$}  \restore}
- [r]
*{ \;  \dots \; }
- [r]
*\cir<2pt>{}
!{\save -<0pt,6pt>*\txt{$_2$}  \restore}
- [r]
*\cir<2pt>{}
!{\save -<0pt,6pt>*\txt{$_2$}  \restore}
- [ru]
*\cir<2pt>{}
!{\save -<0pt,6pt>*\txt{$_1$}  \restore}
[ld]
*\cir<2pt>{}
- [rd]
\bullet{}
!{\save -<0pt,6pt>*\txt{$_1$}  \restore}
}\\
\widetilde{\mathrm{E}}_6 && &
\xygraph{
[]
!{<0pt,0pt>;<20pt,0pt>:}
*\cir<2pt>{}
!{\save -<0pt,-6pt>*\txt{$_1$}  \restore}
- [r]
*\cir<2pt>{}
!{\save -<0pt,-6pt>*\txt{$_2$}  \restore}
- [r]
*\cir<2pt>{}
!{\save -<0pt,-6pt>*\txt{$_3$}  \restore}
- [d]
*\cir<2pt>{}
!{\save -<5pt,0pt>*\txt{$_2$}  \restore}
- [d]
*\cir<2pt>{}
!{\save -<5pt,0pt>*\txt{$_1$}  \restore}
[uu]
*\cir<2pt>{}
- [r]
*\cir<2pt>{}
!{\save -<0pt,-6pt>*\txt{$_2$}  \restore}
- [r]
\bullet{}
!{\save -<0pt,-6pt>*\txt{$_1$}  \restore}
}\\
\widetilde{\mathrm{E}}_7 && &
\xygraph{
[]
!{<0pt,0pt>;<20pt,0pt>:}
*\cir<2pt>{}
!{\save -<0pt,-6pt>*\txt{$_1$}  \restore}
- [r] *\cir<2pt>{}
!{\save -<0pt,-6pt>*\txt{$_2$}  \restore}
- [r]  *\cir<2pt>{}
!{\save -<0pt,-6pt>*\txt{$_3$}  \restore}
- [r]  *\cir<2pt>{}
!{\save -<0pt,-6pt>*\txt{$_4$}  \restore}
- [d] *\cir<2pt>{}
!{\save -<5pt,0pt>*\txt{$_2$}  \restore}
[u] *\cir<2pt>{}
- [r] *\cir<2pt>{}
!{\save -<0pt,-6pt>*\txt{$_3$}  \restore}
- [r] *\cir<2pt>{}
!{\save -<0pt,-6pt>*\txt{$_2$}  \restore}
- [r] \bullet{}
!{\save -<0pt,-6pt>*\txt{$_1$}  \restore}
}\\
\widetilde{\mathrm{E}}_8 && &
\xygraph{
[]
!{<0pt,0pt>;<20pt,0pt>:}
*\cir<2pt>{}
!{\save -<0pt,-6pt>*\txt{$_2$}  \restore}
- [r] *\cir<2pt>{}
!{\save -<0pt,-6pt>*\txt{$_4$}  \restore}
- [r]  *\cir<2pt>{}
!{\save -<0pt,-6pt>*\txt{$_6$}  \restore}
- [d] *\cir<2pt>{}
!{\save -<5pt,0pt>*\txt{$_3$}  \restore}
[u] *\cir<2pt>{}
- [r] *\cir<2pt>{}
!{\save -<0pt,-6pt>*\txt{$_5$}  \restore}
- [r] *\cir<2pt>{}
!{\save -<0pt,-6pt>*\txt{$_4$}  \restore}
- [r] *\cir<2pt>{}
!{\save -<0pt,-6pt>*\txt{$_3$}  \restore}
- [r] *\cir<2pt>{}
!{\save -<0pt,-6pt>*\txt{$_2$}  \restore}
- [r] \bullet{}
!{\save -<0pt,-6pt>*\txt{$_1$}  \restore}
}
\end{align*}
\begin{figure}[ht]
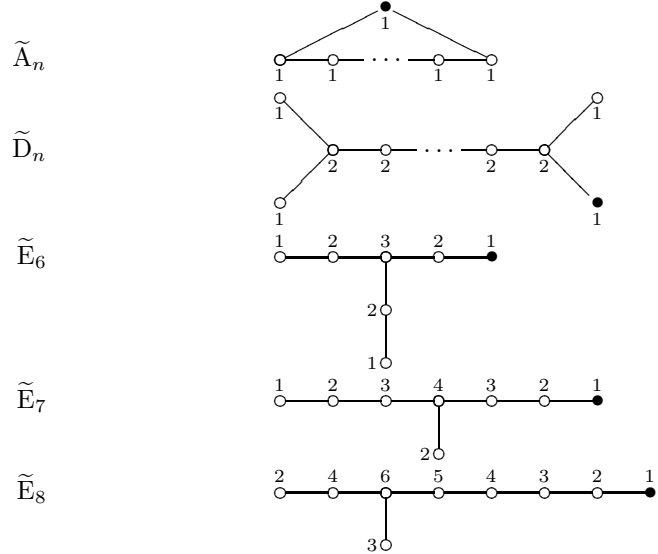

\label{affdyn}
\caption{Affine Dynkin diagrams with extending vertex emboldened}
\end{figure}

 \subsection{$\!\!\!\!\!$} This suggests that we should be able to associate naturally to each non-trivial irreducible representation of $G$ an irreducible component of $\pi^{-1}(0)$: the McKay correspondence. This was achieved by Gonzalez-Springberg--Verdier who gave a $K$-theoretic description; it was ``upgraded" to a derived equivalence by Kapranov--Vasserot $$D^b(\Lmod{\C[x,y]\rtimes G}) \stackrel{\sim} \longrightarrow D^b(\Coh (\widetilde{\C^2/G}))$$  where $D^b$ denotes the bounded derived category. Any representation of $G$ extends to a representation of $\C[x,y]\rtimes G$ by letting $x$ and $y$ act by zero. Under the equivalence the non-trivial irreducible $G$-representations are sent to bundles on $\pi^{-1}(0)$ which are non-trivial on one irreducible component only: this constructs the correspondence between irreducible representations and irreducible components.
 
\subsection{$\!\!\!\!\!$} Inspired by the McKay correspondence and earlier work of Hodges, \cite{CBH} makes two observations: the smash product $\C[x,y]\rtimes G$ is Morita equivalent to the preprojective algebra $\Pi^0(Q)$ of the affine Dynkin diagram $Q$ corresponding to $G$; the preprojective algebra is easy to deform. 

Let us begin by recalling the definition of the preprojective algebras $\Pi^0(Q)$ for any quiver ${\mathrm Q} = (\mathrm{Q}^0, \mathrm{Q}^1)$. Let $\overline{Q}$ be the double quiver of $Q$, obtained by inserting an arrow $a^*$ in the opposite direction to every arrow $a\in Q^1$ in the original quiver. Then the preprojective algebra is the following quotient of the path algebra of $\overline{Q}$:
$$ \Pi^0(Q) \defn \frac{ \C\overline{Q}}{\langle \sum_{a\in Q^1} [a,a^*] \rangle}.$$ Given a vector $\blambda = (\lambda_i)_{i\in Q^0} \in \C^{Q^0}$ we define the deformed preprojective algebra to be $$\Pi^{\blambda}(Q) \defn \frac{\C\overline{Q}}{\langle \sum_{a\in Q^1} [a,a^*] - \sum_{i\in Q^0} \lambda_i e_i \rangle} .$$

\subsection{$\!\!\!\!\!$} Since we are in the two dimensional case, every non-trivial element of $G$ is a symplectic reflection, i.e. $\mathcal{S} = G \setminus \{ 1\}$. Given $t\in \C$ and $c \in \C[\mathcal{S}]^{\ad G}$, define $\blambda(t, c) \in \C^{Q^0}$  by \begin{equation} \label{prpparam} \blambda(t, c)_i \defn \Tr_{S_i} \left(t1 - 2 \sum_{s\in \mathcal{S}} c(s)s\right).\end{equation} 
This produces a linear isomorphism between $\C\times \C[\mathcal{S}]^{\ad G}$ and $\C^{Q^0}$ such that $t = \blambda(t,c) \cdot \delta$ where $\delta\in \mathbb{N}^{Q^0}$ is the vector with $\delta_i = \dim S_i$.

\begin{theorem}[{\cite[Corollary 3.5]{CBH}}] \label{CBHthm} Let $Q$ be a quiver whose underlying graph is an affine Dynkin diagram, and let $G$ be the corresponding finite subgroup of $SL(2,\C)$ under the McKay correspondence. Then there is a Morita equivalence between $\Pi^{\blambda (t,c)}(Q)$ and $H_{t,c}(G)$. 
\end{theorem}
The result proved in \cite{CBH} is actually more precise and follows from McKay's construction of the affine Dynkin graphs via the representation theory of $G$. The theorem states: $$\Pi^{\blambda(t,c)} (Q) \cong f \frac{\C \langle x,y \rangle\rtimes G}{\langle xy -yx - \blambda(t,c) \rangle}f$$ where $\blambda$ is considered as the element $\sum_{i} \lambda_i \id_{S_i} \in \C G$ and $f = \sum_i f_i$ where $f_i$ is an idempotent in $\C G$ such that $\C G f_i \cong S_i$. Since $f\C Gf $ is Morita equivalent to $\C G$, this provides the desired Morita equivalence. It is then a simple calculation to check the correspondence between parameters.

\subsection{$\!\!\!\!\!$} We let $e_0\in \C\overline{Q}$ be the idempotent corresponding to the extending vertex -- a filled in vertex in Figure 1 --, and following \cite{CBH} set $\mathcal{O}^{\blambda} \defn e_0 \Pi^{\blambda} (Q) e_0$. As a simple corollary of the above theorem we find
\begin{corollary}
There is an algebra isomorphism $\mathcal{O}^{\blambda(t,c)} \cong U_{t,c}(G)$. \end{corollary}
Indeed since $\C G e_0 f \cong S_0$ we must have that $e_0 f = e$ and so $$e_0 \Pi^{\blambda(t,c)}(Q) e_0 \cong e_0 f \frac{\C \langle x,y \rangle\rtimes G}{\langle xy -yx - \blambda(t,c) \rangle}f e_0 = e  \frac{\C \langle x,y \rangle\rtimes G}{\langle xy -yx - \blambda(t,c) \rangle}e = U_{t,c}(G).$$ 

\subsection{Crawley-Boevey and Holland's results} Many structural results about $\Pi^{\bf \lambda}(Q)$, and hence about $H_{t,c}(G)$, were proved in \cite{CBH} and a lot of these now have analogues for all symplectic reflection algebras. We will recall them in this special case in the following portmanteau theorem. 
\begin{theorem} \label{defppthm} Let $G$ be a finite subgroup of $SL(2,\C)$  and $R$ the affine root system associated to the affine Dynkin diagram corresponding to $G$.
\begin{enumerate}
\item The centre of $H_{t,c}$ is non-trivial if and only if $t = 0$, in which case $Z(H_{t,c}) \cong U_{t,c}$. 
\item If $t = 0$ then $\Spec U_{t,c}$ is isomorphic to the variety of semisimple $H_{t,c}$-representations which are isomorphic to the regular representation when considered as $G$-representations.
\item If $t \neq 0$ then the isomorphism classes of finite-dimensional simple $H_{t,c}$-modules are in one-to-one correspondence with the set $\Sigma_{\blambda(t,c)}$ of minimal elements of $\{ \alpha \in R: \blambda(t,c) \cdot \alpha = 0 \}$ 
\item There is a Morita equivalence between $U_{t,c}$ and $H_{t,c}$ if and only if $\blambda(t,c) \cdot \alpha \neq 0$ for every Dynkin root $\alpha$.
\item $\mathcal{O}^{\blambda}$ is simple if and only if $\blambda \cdot \alpha\neq 0$ for all non-Dynkin roots $\alpha$.
\item $$\text{\sf gl.dim.} (U_{t,c}) = \begin{cases} 1 & \blambda(t,c) \cdot \alpha \neq 0 \text{ for all roots $\alpha$}, \\  \infty \quad & \blambda(t,c) \cdot \alpha = 0\text{ for some Dynkin root $\alpha$}, \\ 2 & \text{otherwise}. \end{cases}$$
\item The quotient division ring of $U_{t,c}$ is isomorphic to the quotient division ring of $\C \langle x, y : xy-yx = t \rangle.$
\end{enumerate}
\end{theorem}
\subsection{$\!\!\!\!\!$} We will not discuss the proofs of these results as we will present and outline generalisations of many of them later. However, note that Part (1) provides a flat family of deformations $Y_{c} \defn \Spec U_{0,c}$ of the kleinian singularity $Y_0 = \C^2/G$ over the affine space $\C[\mathcal{S}]^{\ad G}$. This is the {\it universal deformation} studied by Slodowy and others, \cite{slodowy}. Part (6) shows that the generic deformation is smooth. When $t=0$, the algebra $H_{t,c}$ is a finite module over its centre, and it follows quickly that every irreducible $H_{t,c}$-representation is finite dimensional. In contrast, Part (3) shows that if $t \neq 0$ then there are only finitely many irreducible finite dimensional $H_{t,c}$-representations.
\section{Basic properties}
There are a few results which are valid for all symplectic reflection algebras. We present them here. Throughout $(G,V,\omega_V)$ will be an indecomposable symplectic triple.

\subsection{Ring theoretic properties} The PBW isomorphisms \eqref{PBW} and \eqref{PBWspher} allow us to deduce some interesting properties for $H_{t,c}$ and $U_{t,c}$ from the same properties for their associated graded algebras. By \cite[Lemma 6.11, Corollary 6.18]{MCR} these include
\begin{proposition}
\label{easyfilter}
\begin{enumerate}
\item
The symplectic reflection algebra $H_{t,c}$ is noetherian, prime and of finite global dimension.
\item
The spherical subalgebra $U_{t,c}$ is a noetherian domain with finite injective dimension.
\end{enumerate}
\end{proposition}
In fact $\text{\sf gl.dim.} H_{t,c} = \dim V - \min  \{ \text{\sf GKdim} (I): I \text{ irreducible $H_{t,c}$-representation}\}$ where {\sf GKdim} stands for Gelfand--Kirillov dimension, cf. Theorem \ref{defppthm}(6).

Note we cannot deduce that $U_{t,c}$ has finite global dimension since $\gr U_{t,c} \cong \C[V]^G$ is the coordinate ring of a singular variety. It is a non-trivial question to decide for which $(t,c)$ $U_{t,c}$ has finite global dimension. 

\subsection{$\!\!\!\!\!$} Since $H_{t,c} e$ is an $(H_{t,c},U_{t,c})$-bimodule, we have a homomorphism \begin{equation} \label{endo}\mu: H_{t,c} \longrightarrow \End_{U_{t,c}} (H_{t,c}e)\end{equation} whose associated graded mapping can be identified with the mapping $\gr \mu: \C[V]\rtimes G \longrightarrow \End_{\C [V]^G}(\C [V])$ given by the natural $\C[V]\rtimes G$-module structure on $\C[V]$. Now $\C[V]$ is a faithful $\C[V]\rtimes G$-module since $G$ acts faithfully on $V$, so $\gr \mu$ is injective. Galois theory ensures that $\gr\mu$ becomes an isomorphism on passing to the quotient field $\C(V)$, so any element of $\End_{\C[V]^G}(\C[V])$ can be written as $\sum a_g \cdot g$ where $a_g \in \C(V)$ is regular on $\vr$, since $G$ acts freely there. But the complement to $\vr$ has codimension at least 2, see  \ref{sympres}, and so by Hartog's theorem each function $a_g$ is in fact regular on $V$. Thus $\gr \mu$ is surjective and hence it is an isomorphism. Standard associated graded techniques then show that $\mu$ is also an isomorphism. Of course, it is immediate that $\End_{H_{t,c}} (H_{t,c}e) \cong U_{t,c}^{op}$.

\subsection{Centres} The following theorem presents a dichotomy in the behaviour of symplectic reflection algebras which percolates through all of their representation theory.
The difficult Part (1) was proved in \cite[Theorem 1.6]{EG}, the easy Part (2) in \cite[Proposition 7.2]{BG}, cf. Theorem \ref{defppthm}(1).

\begin{theorem}
\label{centrethm}
The mapping $Z(H_{t,c}) \longrightarrow U_{t,c}$ which sends $z$ to $ez$ induces an isomorphism onto the centre of $U_{t,c}$. Moreover,
\begin{enumerate}
\item if $t = 0$ then $Z(U_{t,c}) = U_{t,c}$, i.e. $U_{t,c}$ is commutative; \item  if $t\neq 0$ then $Z(U_{t,c}) = \C$.\end{enumerate}
\end{theorem}
The proof begins by observing that the mapping in the statement of the theorem is indeed an algebra homomorphism with image in $Z(U_{t,c})$. Now take $y\in Z(U_{t,c})$ and note that right multiplication by $y$ on $H_{t,c}e$ is a $U_{t,c}$-endomorphism since $y$ is central.  By the preceding paragraph this gives rise to an element of $H_{t,c}$, say $\hat{y}$, and this element must be central since right multiplication by $y$ on $H_{t,c}e$ commutes with left multiplication by $H_{t,c}$. Now $ \hat{y} e = \mu(\hat{y})(e) = e y =y$ so we have constructed the inverse to the mapping.  

The second two statements rely on a study of the properties of the function $m(t,c): \C \times \C[\mathcal{S}]^{\ad G} \longrightarrow \mathbb{N}\cup \{\infty\}$ which measures the noncommutativity of $U_{t,c}$ and the corresponding Poisson bracket, as in  \ref{defpoisson}. 

\subsection{$\!\!\!\!\!$} For the rest of this paper, given $c\in \C[\mathcal{S}]^{\ad G}$ we will set $$Y_{c} \defn \Spec Z(H_{0,c}).$$ By Theorem \ref{centrethm} $Y_c \cong \Spec U_{0,c}.$
\subsection{$\!\!\!\!\!$} \label{sympsing} Fix $c\in \C[\mathcal{S}]^{\ad G}$ and define the $\C[\hbar]$-algebra $H_{\hbar ,c}$ exactly as for $H_{t,c}$ except that $t$ is replaced by the central indeterminate $\hbar$; applying the idempotent $e$ produces the $\C[\hbar]$-algebra $U_{\hbar,c}$. By \eqref{PBWspher} $U_{\hbar,c}$ is a deformation of $U_{0,c}$ and so, by  \ref{defpoisson} and Theorem \ref{centrethm}, $U_{0,c}$ is a Poisson algebra.
\begin{proposition} For any $c\in \C[\mathcal{S}]^{\ad G}$ the variety $Y_{c}$ is a symplectic singularity. The family $(Y_c)_c$ is a flat family of Poisson deformations of $Y_0 = V/G$ over $\C[\mathcal{S}]^{\ad G}$.
\end{proposition}
\begin{proof} We will give brief details since this result has not been written down before.
By \cite[Lemma 3.10]{gordsmi} $H_{0,c}$ is a noncommutative crepant resolution of $Z(H_{0,c})$ and so, by \cite[Theorem 4.3]{SvdB} $Y_c$ has rational singularities. Furthermore, \cite[Theorem 7.8]{BG}, the restriction of the Poisson form to $(Y_c)_{\sf sm}$ is non-degenerate. It then follows from \cite[Theorem 6]{Namikawa} that $Y_c$ has symplectic singularities.

The second statement follows immediately from \eqref{PBWspher} and the fact that the Poisson bracket on $Y_0$ agrees with the one on $V/G$ inherited from the $G$-invariant symplectic form $\omega_{V^*}$ on $V$ since it is induced from the quantisation $T(V^*) / (xy-yx - t\omega_{V^*}(x,y))$ of $\C[V]$.
\end{proof} 

\subsection{$\!\!\!\!\!$} Thus we have constructed a good family $(Y_c)_c$ of deformations of $V/G$ which are, moreover, equipped with two extra pieces of information: a coherent sheaf
corresponding to $H_{0,c}e$ and a quantisation, $U_{t,c}$. It is partly the interplay of these structures that will allow us to prove interesting theorems.

\subsection{$\!\!\!\!\!$} The algebra $H_{0,c}$ is a finite module over its centre $Z(H_{0,c})$ and so a simple lemma, which goes back at least to Kaplansky, shows that every irreducible $H_{0,c}$-module is a {\it finite dimensional} complex vector space, see for example \cite[1.3]{jantzen}. Thus, by Schur's lemma, every central element $z$ must act on an irreducible representation $I$ by a scalar, which we denote $\chi_I(z)$. So there is a {\it central character} mapping \begin{equation} \label{cenchar} \chi: \irr (H_{0,c}) \longrightarrow Y_c,\end{equation} which provides the direct comparison between representation theory and geometry.

\subsection{$\!\!\!\!\!$} At first sight there is not really a good geometry attached to $\irr (H_{0,c})$ and so it is not clear how to play the domain and codomain of $\chi$ against each other. However, the prototype theorem linking the two together is the following result which has a long history going back to Artin, Procesi, LeBruyn, Brown--Goodearl and then proved independently in this context by Etingof and Ginzburg, \cite[Theorem 1.7]{EG}, cf. Theorem \ref{defppthm}(2).

\begin{theorem}
\label{azumayathm}
\begin{enumerate}
\item $|G|$ is an upper bound on the dimension of irreducible $H_{0,c}(G)$-representations.
\item If $I\in \irr (H_{0,c}(G))$ with $\dim_\C I = |G|$ then $I\cong_{\C G} \C G$.
\item Let $\mathcal{A}_c \defn  \{ I\in \irr (H_{0,c}(G)) : \dim_{\C} I = |G| \}$, the so-called {\it Azumaya locus}. Then $\chi$ restricts to a bijection $\mathcal{A}_c \stackrel{\sim}\longrightarrow (Y_c)_{\sf sm}$ onto the smooth locus of $Y_c$.
\end{enumerate}
\end{theorem}

The heart of this theorem is Part (3). Quite generally, there are criteria given by \cite[Theorem 3.8]{BrownGoodearl} which give sufficient conditions for the smooth locus and the image of the Azumaya locus to agree for algebras finite over their centre. Most of the criteria are homological in nature and thus easy to check here because of Proposition \ref{easyfilter}. However, one of the criteria is that $\codim_{Y_c} \chi (\mathcal{A}_c) \geq 2$. As we will see in the next paragraph, this is a consequence of the symplectic structure.

\subsection{$\!\!\!\!\!$} \label{killcen}The central character \eqref{cenchar} allows us to partition the irreducible $H_{0,c}$-representations fibre-by-fibre. In other words it is enough to study the various $\chi^{-1} (y)$ as $y$ varies in $Y_c$. But an irreducible representation in $\chi^{-1}(y)$ is exactly the same as an irreducible representation of the algebra $H_{0,c}(y) \defn H_{0,c}/\mathfrak{m}_y H_{0,c}$ where $\mathfrak{m}_y$ is the maximal ideal of $Z(H_{0,c})$ corresponding to $y$. Thus we have turned the problem of describing the irreducible representations of the infinite dimensional algebra $H_{0,c}$ into the problem of describing the irreducible representations of the infinitely many finite dimensional algebras $H_{0,c}(y)$! 

This is progress. Since $Y_c$ has symplectic singularities, we can stratify $Y_c$ by some special closed sets $\emptyset = Y^r \subset Y^{r-1} \subset \cdots \subset Y^0 = Y_c$ where for any $0\leq i < r$ $Y^{i} \setminus Y^{i+1}$ is smooth and symplectic and moreover $Y^0 \setminus Y^1 = (Y_c)_{\sf sm}$, \cite[Proposition 3.1]{kaledinsurvey}. These strata are obtained as the spaces traced out by all the flows of the hamiltonian vector fields $\{ z, \cdot \}_c$ where $z\in Z(H_{0,c})$ and $\{ \cdot , \cdot \}_c$ is the Poisson bracket on $Z(H_{0,c})$. They are called the {\it symplectic leaves} of $Y_c$.

\begin{theorem}[{\cite[Theorem 4.2]{BG}}] Let $y_1,y_2 \in Y_c$ belong to the same symplectic leaf of $Y_c$. Then there is an algebra isomorphism $H_{0,c}(y_1) \cong H_{0,c}(y_2)$.
\end{theorem}

Since $y_1$ and $y_2$ belong to the same symplectic leaf they can be connected by a sequence of hamiltonian flows; the theorem is proved by lifting these to the algebra $H_{0,c}$. 

\subsection{$\!\!\!\!\!$}
\label{noflatcentre}
There are only finitely symplectic leaves -- the stratification above was finite -- so it follows that there are only finitely many different isomorphism classes of finite dimensional algebras $H_{0,c}(y)$ and hence \begin{center} {\em to understand all the irreducible representations of $H_{0,c}$ it is sufficient to classify \\ the irreducible representations of only a finite number of finite dimensional algebras.}
\end{center}
However, these finite dimensional algebras may be rather difficult to understand. For instance their dimension varies with $y\in Y_c$ since $H_{0,c}$ is {\it not} a flat $Z(H_{0,c})$-module. For any $c\in \C[\mathcal{S}]^{\ad G}$ and a generic choice of $y\in Y_{c}$ the algebra $H_{0,c}(y)$ has dimension $|G|^2$ (this will follow from Theorem \ref{azumayathm}); if we set $c = 0$ so that $Y_c = V/G$ and take $y$ to be the zero orbit then the corresponding finite dimensional algebra is $\co{V}{G}\rtimes G$ whose dimension, in general, is unknown -- for $G = S_n$ acting on $V = \C^n\oplus \C^n$ the dimension is $n!(n+1)^{n-1}$ thanks to a celebrated recent theorem of Haiman, \cite[Proposition 3.6]{haimandiag}.

\subsection{$\!\!\!\!\!$}
We have seen that the dimension of the algebras $H_{0,c}(y)$ vary as we vary $y\in Y_c$, but that they are constant on symplectic leaves. Since $H_{0,c} \cong \End_{U_{0,c}}(H_{0,c}e)$ by \eqref{endo}, this translates to the fact that the dimension of $H_{0,c}e\otimes_{U_{0,c}} \C_y$ is constant on symplectic leaves as we run through the irreducible $U_{0,c}$-representations $\C_y$ which,  since $U_{0,c}\cong Z(H_{0,c})$, are labelled by $y\in Y_c$. Since the smooth locus $(Y_c)_{\sf sm}$ is itself a symplectic leaf we see that there is a dense open set of $Y_c$ over which the coherent sheaf $H_{0,c}e$ has constant rank and is thus a vector bundle, of rank $|G|$ by Theorem \ref{azumayathm}(3). Hence we deduce the following theorem.

\begin{theorem}[{\cite[Theorem 1.7]{EG}}] 
\label{smoothpartcentre}
There is an algebra isomorphism $$H_{0,c}|_{(Y_c)_{\sf sm}} \stackrel{\sim} \longrightarrow \Mat_{|G|} \left( \C[Y_c]|_{(Y_c)_{\sf sm}} \right).$$
\end{theorem}

\begin{remark}This isomorphism does not extend beyond the smooth locus of $Y_c$. Indeed, if the dimension of $H_{0,c}e$ remained constant on an open set $U$ larger than $(Y_c)_{\sf sm}$ then the isomorphism of the theorem would extend and we would deduce that the global dimension of the restrictions of $H_{0,c}$ and $\C[Y_c]$ to $U$ are the same. But $H_{0,c}$ has finite global dimension by Proposition \ref{easyfilter}, whereas $U$ has singularities and so the $\C[Y_c]|_U$ has infinite global dimension.
\end{remark}

\subsection{$\!\!\!\!\!$} \label{GKconj}There is one other general result for symplectic reflection algebras. In fact, it is not quite known to be true generally at the moment, but rather only for $(G,V,\omega_V)$ belonging to Families (1) and (2) in \ref{groups}, cf. Theorem \ref{defppthm}(7).
\begin{theorem}[{\cite[(1.18)]{EG}, \cite[Theorem 1.2]{gordGK}}]
The quotient division ring of $U_{t,c}$ is isomorphic to the quotient division ring of $T_{\C}(V^*)/ \langle xy-yx = t\omega_{V^*}(x,y) : x,y\in V^* \rangle.$
\end{theorem}
In particular the birational equivalence class of $H_{t,c}(G)$ is independent of $c$.

\section{Specific case: rational Cherednik algebras I ($t\neq 0$)}
In \ref{groups} we saw that there were essentially two families of groups generated by symplectic reflections. The corresponding symplectic reflection algebras tend to be treated differently depending on which family we are dealing with. When the triple $(G,V,\omega_V)$ belongs to Family (1), so that $G< GL(\h)$ is a complex reflection group and $V = \h\oplus \h^*$, we call $H_{t,c}(G)$ a rational Cherednik algebra. The reason for this name is that if $G$ is a Weyl group, then $H_{t,c}(G)$ is a degeneration of the {\it double affine Hecke algebra} which was discovered by Cherednik, \cite{cherednikbook}. 

The representation theory of rational Cherednik algebra has a lot in common with the representation theory of simple complex Lie algebras and we will see much of this in the next two chapters. In this chapter we introduce category $\mathcal{O}$ and the ${\sf KZ}$-functor which links $H_{t,c}(G)$ to Hecke algebras. We also discuss the finite dimensional representations of $H_{t,c}(G)$ and explain their application to invariant theory.

\label{t1RCA}
\subsection{$\!\!\!\!\!$}Throughout this section $G<\GL(\h)$ is a complex reflection group of rank $n$ acting on $V = \h\oplus \h^*$. Such groups were classified by Shephard-Todd: there is one infinite family $G(m,d,n)$ which includes the Weyl groups of type $A,B$ and $D$ and the dihedral groups, and there are 34 exceptional groups, see for instance \cite[Appendix 2]{BMR}.   

We will also assume that $t\neq 0$. Since the mapping  $x\mapsto \lambda x, y\mapsto \lambda y, w\mapsto w$ induces an isomorphism $H_{t,c}\cong H_{\lambda^2 t, \lambda^2 c}$ we can and will assume that $t=1$.   

\subsection{Triangular decomposition} It is a trivial but important observation to rewrite the PBW isomorphism for $H_{t,c}$ in this case, using the polarisation $V = \h\oplus \h^*$. This gives $$\C[\h]\otimes \C G\otimes \C[\h^*] \stackrel{\sim}\longrightarrow \gr H_{t,c}.$$ In this decomposition each tensorand is actually a subalgebra of $H_{t,c}$ and thus $H_{t,c}$ is an algebra with triangular decomposition. The prototype for an algebra with a triangular decomposition is the universal enveloping algebra of a simple complex Lie algebra $\mathfrak{g}$ in which case the decomposition is induced from the standard direct sum decomposition $\mathfrak{g} = \mathfrak{n}_+ \oplus \h \oplus \mathfrak{n}_-$. Thus the subalgebras $\C[\h]\rtimes G$ and $\C[\h^*]\rtimes G$ of $H_{t,c}$ play the roles of the enveloping algebras of the positive and negative Borel subalgebras. It turns out that many of the properties of the representation theory of $H_{t,c}$ are analogous to similar properties for $U(\g)$, or more precisely for minimal primitive quotient rings of $U(\g)$.

\subsection{Category $\mathcal{O}$} The first analogue we will see of this phenomenon is the construction of $\widehat{\mathcal{O}}_c$ for $H_{1,c}$. This interesting category, which contains all the finite dimensional representations of $H_{1,c}$ for instance, was introduced and studied in \cite{dunklopdam}, \cite{BEGqi}, \cite{guay} and \cite{GGOR}. If we ever need to refer to the group $G$ in question, we will write $\widehat{\mathcal{O}}_c^G$ for this category, and if we ever need to make the reflection representation explicit too we will write $\widehat{\mathcal{O}}_c^G(\mathfrak{h})$ -- a complex reflection group can have more than one reflection representation, for instance possibly $\h^*$. 
\begin{definition}
\begin{enumerate}
\item $\widehat{\mathcal{O}}_c$ is the full subcategory of finitely generated $H_{1,c}$-modules on which $\h\subset \C[\h^*]$ acts locally finitely, i.e. if $M\in \widehat{\mathcal{O}}_c$ then for all $m\in M$ $\dim (\C[\h^*]\cdot m) < \infty$.
\item An object $M\in \widehat{\mathcal{O}}_c$ has {\it type} $\lambda \in \h^*/G = \Spec (\C[\h^*]^G)$ if  for any $P\in \C[\h^*]^G$ the action of $P-P(\lambda)$ is locally nilpotent, i.e. for all $m\in M$ $(P-P(\lambda))^N\cdot m = 0$ for large enough $N$. We set $\widehat{\mathcal{O}}_c(\lambda)$ to be the full subcategory of $\widehat{\mathcal{O}}_c$ consisting of objects of type $\lambda$. 
\item In the special case $\lambda = 0 \in \h^*/G$ we set $\mathcal{O}_c \defn \widehat{\mathcal{O}}_c(0)$. 
\end{enumerate}
\end{definition}
So by definition $\mathcal{O}_c$ is the full subcategory of finitely generated $H_{1,c}$-modules on which every  $G$-invariant polynomial without constant term acts locally nilpotently. Since the coinvariant ring $\C[\h^*]^{\coinv G} = \C[\h^*]/\langle \C[\h^*]^G_+\rangle$ is a positively graded finite dimensional algebra, the image of every element in $\C[\h^*]_+$ is nilpotent in $\co{\h^*}{G}$ and so we have an equivalent definition:
\begin{definition}
$\mathcal{O}_c$ is the full subcategory of finitely generated $H_{1,c}$-modules on which every $P\in\C[\h^*]_+$ acts locally nilpotently.
\end{definition}
\subsection{$\!\!\!\!\!$} It is possible to reduce the study of $\widehat{\mathcal{O}}_c$ to the study of $\mathcal{O}_c$ thanks to recent work of Bezrukavnikov and Etingof, \cite{BE}. Given $\lambda \in \h^*/W$ we let $G_{\lambda}$ be the $G$-stabiliser of any lift of $\lambda$ to $\h^*$: this is of course only well-defined up to conjugacy, but \cite[Theorem 1.5]{Steinberg} shows that $G_{\lambda}$ is again a complex reflection group under its action on $\h/\fix (G_{\lambda})$, generated by a subset of the reflections in $G$. The first part of the following theorem is quite straightforward, whereas the second part is deeper.

\begin{theorem}
\begin{enumerate}
\item There is a decomposition of categories $\widehat{\mathcal{O}}_c = \bigoplus_{\lambda \in \h^*/G} \widehat{\mathcal{O}}_c(\lambda)$.
\item There is an equivalence of categories $\widehat{\mathcal{O}}_c^G(\lambda) \stackrel{\sim} \rightarrow \widehat{\mathcal{O}}_{c'}^{G_{\lambda}}(0) = \mathcal{O}_{c'}^{G_{\lambda}}$ where $c'$ is the restriction of $c$ to the reflections in $G_{\lambda}$, \cite[Corollary 3.3]{BE}.
\end{enumerate}
 \end{theorem}
Thus, without loss of generality, we need only study the category $\mathcal{O}_c$ in order to understand all of $\widehat{\mathcal{O}}_c$.
\begin{remark} In studying a simple complex Lie algebra, $\mathfrak{g}$, one defines the Bernstein--Gelfand--Gelfand category $\mathcal{O}$ by looking at finitely generated $U(\mathfrak{g})$-modules on which $\mathfrak{n}_-$ acts locally nilpotently {\it and} $\h$ acts completely reducibly. In the rational Cherednik algebra situation we do not need to insist on this second condition since $\C G$ is a finite dimensional semisimple algebra so any object in $\mathcal{O}_c$ will automatically be completely reducible for the action of this subalgebra of $H_{1,c}$. So our definition of $\mathcal{O}_c$ is indeed analagous to the Lie theory setting.
\end{remark} 
\subsection{$\!\!\!\!\!$} $\mathcal{O}_c$ has a very nice homological structure: it is a highest weight category, in the sense of \cite{CPS}. In particular this means that $\mathcal{O}_c$ is an abelian category with enough projectives and a distinguished set of {\it standard objects}. The standard objects are particularly easy to define and are analogues of Verma modules for simple complex Lie algebras. 
\begin{definition}
Let $\lambda \in \irr(G)$ be an irreducible (complex) representation of $G$. We define the {\it standard module associated to $\lambda$} to be $$\Delta_c(\lambda) \defn H_{1,c}\otimes_{\C[\h^*]\rtimes G} \lambda ,$$ where $\C[\h^*]$ acts on $\lambda$ via $p \cdot v = p(0) v$ for any $p\in \C[\h^*]$ and $v\in \lambda$. 
\end{definition} 
By the PBW isomorphism \eqref{PBW} we see that as a $\C[\h]\rtimes G$-module $\Delta_c(\lambda)$ is isomorphic to $\C[\h]\otimes \lambda$ with $\C[\h]$ acting on the first tensorand and $G$ acting diagonally. 

In order to check that $\Delta_c(\lambda) \in \mathcal{O}_c$ we will introduce a grading on $\Delta_c(\lambda)$ which will prove to be extremely useful in what follows. Let $\{ y_i : 1\leq i \leq n\}$ be a basis for $\h$, and let $\{ x_i \}$ be a dual basis. Define \begin{equation} \label{hdefinition} {\bf h} = -\frac{1}{2}\sum_{i=1}^n (x_iy_i + y_i x_i) \in H_{1,c}. \end{equation}
It's a straightforward calculation to see that this element is independent of the choice of basis of $\h$ and thus is invariant under conjugation by $G$. Moreover, it is proved in \cite[Lemma 2.5]{BEGqi} that $$[{\bf h}, x] = -x  \text{ for all $x\in \h^*$ and $[{\bf h}, y ] = y$ for all $y\in \h$,}$$ cf.  \ref{firsttoyex}. Now if we let ${\bf h}$ act on $1\otimes \lambda \in \Delta_c(\lambda)$ we find that \begin{equation}\label{centralelement}{\bf h} \cdot (1\otimes \lambda)  =    -\frac{1}{2}\sum_{i=1}^n (x_iy_i + y_i x_i) \cdot 1\otimes \lambda  = - \frac{1}{2}\sum_{i=1}^n y_ix_i \cdot (1\otimes \lambda)  = \left(-\frac{n}{2} +  \sum_{s\in \mathcal{S}}  c(s) s \right) \cdot (1\otimes \lambda) \end{equation} where we used the defining commutation relation of Theorem \ref{PBWthm} for $H_{1,c}$ to get the last equality. Since $-n/2 + \sum_{s\in \mathcal{S}} c(s) s\in \C G$ is invariant under conjugation, it acts on the irreducible $G$-representation $\lambda$ by a scalar which we denote by $\kappa (c, \lambda)$. Hence if $p\in \C[\h]$ is a homogeneous element of degree $m$ then $p\otimes \lambda\in \Delta_c(\lambda)$ is an eigenvector for ${\bf h}$ with eigenvalue $\kappa(c,\lambda) - m$. Thus it makes sense to talk about the decomposition of $\Delta_c(\lambda)$ into ${\bf h}$-eigenspaces and to say that $\Delta_c(\lambda)$ is a highest weight module with highest weight space $\C \otimes \lambda$ of weight $\kappa(c,\lambda)$. We will write the weight space decomposition as \begin{equation} \label{wtspchar} [ \Delta_c(\lambda)]_{\bf h} = \dim (\lambda)\frac{u^{\kappa(c,\lambda)}}{(1-u^{-1})^n}.\end{equation} Since the action of $\h\subset \C[\h^*]$ on $\Delta_c(\lambda)$ increases degree, it follows that the action of $\h$ is locally nilpotent and so $\Delta_c(\lambda)\in \mathcal{O}_c$.

\subsection{$\!\!\!\!\!$} \label{orderonG} We define an ordering on $\irr(G)$ which depends on $c$. For $\lambda , \mu \in \irr(G)$ $$
\lambda <_c \mu \text{ if and only if } \kappa (c,\mu) - \kappa(c,\lambda) \in \mathbb{Z}_+.$$
So, since $c$ is a complex-valued function, for generic choices of $c$ (including $c  =0$) all $\lambda\in \irr(G)$ are incomparable.

\subsection{$\!\!\!\!\!$} It is now clear that every $\Delta_c(\lambda)$ has a unique irreducible quotient $L_c(\lambda)$. Indeed, any proper submodule of $\Delta_c(\lambda)$ must be contained in the ${\bf h}$-eigenspaces $\Delta_c(\lambda)_{< \kappa(c,\lambda)}$ and so there exists a unique maximal proper submodule, namely the sum of all proper submodules. Furthermore any irreducible object $I$ in $\mathcal{O}_c$ must have a $G$-stable subspace annihilated by $\C[\h^*]_+$; fixing an irreducible summand inside this subspace induces a non-zero $H_{1,c}$-homomorphism from $\Delta_c(\lambda)$ to $I$, which by irreducibility must be surjective. Thus $I\cong L_c(\lambda)$.   This confirms the second bullet point of the following theorem.
\begin{theorem}\label{HWthm} $(\mathcal{O}_c, <)$ is a highest weight category with
\begin{itemize}
\item standard objects $\Delta_c(\lambda)$ for $\lambda \in \irr(G)$,
\item irreducible objects $L_c(\lambda)$ for $\lambda \in \irr(G)$,
\item projective objects $P_c(\lambda)$ for $\lambda \in \irr(G)$,
\item ordering $L_c(\lambda) < L_c(\mu)$ if and only if $\lambda <_c \mu$.
\end{itemize}
\end{theorem}
The projective objects are constructed by induction from ``big enough" $\C[\h^*]\rtimes G$-representations. It is crucial for the proof of this theorem that the grading in $\mathcal{O}_c$ is internal, i.e. constructed from the action of the element ${\bf h}$: this ensures all homomorphisms are graded. 

\subsection{$\!\!\!\!\!$} \label{Ocomments}
There are a few easy but interesting observations which we can make immediately. 

\medskip
\noindent
(1) A version of {\it BGG reciprocity} holds for $\mathcal{O}_c$:
\begin{equation}\label{BGGrecip}[P_c(\lambda) : \Delta_c(\mu) ] = [\Delta_c(\mu): L_c(\lambda)] \text{ for all $\lambda, \mu\in \irr(G)$}.\end{equation} Here $[P_c(\lambda) : \Delta_c(\mu)]$ counts the number of copies of $\Delta_c(\mu)$ in some filtration of $P_c(\lambda)$ whose sections are all standard modules: it is part of the axiomatics of highest weight categories that such a filtration exists and that this number is independent of the choice of filtration. This result is not quite immediate from the definition of highest weight category, but follows from a specialisation argument, see \cite[Proposition 3.3]{GGOR}.

\medskip
\noindent
(2) Every finite dimensional representation of $H_{1,c}$ belongs to $\mathcal{O}_c$, \cite[Th\'eor\`eme 4.1]{De}. Indeed if $I$ is a finite dimensional $H_{1,c}$-representation then $I$ must have a bounded ${\bf h}$-spectrum; thus, for $N$ large enough, all elements of $\C[\h^*]_{\geq N}$ must annihilate $I$, as required.

\medskip
\noindent
(3) For a generic choice of $c$ all $L_c(\lambda)$'s are incomparable by \ref{orderonG}. Thus by \eqref{BGGrecip} \begin{equation} \label{genericss} \mathcal{O}_c \text{ is semisimple for almost all $c\in \C[\mathcal{S}]^{\ad G}$}.\end{equation}

\medskip
\noindent
(4)  \label{involsec} There are two involutions on $H_{1,c}$.  
Let $\dagger$ be the involution on $\C[\mathcal{S}]^{\ad G}$ that takes $c$ to the function $c^{\dagger} : s\mapsto c(s^{-1})$: abusing notation, there is an anti--isomorphism $$\dagger: H_{1,c} \longrightarrow H_{1,c^{\dagger}}\qquad x\mapsto x, y\mapsto -y\text{ and }w\mapsto w^{-1};$$
similarly there is another anti--isomorphism $$\phi: H_{1,c}(\h\oplus \h^*) \longrightarrow H_{1,c^{\dagger}}(\h^*\oplus \h) \qquad  x\mapsto x, y\mapsto y \text{ and } w\mapsto w^{-1}.$$ These induce two dualities on $\mathcal{O}_c$.

(a) {\it Na\"ive duality}, \cite[Proposition 4.7]{GGOR}, $
d_{\h ,c} : \mathcal{O}_c^G(\h\oplus \h^*) \longrightarrow \mathcal{O}_{c^{\dagger}}^G(\h^*\oplus \h)$
is obtained by sending an object $M$ to the largest submodule $d_{\h,c}(M)$ of $\Hom_{\C}(M, \C)$ which is locally $\C[\h^*]$-nilpotent. This 
is a right $H_{1,c}(\h\oplus \h^*)$-module and becomes a left $H_{1,c^{\dagger}}(\h^*\oplus \h)$-module via $\phi$. It sends $L_c(\lambda)$ to 
$L_{c^{\dagger}}(\lambda^*)$, the irreducible object indexed by the dual representation $\lambda^*$ of $\lambda\in \irr(G)$; projective objects are sent to injectives and standard objects to costandard objects.

(b) {\it Homological duality}, \cite[Proposition 4.10]{GGOR}, $$ 
D_{\h ,c} \defn \Ext^{n}_{H_{1,c}}(- , H_{1,c}) : \mathcal{O}_{c}^G \longrightarrow \mathcal{O}_{c^{\dagger}}^G,$$ where we have used $\dagger$ to identify right $H_{1,c}$-modules and left $H_{1,c^{\dagger}}$-modules.

The composition of these two dualities can be used to show that $\mathcal{O}_c^G(\h^*\oplus \h)$ is Ringel dual to $\mathcal{O}_c^G(\h\oplus \h^*)$, \cite[Corollary 4.14]{GGOR}.

\medskip
\noindent
(5) Recall that a {\it Serre functor} on a finite triangulated $k$-category $\mathcal{C}$ is a functor $S: \mathcal{C} \longrightarrow \mathcal{C}$ which satisfies $\Hom_{\mathcal{C}}(A,B) \cong \Hom_{\mathcal{C}}(B, SA)^*$ naturally. In \cite[Conjecture 4.12]{MaSt} Mazorchuk and Stroppel conjecture that the left derived functor of $d_{\h^*,c^{\dagger}}\circ D_{\h^*, c}\circ d_{\h, c^{\dagger}}\circ D_{\h,c}$ is a Serre functor for $D^b(\mathcal{O}_c)$, the bounded derived category of $\mathcal{O}_c$. A proof of this would, in particular, show that Hecke algebras are symmetric -- see \ref{KZfunctorintro} and beyond for the appearance of  Hecke algebras. 

\subsection{Primitive ideals} \label{ginzprimthm} For quite general  algebras with triangular decomposition, Ginzburg has proved an analogue of Duflo's theorem for simple complex Lie algebras. Recall that an ideal is called {\it primitive} if it is the annihilator of a simple module: to some extent they play the r\^ole of maximal ideals for noncommutative algebras, and indeed each maximal ideal is primitive. For rational Cherednik algebras the theorem states

\begin{theorem}[{\cite[Theorem 2.3]{ginz}}] An ideal $I$ of $H_{1,c}$ is primitive if and only if $I = \Ann_{H_{1,c}}(L_c(\lambda))$ for some $\lambda \in \irr (G)$.\end{theorem}

The non-trivial part of the theorem is the only if direction: there are many simple $H_{1,c}$-modules which have nothing to do with $\mathcal{O}_c$, but to classify the primitive ideals it is enough to restrict to this small and relatively tractable category. We will illustrate its use in \ref{simplicity2}.

\subsection{Toy example} In the example of \ref{TEdef} we have $\irr(G) = \{+, - \}$ with $$\kappa(c ,+ ) = -\frac{1}{2} + c \quad \text{ and } \quad \kappa(c, -) = - \frac{1}{2} - c.$$ Thus $<_c$ is trivial unless $c\in \frac{1}{2}\Z_+$ or $c\in \frac{1}{2}\Z_-$: in the first case $- <_c +$, whilst in the second $+ <_c -$. Now $\Delta_c(\pm) = \C[x]\otimes \pm$ and the commutation relation \eqref{commmreln} shows that $$y\cdot (x^i\otimes \pm) = \begin{cases} (i\mp 2c)x^{i-1}\otimes \pm \qquad & \text{if $i$ is odd} \\ ix^{i-1}\otimes \pm & \text{if $i$ is even.} \end{cases}$$ Thus $\mathcal{O}_c$ is semisimple unless $c \in \frac{1}{2} +\mathbb{Z}$. If $c = \frac{1}{2} + m$ where $m\geq 0$ then $\Delta_c(-) = L_c(-)$ and there is a short exact sequence \begin{equation} \label{A1SES} 0 \longrightarrow \Delta_c(-) \stackrel{d}\longrightarrow \Delta_c(+) \longrightarrow L_c(+) \longrightarrow 0,\end{equation} where $d(x^i\otimes - ) = x^{2m+1+i}\otimes +.$ In particular $[L_c(+)]_{\bf h} = \sum_{j=-m}^mu^{j}$. Analogous statements hold for $L_c(-)$ when $c= -\frac{1}{2}-m$. 

\subsection{The Dunkl embedding} So far, the manipulation of $\mathcal{O}_c$ has been rather formal. The crucial input which allows to make non-trivial observations and connections is via the ${\sf KZ}$-functor which was studied by Opdam and others, and introduced systematically to the study of rational Cherednik algebras in \cite{GGOR}. In order to define this we must first introduce the Dunkl embedding.

Given $s\in \mathcal{S}$ let $\alpha_s \in \h^*$ be a functional such that $\alpha_s^{-1}(0)$ is the reflecting hyperplane associated to $s$. Similarly to \ref{sympres} we set $\hr = \h \setminus\cup_{s\in \mathcal{S}}\alpha_s^{-1}(0)$ so that $G$ acts freely on $\hr$. Let ${\D}(\hr)$ denote the ring of differential operators on $\hr$ (this can be identified with the localisation of the Weyl algebra ${\D}(\h)$ at the multiplicative set $\{ \delta^i \}_{i> 0}$ where $\delta = \prod_{s\in \mathcal{S}} \alpha_s \in \C[\h]$ is the discriminant). The action of $G$ on $\hr$ extends to an action by algebra automorphisms on $\D (\hr)$, allowing us to form the smash product $\D (\hr)\rtimes G$.
\begin{theorem}[{\cite[Proposition 4.5]{EG}}] \label{dunklemb} There is an injective algebra homomorphism $$\theta_c : H_{1,c} \longrightarrow \D(\hr)\rtimes G$$ induced by $$x\mapsto x, y\mapsto \frac{\partial}{\partial y} + \sum_{s\in \mathcal{S}} c(s) \frac{\alpha_s(y)}{\alpha_s} (s -1), w\mapsto w.$$ On localising at $\{ \delta^i \}_{i > 0}$ we get an isomorphism $H_{1,c}{[}\delta^{-1}{]} \cong \D (\hr) \rtimes G$.
\end{theorem}
This algebra homomorphism is simply describing the action of $H_{1,c}$ on $\C[\h] \cong \Delta_c (\triv)$ where $\triv$ is the trivial representation of $G$. Injectivity can be proved by considering the filtration on $H_{1,c}$ where $\deg \h^* = \deg G = 0$ and $\deg \h = 1$, and the differential operator filtration on $\D (\hr) \rtimes G$ and showing that the associated graded homomorphism is injective.

The second claim is obvious, since inverting $\delta$ has the same effect as inverting all the $\alpha_s$ and so allows us to realise $\partial /\partial y$ in $\im \theta_c$ as $\theta_c ( y - \sum_{s\in \mathcal{S}} c(s) \alpha_s^{-1}\alpha_s(y)(s-1))$.

\subsection{$\!\!\!\!\!$} \label{dunklexp} Theorem \ref{dunklemb} is very important for a number of reasons. The operators $$ \frac{\partial}{\partial y} + \sum_{s\in \mathcal{S}} c(s) \frac{\alpha_s(y)}{\alpha_s} (s -1)$$ are called Dunkl operators; it is an immediate consequence of the theorem that the Dunkl operators commute with one another (since the $y$'s do), a non-obvious fact proved originally by Dunkl. Moreover, the theorem shows that all of the deformation in $H_{1,c}$ lives on the reflecting hyperplanes for the action of $G$ on $\h$. As we shall see, this is a crucial observation.
 
 \begin{remark} $\theta_c$ restricts to an embedding $U_{1,c} \hookrightarrow \D(\hr)^G \cong \D(\hr/G)$ and so, in particular, strengthens Theorem \ref{GKconj} for rational Cherednik algebras.
 \end{remark}
 
 \subsection{$\!\!\!\!\!$}  \label{simplicity2} Now we come to promised illustration of the importance of Theorem \ref{ginzprimthm}: semisimplicity of $\mathcal{O}_c$ implies the simplicity of $H_{1,c}$. To prove this observe that any proper two-sided ideal must be contained in a maximal ideal and hence, by Theorem \ref{ginzprimthm}, must annihilate some $L_c(\lambda)$ for $\lambda \in \irr(G)$. But by semisimplicity, $L_c(\lambda) \cong \Delta_c(\lambda)$. Now if $I = \Ann_{H_{1,c}} (\Delta_c(\lambda))$ is non-zero then $\delta^k \in I$ for some $k\geq 0$ since $H_{1,c}[\delta^{-1}]\otimes_{H_{1,c}}I$ is a two-sided ideal of the simple ring $\D(\hr)\rtimes G$ and hence zero. This is absurd since each $\Delta(\lambda)$ is a free $\C[\h]$-module and thus each $\Delta_c(\lambda)$ is faithful and so the only primitive ideal is $0$. 

\subsection{$\!\!\!\!\!$} \label{simplicity1}
Furthermore simplicity of $H_{1,c}$ implies that $H_{1,c}$ and $U_{1,c}$ are Morita equivalent. Indeed the functors \begin{equation} \label{idemp0}\begin{CD} {\sf I_c} : \Lmod{H_{1,c}} @>>> \Lmod{U_{1,c}}, \qquad M\mapsto eM\end{CD}\end{equation} and \begin{equation} \label{idemp1} \begin{CD}{\sf J_c} : \Lmod{U_{1,c}} @>>> \Lmod{H_{1,c}}, \qquad N \mapsto H_{1,c}e\otimes_{U_{1,c}}N\end{CD}\end{equation} are inverse to another since $H_{1,c}eH_{1,c} = H_{1,c}$ because $H_{1,c}$ is simple and the left hand side is a two-sided ideal, cf. Theorem \ref{defppthm}(4,5).
\subsection{KZ functor} \label{KZfunctorintro}Using the Dunkl embedding we can construct a functor $$\KZ_c : \mathcal{O}_c \longrightarrow \Lmod{\mathcal{H}_q(G)}$$ where $\mathcal{H}_q(G)$ is a Hecke algebra, i.e. the quotient of the group algebra of a generalised braid group, $B_G \defn \pi_1( \hr / G, *)$, by a set of ``Hecke" relations with parameters $q(s)  = \exp (2\pi \sqrt{-1}c(s))$, as presented in \cite{BMR}. If $G$ is a finite Coxeter group then $\mathcal{H}_q(G)$ is the Iwahori--Hecke algebra associated to $G$, a much loved object in Lie theory. 

The functor $\KZ_c$ has similar properties to the Schur functor and to Soergel's functor $\mathbb{V}$, \cite{soe}.
\begin{itemize}
\item $\KZ_c$ is exact: there exists a projective object $P_{\KZ} \in \mathcal{O}_c$ such that $\KZ_c = \Hom_{\mathcal{O}_c}(P_{\KZ}, -)$.  
\item Double centraliser property: $\mathcal{H}_q(G) \cong \End_{\mathcal{O}_c}(P_{\KZ})^{op}$ -- so that $P_{\KZ}$ is in particular a right $\mathcal{H}_q(G)$-module -- and $\mathcal{O}_c$ is equivalent to $\End_{\mathcal{H}_q(G)} (P_{\KZ})$.
\item If $\mathcal{O}_c^{\tor}$ is the full subcategory of $\mathcal{O}_c$ whose objects are annihilated by $\KZ_c$, then $\KZ_c$ induces an equivalence ${\mathcal{O}_c}/{\mathcal{O}_c^{\tor}} \stackrel{\sim}\longrightarrow \Lmod{\mathcal{H}_q(G)}.$
\end{itemize}
So $\mathcal{O}_c$ is a ``highest weight cover" of $\Lmod{\mathcal{H}_q(G)}$ and we can pass properties from one of these categories to the other. At the time of writing, this is the only non-formal input into $\mathcal{O}_c$ that exists for all complex reflection groups and many of the results known about $\mathcal{O}_c$ are consequences of theorems on the representation theory of Hecke algebras. Several examples of this will be given below.

\subsection{$\!\!\!\!\!$} We will only outline the construction of $\KZ_c$ since the details can be found in the final section of Ariki's article in this volume, \cite[Section 4.4]{Ariki}. Let $\lambda \in \irr(G)$. Then $\Delta_c(\lambda)$ extends to a $\D(\hr)\rtimes G \cong H_{1,c}[\delta^{-1}]$-module $\Deltar_c(\lambda) \defn H_{1,c}[\delta^{-1}]\otimes_{H_{1,c}} \Delta_c(\lambda)$. As a $\C[\hr]\rtimes G$-module this is isomorphic to $\C[\hr]\otimes \lambda$, and so is a $G$-equivariant (trivial) vector bundle on $\hr$ of rank $\dim \lambda$. There is, however, still a residue of the action of $\C[\h^*]$ on $\Delta_c(\lambda)$ to be seen on $\C[\hr]\otimes \lambda$, namely a flat $G$-equivariant connection $\nabla_{c,\lambda}$ arising from the various $\partial_y \in \D(\hr)$ for each $y\in \h$, given explicitly by the Dunkl embedding as \begin{equation} \label{connection} \nabla_{c,\lambda} = d -  \sum_{s\in \mathcal{S}} c(s) \frac{d\alpha_s}{\alpha_s}\otimes (s-1)\in \End(\C[\hr]\otimes \lambda).\end{equation}
This is called the {\it Knizhnik-Zamalodchikov connection}, or just {\it $\KZ$-connection}. It has regular singularities and if we take the monodromy of this connection with respect to some basepoint $*\in \hr$ we produce a complex representation of $B_G = \pi_1(\hr/G)$-module. By a rank one calculation, this representation satisfies a Hecke relation and so the resulting monodromy representation is in fact an $\mathcal{H}_q(G)$-representation, to be denoted $\KZ_c(\Delta_c(\lambda))$.

For generic choices of $c$ we have seen already in \eqref{genericss} that $\mathcal{O}_c$ is semisimple and so every object is a direct sum of standard modules. Thus, for such $c$, the above construction can be applied to any object of $\mathcal{O}_c$ and so produces a functor to $\KZ_c : \mathcal{O}_c \longrightarrow \mathcal{H}_q(G)$. The existence of $\KZ_c$ on arbitrary objects in the general case follows by a deformation argument, see \cite[Theorem 5.13]{GGOR}. 

\subsection{$\!\!\!\!\!$} We have factored $\KZ_c$ into two steps:
$$\begin{CD} \KZ_c : \mathcal{O}_c  @>{\sf localisation} >> \mathcal{O}_c^{\sf reg} @> {\sf monodromy}>>\mathcal{H}_q(G)-\modu.\end{CD}$$ The second functor is an equivalence of categories, \cite[Theorem 5.14]{GGOR}, and the localisation is easy to describe formally: it kills the objects in $\mathcal{O}_c$ which are annihilated by some power of $\delta$. In other words, $\mathcal{O}_c^{\tor}$ consists of the objects supported as $\C[\h]$-modules entirely on $V(\delta) \subset \h$, i.e. on the reflecting hyperplanes in $\h$. This explains the third claim in \ref{KZfunctorintro}. The double centraliser property of the second claim is essentially equivalent to the monodromy part of the factorisation being an equivalence, \cite[Theorem 5.15]{GGOR}; finally exactness follows since localisation is exact.

\subsection{$\!\!\!\!\!$} The double centraliser property gives us a much more accurate criterion for the semisimplicity of $\mathcal{O}_c$, and hence by \ref{simplicity2} for the simplicity of $H_{1,c}$.
\begin{corollary} Let $c\in \C[\mathcal{S}]^{\ad G}$ and set $q = \exp (2\pi \sqrt{-1} c) \in \C^*[\mathcal{S}]^{\ad G}$. Then $\mathcal{O}_c$ is semisimple if and only if $\mathcal{H}_q(G)$ is semisimple.
\end{corollary}
If $G$ is a Weyl group then there a good conditions which ensure the semisimplicity of $\mathcal{H}_q(G)$, see for instance \cite[Chapter 9]{geckpfeiffer} and \cite[Theorem 3.29]{Ariki}. Conversely, for arbitrary $G$ the corollary combined with \eqref{genericss} shows that if if the subgroup of $\C^*$ generated by the $q(s)$ is torsion-free then $\mathcal{H}_q(G)$ is semisimple, \cite[Theorem 3.4]{RouquierqSch}. 
 
 \subsection{Toy example} Continuing the example of \ref{TEdef} we have that $B_G = \langle T \rangle$ where we take the basepoint $*$ of $\hr/G$ to be the $G$-orbit of $1$ and $T$ corresponds to the path $\exp(\pi\sqrt{-1}t)$ with $t: 0 \rightarrow 1$ in $\hr = \C^*$. From \eqref{connection} the connection corresponding to $+$ is just $d/dx$ -- its solutions are the constant functions and they have no monodromy --; for $-$ the connection is $d/dx + 2c/x$ -- its solutions are $x^{-2c}$. This last function has monodromy $-\exp(\pi\sqrt{-1}c)$. Thus, setting $q = \exp(2\pi \sqrt{-1}c)$, we see that $\KZ_c(\Delta_c(+))$ has $T$ acting as $1$ and $\KZ_c(\Delta_c(-))$ has $T$ acting as $-q$: both are representations of the Hecke algebra $\mathcal{H}_q(G) = \C[T] /\langle (T-1)(T+q) \rangle$. If $c\notin \frac{1}{2}+\Z$ then $q\neq 1$ and $\KZ_c(\Delta_c(+))$ and $\KZ_c(\Delta_c(-))$ are distinct; if $c\in \frac{1}{2} + \mathbb{Z}$ then $q = -1$ and $\KZ_c(\Delta_c(+)) = \KZ_c(\Delta_c(-))$ and applying $\KZ_c$ to the sequence \eqref{A1SES} shows that $\KZ_c(L_c(+)) = 0$.
 
\subsection{Shifting} \label{rouquier} The close relationship between $\mathcal{O}_c$ and $\mathcal{H}_q(G)$, particularly the double centraliser property, begs an obvious question:

\medskip

\begin{center} Suppose $c, c' \in \C[\mathcal{S}]^{\ad G}$ are such that $q = \exp (2\pi \sqrt{-1} c) = \exp (2\pi \sqrt{-1} c') = q'$, \\so that $\mathcal{H}_q(G) = \mathcal{H}_{q'}(G)$. Is it true that $\mathcal{O}_c$ and $\mathcal{O}_{c'}$ are equivalent?\end{center}
 
 \medskip
 
 \noindent
The answer is, unfortunately, no. In fact, it is easy to see why we could expect a difference between $\mathcal{O}_c$ and $\mathcal{O}_{c'}$. The ordering on $\mathcal{O}_c$ given in Theorem \ref{HWthm} depends on the ordering $<_c$ on $\irr(G)$ which varies with $c\in \C[\mathcal{S}]^{\ad G}$. The variation is linear and so it is possible to choose $c, c'$ such that $q = q'$, but for which the associated orderings are different. Thus, the best we could hope is the following wonderful theorem, \cite[Theorem 5.5]{RouquierqSch}, which is explained in detail in \cite[Section 4]{Ariki}.

\begin{theorem}[\cite{RouquierqSch}] Suppose that $c\in \C[\mathcal{S}]^{\ad G}$ and let $d\in \Z[\mathcal{S}]^{\ad G}$. Then there is an equivalence of categories between $\mathcal{O}_c$ and $\mathcal{O}_{c+d}$ if the orderings $<_c$ and $<_{c+d}$ on $\irr (G)$ are the same.
\end{theorem}

(There are actually some mild restrictions on $c$ -- see \cite{RouquierqSch} for details.) This theorem is actually a special case of a more general result which Rouquier proves concerning highest weight covers of Hecke algebras. Using that general result Rouquier also proves that $\mathcal{O}_c^{S_n}$ is equivalent to the category of representations of the $q$-Schur algebra $S_q(n,n)$, \cite[Theorem 6.8]{RouquierqSch}.

\subsection{$\!\!\!\!\!$} Since $\mathcal{O}_c$ controls much of the structure of $H_{1,c}$, for instance through Theorem \ref{ginzprimthm}, we might hope that an equivalence between $\mathcal{O}_c$ and $\mathcal{O}_{c'}$ would be the shadow of an equivalence between $H_{1,c}$ and $H_{1,c'}$. There is much evidence to suggest this, but at the moment only special cases are understood, relying on remarkable work of Opdam and Heckman which we now explain.

\subsection{$\!\!\!\!\!$} Assume that $G$ is a finite Coxeter group. Recall the $e\in \C G$ is the trivial idempotent. Set $e_- = |G|^{-1}\sum_{g\in G} (-1)^{\ell (g)} g \in \C G$, the sign idempotent, and $U_{1,c}^- \defn e_- H_{1,c} e_-$. 

The Dunkl embedding $\theta_c: H_{1,c} \longrightarrow \D(\hr)\rtimes G$ embeds each $H_{1,c}$ into the same algebra of differential operators on $\hr$ and so we can compare rational Cherednik algebras for different $c$. From now on we will identify $H_{1,c}$ with its image in $\D(\hr)\rtimes G$ under $\theta_c$. 

Since $G$ is a finite Coxeter group, there is a $G$-invariant element of degree two in $\C[\h^*]$, namely $\sum_{i=1}^n y_i^2$ where the $y_i$ form a basis of $\h$. Therefore, under the Dunkl embedding, we find elements deforming the Laplacian $${\sf L}_c = \sum_{i=1}^n \theta_c(y_i^2e) \in e(\D(\hr)\rtimes G)e \text{ and } {\sf L}^-_c = \sum_{i=1}^n \theta_c(y_i^2e_-) \in e_-(\D(\hr)\rtimes G)e_-.$$ The crucial observation of \cite{opdam} and \cite{heckman} is that \begin{equation} \label{HOobs}{\sf L}_{c}  = \delta^{-1} {\sf L}_{c+1}^- \delta \in e(\D(\hr)\rtimes G)e.\end{equation} Thus the two algebras $U_{1,c}$ and $\delta^{-1} U^-_{1,c+1}\delta$ have many elements in common, namely ${\sf L}_c$ and the elements of $\C[\h]^G e$. If $U_{1,c}$ is simple then \cite[Theorem 5]{LevSta} shows that $U_{1,c}$ is generated by $\C[\h]^G e$ and ${\sf L_c}$ and thus for generic choices of $c$ thanks to  \ref{simplicity2}, it follows that $U_{1,c} = \delta^{-1}U^-_{1,c+1}\delta$ in $\D(\hr)\rtimes G$ . A deformation argument \cite[Proposition 5.4]{BEGfd} extends this equality to all parameters. It follows that $eH_{1,c+1} \delta e = eH_{1,c+1} e_- \delta$ is a $(U_{1,c+1}, U_{1,c})$-bimodule and we have a functor
 \begin{equation} \label{sphershift} {\sf S}_c : U_c - \modu \longrightarrow U_{c+1} -\modu , \qquad M \mapsto eH_{1,c+1} \delta e \otimes_{U_{c}} M.\end{equation} Combining this with functors ${\sf I}_c$ and ${\sf J}_{c+1}$ of \eqref{idemp0} and \eqref{idemp1} we also find \begin{equation} \label{shift} {\sf S}_c = {\sf J}_{c+1}\circ {\sf S}_c \circ {\sf I}_{c} : H_{1,c}-\modu \longrightarrow H_{1,c+1}-\modu .\end{equation} These are called the {\it Heckman-Opdam shift functors}. 
 
 \subsection{$\!\!\!\!\!$} \label{shiftequiv} Now these functors cannot be equivalences in general because of the comments made in  \ref{rouquier}. However, there are still many examples where the orderings $<_c$ and $<_{c+1}$ on $\irr (G)$ are the same. In these cases one should ask whether ${\sf S}_c$ is an equivalence of categories. The answer is unknown in general, but there are important special cases for which we have a positive answer:
\begin{itemize}
\item generic $c$ (by arguments analogous to  \ref{simplicity1});
\item $c$ is ``large enough" (see \cite[Proposition 4.3]{BEGfd});
\item $c =m + 1/h $ where $h$ is the Coxeter number of $G$ and $m\in \Z$ (see \cite[Section 4]{gordc} and \cite[Proposition 4.3]{BEGfd});
\item $G = {S}_n$ for $c\notin (-1,0)$ (by \cite[Theorem 3.3]{gordst1} and \cite[Corollary 4.2]{BE});.
\end{itemize}
\subsection{$\!\!\!\!\!$} \label{generalshift}The shift functors can be constructed in more generality: instead of adding $1$ to $c$ we can add $1$ to a given conjugacy class of reflections only, \cite[Section 4]{BEGqi}; we can also deal with the infinite family $G(m,d,n)$ of complex reflection groups, \cite[Section 3]{dunklopdam}. There is, however, less known about when they are equivalences, see \cite[Chapter 4]{vale}.

\subsection{Finite dimensional representations} To illustrate some of the power of the $\KZ$-functor and the Heckman-Opdam shift functors we will construct finite dimensional representations of the rational Cherednik algebras associated to Coxeter groups.

\subsection{$\!\!\!\!\!$} \label{fd1}We look first for the easiest possible representation of $H_{1,c}$ which could possibly exist, namely a one-dimensional space $\C$ which carries the trivial action of $G$. By \ref{Ocomments}(2) $\h$, and by symmetry $\h^*$, must act trivially on $\C$ and thus we need only make sure that the defining relation of Theorem \ref{PBWthm} is satisfied on $\C$, i.e. that for all $v,w\in V$ $$0 = \omega_{V^*}(v,w) + \sum_{s\in \mathcal{S}} c(s)\omega_s(v,w)s .$$ If $c$ is constant then it is straightforward to check that this is satisfied if and only if $2|\mathcal{S}|c =n$, in other words $c = \rank(G)/\#\text{roots}.$ This equals $1/h$, where $h$ is the Coxeter number of $G$. Thus we have proved that $L_{1/h}(\triv) = \C$ is one-dimensional.

If we apply the shift functor $m$ times to $L_{1/h}(\triv)$ then, thanks to  \ref{shiftequiv}, we will produce a finite dimensional irreducible representation ${\sf S}_{m-1 +1/h}\circ \cdots \circ {\sf S}_{1+ 1/h}\circ {\sf S}_{1/h} (L_{1/h}(\triv))$ and it's not hard, by considering $\bf h$-weight spaces, to check that this is $L_{m +1/h}(\triv)$.

\subsection{$\!\!\!\!\!$} \label{fd2} For the moment, set $c = m+1/h$. In order to say more about $L_{c}(\triv)$ we need to use the $\KZ$-functor. This sends $\mathcal{O}_{c}$ to $\Lmod{\mathcal{H}_{q}(G)}$ where $q$ is the constant function taking value $\exp(2\pi\sqrt{-1}/h)$. Now $\KZ_{c}$ cannot possibly be an equivalence since $\KZ_{c}(L_{c}(\triv)) =0$ since $L_{c}(\triv)$ is finite dimensional. However, this turns out to be the only information that is lost. This is because the parameter $q = \exp(2\pi \sqrt{-1}/h)$ is very special for Iwahori--Hecke algebras: it produces the ``first" and simplest examples of $\mathcal{H}_q(G)$ which are not semisimple. They have been studied in detail in \cite{BGK} and \cite{mueller} leading to a description of their decomposition matrices. It is shown in \cite[Lemma 4.3]{gordc} and \cite[Proposition 3.8]{BEGfd} that the simple nature of these decomposition matrices is enough to insure that (up to scalar multiplication) there is a unique non-trivial $H_{1,c}$-homomorphism $d_k: \Delta_{c}(\wedge^k \h) \longrightarrow \Delta_{c}(\wedge^{k-1}\h)$ for any $1\leq k \leq n = \dim \h$ and no other homomorphisms between standard modules. It follows that we have a ``BGG-resolution" of $L_c(\triv)$ \begin{equation}\label{BGGresol}0 \longrightarrow  \Delta_{c}(\wedge^{n} \h)\stackrel{d_{n}}\longrightarrow \cdots \stackrel{d_2}\longrightarrow \Delta_{c}( \h) \stackrel{d_1}\longrightarrow \Delta_{c}(\triv) \longrightarrow L_c(\triv) \longrightarrow 0,\end{equation} and that $\Delta_c(\lambda)$ is irreducible if $\lambda \neq \wedge^k \h$ for some $k$. Since ${\bf h}\in H_{1,c}$ each $d_k$ respects decompositions into $\bf h$--eigenspaces and so we see from \eqref{BGGresol} and \eqref{wtspchar} $$[L_c(\triv)]_{\bf h} = \sum_{k=0}^n (-1)^k [\Delta_c(\wedge^k \h)]_{\bf h} = \sum_{k=0}^n (-1)^k \binom{n}{k}\frac{u^{\kappa(c, \wedge^k \h)}}{(1-u^{-1})^nm}.$$
An elementary calculation gives $\kappa(c, \wedge^k\h) =   mh(k + n/2) + k$. Thus $$[L_c(\triv)]_{\bf h}  = u^{-nmh/2} (1+u+ \cdots + u^{mh})^n.$$ In particular, setting $u=1$ we see that $\dim L_{m+1/h} (\triv) =  (mh+1)^n.$

\subsection{$\!\!\!\!\!$} In this example the entire quiver and relations corresponding to $\mathcal{O}_c$ can be calculated. Ignoring the simple blocks of $\mathcal{O}_c$ (which correspond to $\lambda \neq \wedge^k \h$) this produces the algebra corresponding to the category of perverse sheaves on $\mathbb{P}^{n}$ which are constructible with respect to the standard stratification $(\C^0 \cup \C^1 \cup \C^2 \cup \cdots \cup \C^{n})$, or equivalently to parabolic category $\mathcal{O}$ for the Lie algebra $\mathfrak{gl}_{n+1}(\C)$ and maximal parabolic with Levi $\mathfrak{gl}_n(\C)$. A general version of these results is given in \cite[5.2.4]{RouquierqSch}.

\subsection{$\!\!\!\!\!$} There has been considerable work towards the classification of finite dimensional $H_{1,c}$-representations, starting with \cite{BEGfd} which showed that for $G=S_n$, there existed a (unique) finite dimensional representation if and only if $c= r/n$ with $(r,n) = 1$. Chmutova studied the case where $G$ is a dihedral group, \cite{chmut}, using systematically the $\KZ$-functor and representations of the corresponding Hecke algebras. The most general results for Weyl groups have been found recently by Varagnolo--Vasserot, \cite{VV}: they prove in the equal parameter case that $H_{1,c}$ has a finite dimensional representation if and only if $c = r/s$ where $(r,s)=1$ and $s$ is an {\it elliptic} number, see \cite[p.9]{VV} for the list of these. The Coxeter number is always elliptic, and in case $G=S_n$ it is the only elliptic number. Varagnolo--Vasserot's proof is geometric: they show that any irreducible finite dimensional $H_{1,c}$-representation can be considered as an irreducible ``spherical" representation of a double affine Hecke algebra; they then construct such representations geometrically using the equivariant $K$-theory of affine Springer fibres.

Finite dimensional representations for unequal parameters are not yet understood; nor are they for the majority of complex reflection groups, but see \cite{Gr1} and \cite{Gr2}.

\subsection{Application: Coinvariants} \label{coinvariants} The work done in \ref{fd1} and \ref{fd2} has an application in combinatorics and invariant theory and points to a geometric interpretation of the rational Cherednik algebras which is quite different from the $K$-theoretic construction of Varagnolo--Vasserot.

Recall from \eqref{coinvariant} the definition of the coinvariant ring $\co{\h\oplus \h^*}{G}$. Since $G$ does {\it not} act by complex reflections on $\h\oplus \h^*$ we expect that $\co{\h\oplus \h^*}{G}$ is larger than $|G|$. This is borne out by the following theorem.
\begin{theorem}[\cite{gordc}] Let $G$ be a finite Coxeter group of rank $n$ and with Coxeter number $h$. Then there is a graded $\C$-algebra $R$ with an action of $G$ and a surjective $G$-equivariant algebra homomorphism $$\pi_G : \co{\h\oplus \h^*}{G}\longrightarrow R$$ such that the Poincar\'e polynomial of $R$ is $u^{-nh/2} (1+u+\cdots + u^h)^n$.
\end{theorem}
We have enough to prove this here. There is an isomorphism $L_{1+1/h}(\triv) \cong {\sf S}_{1/h}(\C) \cong H_{1, 1+1/h} e_- \otimes_{U_{1/h}} \C$. Give $L_{1+1/h}(\triv)$ the filtration $F^p(L_{1+1/h}(\triv)) \defn F^p(H_{1,1+1/h})\cdot e_-\otimes \C$. Then there is a surjective mapping $$\gr (H_{1,1+1/h}) e_- \otimes_{\gr U_{1/h}} \gr \C \longrightarrow \gr( H_{1, 1+1/h} e_- \otimes_{U_{1/h}} \C) = \gr L_{1+1/h}(\triv).$$ By  \eqref{PBW}, the left hand side is isomorphic to $(\C[\h  \oplus \h^*]\rtimes G)e_-\otimes_{\C[\h\oplus \h^*]^G} \C$. Now $(\C[\h  \oplus \h^*]\rtimes G)e_- \cong \C[\h  \oplus \h^*]\otimes \sign$ and so we further refine the left hand side to $(\C[\h  \oplus \h^*]\otimes_{\C[\h\oplus \h^*]^G} \C) \otimes \sign$. But, by definition, this is just $\co{\h\oplus \h^*}{G}\otimes \sign$ and so it follows that $R = \gr L_{1+1/h}(\triv)\otimes \sign$ is a quotient of the commutative algebra $\co{\h\oplus \h^*}{G}$ and hence the algebra we are looking for.

\subsection{$\!\!\!\!\!$} A stronger version of this theorem in the case $G = S_n$ goes under the name of the $(n+1)^{n-1}$ theorem, and was proved by Haiman in \cite{haimandiag}: it asserts that the kernel of $\pi_{S_n}$ is zero. Haiman's proof is geometric, relying on Hilbert schemes, the $n!$ theorem and work by Bridgeland--King--Reid on the higher dimensional homological McKay correspondence. The result above was conjectured in \cite{haimanconj} for all Coxeter groups. It is known that outwith the symmetric group, the mapping $\pi_G$ can fail to be injective and that geometrical reasoning along the lines of the $n!$ theorem will also fail since, as we will see in the next chapter, symplectic resolutions do not always exist for $\h\oplus \h^*/G$. 

Generalisations of the theorem exist for other complex reflection groups, \cite{richard}, and another proof for Weyl groups -- using double affine Hecke algebras -- was given in \cite{cherednikcoinv}.
\section{Specific case: rational Cherednik algebras II ($t=0$)} \label{t0RCA}
Throughout this section we continue to focus on the case when $G< GL(\h)$ is a complex reflection group of rank $n$ and $V = \h\oplus \h^*$. However, we will assume that $t = 0$ so that by Theorem \ref{centrethm}  $H_{0,c}$ is a finite module over its centre and there is a link between the geometry of $Y_c \defn \Spec Z(H_{0,c})$ and the representation theory of $H_{0,c}$. We will exploit this to determine when $Y_c$ is smooth and when $V/G$ admits a symplectic resolution. We also discuss a conjectural connections with Rouquier families for complex reflection groups.

\subsection{Reduced rational Cherednik algebras}
Recall from \ref{noflatcentre} the family of algebras $H_{0,c}(y)$ and that its members do not vary continuously with $y\in Y_c$. This problem can be circumvented for rational Cherednik algebras thanks to the useful observation of  \cite[Proposition 4.15]{EG} that $\C[\h]^G\otimes \C[\h^*]^G \subset Z(H_{0,c}).$ Remarkably this subalgebra does not depend on $c\in \C[\mathcal{S}]^{\ad G}$. Moreover, by \eqref{PBW} and  \ref{STC} $H_{0,c}$ is a free $\C[\h]^G\otimes \C[\h^*]^G$-module of rank $|G|^3$. Thus we have a finite morphism \begin{equation} \label{ups} \Upsilon_c : Y_c \longrightarrow \h/G\times \h^*/G,\end{equation} where, thanks to Theorem \ref{STC}, the codomain is an affine space of dimension $2n$. We also have a flat family of finite dimensional algebras $$H_{0,c}(p,q) \defn \frac{H_{0,c}}{\mathfrak{m}_{p,q} H_{0,c}}$$ where $(p,q)$ runs over $\h/G\times \h^*/G$ and $\mathfrak{m}_{p,q}$ is the corresponding maximal ideal of $\C[\h]^G \otimes \C[\h^*]^G$. These $|G|^3$-dimensional algebras are called {\it reduced rational Cherednik algebras}.

\subsection{$\!\!\!\!\!$} The irreducible representations of $H_{0,c}(p,q)$ consist of the irreducible representations of all $H_{0,c}(y)$ where $y\in \Upsilon_c^{-1}(p,q)$. So, following \ref{killcen}, we can understand the irreducible representations of $H_{0,c}$ by study those of the flat family of reduced rational Cherednik algebras.

\subsection{$\!\!\!\!\!$}
There is a $\C^*$-action on $H_{0,c}$ with $\deg \h = - \deg \h^* =1$ and $\deg G = 0$. This action passes to $Z(H_{0,c})$, $Y_c$ and $\C[\h]^G\otimes \C[\h^*]^G$. If we want to study the geometry of $Y_c$ or the representation theory of $H_{0,c}$ then it makes sense to look for $\C^*$-fixed points in $Y_c$ and graded irreducible $H_{0,c}$-representations since much interesting structure will be attracted toward these. The $\C^*$-action on $\h/G\times \h^*/G$ is given by $\lambda \cdot (p,q) = (\lambda p, \lambda^{-1}q)$ for all $\lambda \in \C^*$ and so the only fixed point in $\h/G\times \h^*/G$ is the origin $(0,0)$. Thus the {\it restricted rational Cherednik algebra} $$\overline{H}_{0,c} \defn H_{0,c}(0,0)$$ will be rather important in the study of $H_{0,c}$. It is a finite dimensional symmetric graded algebra with $G$-action and triangular decomposition $$\overline{H}_{0,c} \stackrel{\sim}\longrightarrow \co{\h}{G}\otimes\C G\otimes \co{\h^*}{G}.$$

\begin{remark} It is not true that $H_{0,c}$ is the only interesting reduced rational Cherednik algebra, but it is the only {\it canonically defined} interesting quotient. It can happen that the ramification locus of $\Upsilon_c$ is a non-empty $\C^*$-stable subvariety of $\h/G\times \h^*/G$ which does not meet $(0,0)$, producing nilpotence in some of the $H_{0,c}(p,q)$, but not in $\overline{H}_{0,c}$; the ramification locus will, however, vary with $c$. 
\end{remark}

\subsection{``Category $\mathcal{O}$"}
We can mimic the construction of standard modules in $\mathcal{O}_c$ to construct a family of $H_{0,c}$-modules.

\begin{definition}
Let $\lambda \in \irr(G)$. The {\it baby Verma module associated to $\lambda$} is  $$\overline{\Delta}_c(\lambda) \defn H_{0,c}\otimes_{\C[\h]^G\otimes \C[\h^*]\rtimes G} \lambda ,$$ where $\C[\h]^G\otimes \C[\h^*]$ acts on $\lambda$ via $p \cdot v = p(0,0) v$ for any $p\in \C[\h]^G\otimes C[\h^*]$ and $v\in \lambda$. 
\end{definition}
Since the maximal ideal $\mathfrak{m}_{0,0}$ of $\C[\h]^G\otimes \C[\h^*]^G$ annihilates $\overline{\Delta}_c(\lambda)$ we see that $\overline{\Delta}_c(\lambda)$ is an $\overline{H}_{0,c}$-module. As a graded $\C[\h]\rtimes G$-module $\overline{\Delta}_c(\lambda)$ is isomorphic to $\co{\h}{G}\otimes \lambda$ by the PBW isomorphism \eqref{PBW}; by Theorem \ref{STC}, it is isomorphic to $\C G^{\oplus \dim \lambda}$ as a $G$-representation.

\subsection{$\!\!\!\!\!$} \label{HNprop} Standard arguments on graded algebras with triangular decompositions, \cite{HolNak}, show that 
\begin{enumerate}
\item Each $\overline{\Delta}_c(\lambda)$ has a simple head $\overline{L}_c(\lambda)$ which is graded.
\item The set $\{ \overline{L}_c(\lambda) : \lambda\in \irr(G)\}$ is a complete set of irreducible $\overline{H}_{0,c}$-modules up to isomorphism.
\end{enumerate}
Thus it becomes a crucial to understand the composition multiplicities $[\overline{\Delta}_c(\lambda): \overline{L}_c(\lambda)]$. This is answered only in a few simple cases, but may again be related to the Hecke algebras $\mathcal{H}_q(G)$.

\subsection{Singularities in the centre} Theorem \ref{smoothpartcentre} tells us that it is important to know whether $Y_c$ is smooth. If it is, then describing $H_{0,c}$ is equivalent to describing $Y_c$. Now Theorem \ref{azumayathm}(3) tells us that $Y_c$ is smooth if and only if all irreducible representations of $H_{0,c}$ carry the regular representation as a $G$-module. We now have a source of irreducible $H_{0,c}$-representations, the $\overline{L}_c(\lambda)$'s: they are quotients of baby Verma modules whose $G$-structure we are able to understand perfectly. 

\begin{theorem}[{\cite[Corollary 1.14]{EG}, \cite[Proposition 7.3]{babyv}, \cite[Theorems 3.1 and 4.1]{gwyn}}] \label{singularities} Let $G$ be an irreducible complex reflection group. There exists a choice of parameter $c\in \C[\mathcal{S}]^{\ad G}$ such that $Y_c(G)$ is smooth if and only if $G = \mu_{\ell}\wr S_n$ or $G$ is the binary tetrahedral group. 
\end{theorem}
Here $\mu_{\ell}\wr S_n$ is called $G(\ell, 1, n)$ in the Shephard--Todd classification while the binary tetrahedral group is the exceptional group $G_4$. The binary tetrahedral group also appears as the kleinian group associated to $\tilde{E}_6$ in Figure 1 -- the imaginary root $\delta$ for $\tilde{E}_6$ shows that this group has three two-dimensional irreducible representations: one of these realises it as a subgroup of $SL(2,\C)$; the other two give its action on $\h$ and $\h^*$ as a complex reflection group.

\subsection{$\!\!\!\!\!$}
The proof of the theorem proceeds in two steps. If $G = G(\ell, 1,n)$ -- or more generally $G = \Gamma \wr S_n$ for some $\Gamma < SL(2,\C)$ -- Etingof and Ginzburg present $Y_c(G)$ as a smooth quiver variety for generic values of $c$, thus proving most of the ``if" direction. (This will be discussed in Chapter \ref{QUIVERS}.) The binary tetrahedral example is dealt with by Bellamy by showing that all irreducible $H_{0,c}$ have dimension $|G|$. 

To prove the ``only if" direction it is only necessary to show that there exists some $\lambda\in\irr(G)$ such that $\dim \overline{L}_c(\lambda) < |G|$. Well, if all $\overline{L}_c(\lambda)$'s were $|G|$-dimensional then it would follow quickly from Theorem \ref{azumayathm} that there would be no non-split extension between the different $\overline{L}_c(\lambda)$'s. Since each baby Verma module is indecomposable this would mean that $\overline{L}_c(\lambda)$ would be the only composition factor of $\overline{\Delta}_c(\lambda)$ for each $\lambda \in \irr(G)$. Looking at the graded $G$-action on $\overline{\Delta}_c(\lambda)$ and on $\overline{L}_c(\lambda)$ one sees, case-by-case, that this is impossible for all the complex reflection groups except $G(\ell , 1 ,n)$ and $G_4$.

\subsection{Application: Symplectic resolutions} Theorem \ref{singularities} has a pleasant consequence.
\begin{theorem}[\cite{GK}, \cite{wang}, \cite{gwyn}] \label{weresolve} Let $G\leq GL(\h)$ be a complex reflection group acting on $V = \h\oplus \h^*$. Then $V/G$ has a symplectic resolution if and only if $G = G(\ell , 1 ,n)$ or $G = G_4$.
\end{theorem}

In \cite[Sections 1.3 and 1.4]{wang} Wang observes that if $G = G(\ell, 1,n)$ (or more generally $G= \Gamma \wr S_n$ for some $\Gamma < SL_2(\C)$) then $V/G$ has a symplectic resolution given by the Hilbert scheme of $n$ points on the minimal resolution of the Kleinian singularity $\C^2 / \Gamma$. In \cite[Corollary 4.2]{gwyn}, Bellamy uses work of Kawamata and the minimal model programme to show the existence of a symplectic resolution for $G = G_4$. This completes the ``if" direction.

The ``only if" direction is proved in \cite{GK} by deformation theoretic arguments. It is crucial that the $\C^*$-action on $V/G$ induced from dilation on $V$ contracts all orbits to the zero orbit and rescales the symplectic form $\omega_{V/G}$ by a positive weight. This allows formal deformation theoretic arguments to be globalised, as in Remark \ref{gradeddef}.
There are two key results in \cite{GK} for us. The first, \cite[Theorem 1.13]{GK}, states that if there is a symplectic resolution $\pi : X \longrightarrow V/G$ then there is a family of symplectic resolutions $\pi_b : X_b \longrightarrow \pi_b(X_b)$ over an affine base space $B = H^2(X,\C)$ whose special fibre $\pi_0$ is  $\pi$ and whose generic fibre is an isomorphism. In particular, this implies that the generic deformation $\pi_b(X_b)$ of $V/G$ must be smooth since $X_b$ is. Then \cite[Theorem 1.18]{GK} states that the family of deformation $(Y_c)_c$ of $V/G$ constructed by rational Cherednik algebras is ``big enough" to meet the generic $\pi_b(X_b)$.
 
Putting these together we see that if $V/G$ has a symplectic resolution, then $Y_c$ must be smooth for some choice of $c$. Combined with Theorem \ref{singularities}, this proves ``only if" direction.

\subsection{$\!\!\!\!\!$} \label{res=def} The theorem can be paraphrased by the five-word-phrase ``the resolution is the deformation". For example, if $G=G (\ell, 1, n)$ then there is a resolution with a hyper-K\"ahler structure -- it is a quiver variety -- and a rotation of the complex structure produces the generic deformation $Y_c$. Thus, the resolution and the deformation are diffeomorphic. The case of $G_4$ is not yet understood in this context. In particular, it would be a good idea to describe the corresponding symplectic resolution in terms of $\C[V]\rtimes G_4$-representations and also to find an explicit description of a  hyper-K\"ahler structure on $Y_c(G_4)$.  

\subsection{Application: Families for complex reflection groups}
 The blocks of $\overline{H}_{0,c}$ correspond to the decomposition of the centre $Z(\overline{H}_{0,c})$ into local algebras. We have a natural homomorphism $$ Z(H_{0,c}) \hookrightarrow Z(\overline{H}_{0,c}),$$ and its kernel is $\mathfrak{m}_{0,0}$. In other words we have an embedding $\C [\Upsilon_c^{-1} (0,0)]\longrightarrow Z(\overline{H}_{0,c})$; however, in general this is not surjective. But there is a ring theoretic result of M\"uller which helps us, see \cite[Chapter III.9]{ken^2}.
\begin{theorem} The image of $\C[\Upsilon_c^{-1}(0,0)]$ in $Z(\overline{H}_{0,c})$ determines the block decomposition of $\overline{H}_{0,c}$. More precisely, each primitive idempotent of $\C[\Upsilon_c^{-1}(0,0)]$ lifts to a unique primitive central idempotent of $\overline{H}_{0,c}$.
\end{theorem} 
This theorem can be given a geometric interpretation. Since $\Upsilon_c$ is a finite $\C^*$-equivariant morphism and $(0,0)$ is the unique $\C^*$-fixed point of $\h/G\times \h^*/G$, the points of $\Upsilon_c^{-1}(0,0)$ are the $\C^*$-fixed points of $Y_c$. Thus the $\C^*$-fixed points of $Y_c$ label the blocks of $\overline{H}_{0,c}$. 

\subsection{$\!\!\!\!\!$} By \ref{HNprop}(1) each $\overline{\Delta}_c(\lambda)$ is indecomposable, so it belongs to a block of $\overline{H}_{0,c}$. Hence we have a surjective mapping $$\Theta_c: \irr(G) \longrightarrow \Upsilon_c^{-1}(0,0)$$ associating a fixed point of $Y_c$ to each irreducible representation of $G$. This induces a partition of $\irr(G)$ which depends on $c\in\C[\mathcal{S}]^{\ad G}$: $\lambda, \mu \in \irr(G)$ belong to the same class if and only if $\Theta_c(\lambda) = \Theta_c(\mu)$. We call this partition the {\it $CM_c$-partition} of $\irr (G)$.

\subsection{$\!\!\!\!\!$} \label{gormartconj} We refer to the survey of Martino, \cite{mosurvey}, for the definitions of Rouquier families -- natural partitions of $\irr(G)$ which appear in Kazhdan--Lusztig theory -- and for the precise parametrisation of these families.

\begin{conjecture}[\cite{gordcells}] \begin{enumerate}
\item The $CM_c$-partition of $\irr (G)$ agrees with the partition of $\irr (G)$ into Rouquier families. 
\item Let $y\in \Upsilon_c^{-1}(0,0)$. Then the dimension of $\C[\Upsilon_c^{-1}(0,0)_y]$ equals the dimension of the corresponding Rouquier block.
\end{enumerate}
\end{conjecture}

There is quite a lot of evidence for this conjecture. Part (2) holds for all smooth points in $\Upsilon^{-1}(0,0)$ and Part (1) holds for $G(\ell , 1, n)$ and any choice of parameter, \cite{gordquiv}, \cite{gordcells} and \cite{mosurvey}. The proof of Theorem \ref{singularities} shows that this partition is trivial if and only if $Y_c$ is smooth, so Part (1) of the conjecture predicts that $G_4$ should be the only exceptional complex reflection group for which there exists a choice of parameters which produces trivial Rouquier families: this is true thanks to recent work of Chlouveraki, \cite{chlou}.

\subsection{$\!\!\!\!\!$} There is no conceptual confirmation of this conjecture in sight at the moment. Rouquier families contain information on the integral representation theory of Hecke algebras for complex reflection groups; it is not clear how the representation theory of $H_{0,c}$ or the geometry of $Y_c$ sees this. One approach for Weyl groups may be to consider an appropriate semiclassical limit of the $\KZ_c$ functor, relating certain $H_{0,c}$-representations to the representations of Lusztig's aymptotic Hecke algebra, but this is unclear at the moment.

\subsection{Toy example} For $G = \mu_2$ it is immediate from the relation $[y,x] = -2cs$ that $x^2$ and $y^2$ are central. The baby Verma modules $\overline{\Delta}_c(\pm) = \C[x]/\langle x^2 \rangle \otimes \pm$ have the action $y\cdot (1\otimes \pm) = 0$ and $y \cdot (x\otimes \pm ) = \mp 2c\otimes \pm$ and are thus simple if and only if $c\neq 0$. If $c = 0$ then $\overline{\Delta}_c(\pm)$ are filtered by one copy of $\overline{L}_c(+) = + $ and one copy of $\overline{L}_c(-) = -$ with $x$ and $y$ acting as zero on both. The final generator of $Z(H_{0,c})$ is ${\bf h}$ and by \ref{firsttoyex} the defining relation is $x^2y^2 = (h-4c)(h+2c)$. Thus this is smooth if and only if $c\neq 0$, as predicted.

\section{Specific case: quivers and hamiltonian reduction}
\label{QUIVERS}
We move on to study the second family of \ref{groups}, $G = \Gamma^n\rtimes  S_n = \Gamma \wr S_n$, where $\Gamma$ is a finite subgroup of $SL(2,\C)$, acting on $V = (\C^2)^n$. The McKay correspondence relates these groups to affine Dynkin diagrams and this allows us to describe $H_{t,c}(G)$ as a deformation of a tensor product of deformed preprojective algebras. The representation at $t = 0$ is quite well understood in terms of quiver varieties; for $t\neq 0$ this is still under investigation via a process called quantum hamiltonian reduction.

\subsection{$\!\!\!\!\!$} When $G = \Gamma\wr S_n$ the conjugacy classes of symplectic reflections in $G$ fall into two families:
\begin{enumerate}

\item[$(S)$] The elements $s_{ij}\gamma_i\gamma_j^{-1}$ where $1\leq i<j\leq n$ and $\gamma\in \Gamma$. Here $\gamma_i$ indicates the element of $\Gamma^n\leq G$ whose $i$th component is $\gamma$ and all of whose other components are the identity.
\item[($\Gamma$)] The elements $\gamma_i$ for $1\leq i \leq n$ and $\gamma\neq 1$.
\end{enumerate}
Here the elements of $(S)$ comprise one conjugacy class, while conjugacy classes in $(\Gamma)$ correspond to the non-trivial conjugacy classes in $\Gamma$.

We will write $\underline{c} = (c_1, {c})$ for an element of $\underline{c}\in \C[\mathcal{S}]^{\ad \Gamma_n}$ where $c_1$ is the value of $\underline{c}$ on $(S)$ and ${c}$ is the restriction of $\underline{c}$ to $(\Gamma)$.

\subsection{$\!\!\!\!\!$}
\label{changeofvars}
Recall from  \ref{affdyn} that there is an affine Dynkin quiver $Q$ -- choose any orientation of the diagram --  corresponding to $\Gamma$. The $k$-th vertex corresponds to an irreducible representation $S_k$ whose dimension is $\delta_k$. Recall too the definition of $\blambda (t,c) \in \C^{Q^0}$ from \eqref{prpparam} which, together with setting $\nu = c_1 |\Gamma| /2$, produces a linear isomorphism between $\C\times \C[\mathcal{S}]^{\ad G}$ and $\C^{Q^0} \times \C$ such that $t = \blambda(t,c) \cdot \delta$.
\subsection{Generalised deformed preprojective algebras}
There is a quiver theoretic description of $H_{t,c}(G)$, at least up to Morita equivalence. This is easy to explain heuristically now; the explicit description will be given in the next section. 

We are trying to deform the algebra $\C [V]\rtimes G$ where $V = (\C^2)^n$. Rewriting this algebra as $(\C[\C^2]\rtimes \Gamma)^{\otimes n} \rtimes S_n$ suggests that we can construct a family of deformations simply by simultaneously deforming each of the tensorands $\C[\C^2]\rtimes \Gamma$. However, $\Gamma$ is a finite subgroup of $SL(2, \C)$ so we know by Theorem \ref{PBWthm} that these deformations correspond to the non-trivial conjugacy classes of $\Gamma$. Now by Theorem \ref{CBHthm} each of these deformations is Morita equivalent to a deformed preprojective algebra, so we have produced a family of deformations depending on the conjugacy classes of type $(\Gamma)$, namely $(\Pi^{\blambda(t,c)} (Q))^{\otimes n} \rtimes S_n.$

The deformations corresponding to the class of type $(S)$ arise now exactly as the Dunkl operators arose for $S_n$ acting on $\C[\C^{2n}] = \C[\C^2]^{\otimes n}$. In fact, since a generic choice of parameter $\blambda(t,c)$ produces a simple $\Pi^{\blambda(t,c)} (Q)$ by Theorem \ref{defppthm}(5), a homological theorem of \cite[Theorem 6.1]{etob}, inspired by \cite{AFLS}, assures us of the existence of such deformations.
\subsection{$\!\!\!\!\!$}
Here is the quiver-theoretic algebra defined in \cite[Definition 1.2.3]{ganginz}.
Let $k \defn \oplus_{i\in Q^0} \C e_i$ and $E$ be the vector space over $\C$ whose basis is given by the set of edges $\{a, a^*: a\in Q\}$. Set $\B \defn k^{\otimes n}$ and for any $1\leq \ell\leq n$, define the 
$\B$-bimodules
$$ \E_{\ell} \defn k^{\otimes (\ell-1)}\otimes E\otimes k^{\otimes (n-\ell)}
\qquad \mathrm{and}\qquad \E \defn
\bigoplus_{1\leq\ell\leq n} \E_\ell\,.$$
Given two elements $\ve\in\E_\ell$ and $\ve'\in\E_m$ of the form
\begin{equation} \label{eqnve}
\ve= e_{v_1}\otimes e_{v_2}\otot a\otot e_{h(b)} \otot e_{v_n}\,,
\end{equation}
\begin{equation} \label{eqnve'}
\ve'= e_{v_1}\otimes e_{v_2}\otot e_{t(a)}\otot b\otot e_{v_n}\,,
\end{equation}
where $\ell\neq m$, $a,b\in \QQ$ and $v_1, \ldots, v_n \in V$,
define
\begin{align}
\lf\ve,\ve'\rf &:=& 
(e_{v_1}\otot a\otot e_{h(b)} \otot e_{v_n})
(e_{v_1}\otot e_{t(a)}\otot b \otot e_{v_n}) \nonumber\\ & &
-(e_{v_1}\otot e_{h(a)}\otot b \otot e_{v_n})
(e_{v_1}\otot a\otot e_{t(b)} \otot e_{v_n}) . \nonumber
\end{align}
For $(\blambda , \nu) \in \C^{Q_0}\times \C$ define the algebra
$\mathsf{A}_{n,\blambda,\nu}$ to be the quotient of $T_{\B}\E\rtimes {S}_n$
by the following relations.
\begin{itemize}
\item[\vi] 
For any $v_1,\ldots, v_n\in V$ and $1\leq \ell\leq n$:
$$e_{v_1}\otot (r_{v_\ell}-\la_{v_\ell}e_{v_\ell}) 
\otot e_{v_n} = \nu \sum_{\{ j\neq\ell \,|\, v_j=v_\ell\}}
(e_{v_1}\otot e_{v_\ell} \otot e_{v_n})( j \, \ell) \,.$$
\item[\vii]
\label{quivercommutator}
For any $\ve, \ve'$ of the form \eqref{eqnve}--\eqref{eqnve'}:
\[  \lf\ve,\ve'\rf = \left\{ \begin{array}{ll}
\nu(e_{v_1}\otot e_{h(a)}\otot e_{t(a)}\otot e_{v_n})(\ell \, m)
& \textrm{if $a = b^*$ \text{ and } $b\in Q$}\,,\\
- \nu(e_{v_1}\otot e_{h(a)}\otot e_{t(a)}\otot e_{v_n})(\ell \, m)
& \textrm{if $b = a^*$ \text{ and } $a\in Q$}\,,\\
0 & \textrm{otherwise}\,.  \end{array} \right.  \]
\end{itemize}
As we mentioned, setting $\nu = 0$ produces $\Pi^{\blambda} (Q)^{\otimes n}\rtimes S_n$; taking $\Gamma = \{ 1\}$ recovers the relations for the rational Cherednik algebra associated to $S_n$ with $V = (\C^n)^2$.

Generalising Theorem \ref{CBHthm} we have the following theorem of Gan and Ginzburg.
\begin{theorem}[{\cite[Theorem 3.5.2]{ganginz}}] \label{GGmor}
There is a Morita equivalence between $H_{t,\underline{c}}(\Gamma_n)$ and $\mathsf{A}_{n,\blambda(t,c),\nu}$.
\end{theorem}
\subsection{$\!\!\!\!\!$} This Morita equivalence has proved very useful in the study of finite dimensional irreducible $H_{1,c}(G)$-representations. 

If $\nu = 0$ then $\mathsf{A}_{n,\blambda,\nu} = \Pi^{\blambda} (Q)^{\otimes n}\rtimes S_n$ and the finite dimensional irreducible representations can be described in the following standard manner. Take a partition $\mu = (\mu_1, \ldots , \mu_r)$ of $n$, set $S_{\mu} = S_{\mu_1}\times \cdots \times S_{\mu_r}$ be the corresponding Young subgroup, and let $W = W_1\otimes \cdots \otimes W_r$ where each $W_i$ is an irreducible representation of $S_{\mu_i}$. Now choose a collection $Y_1, \ldots , Y_r$ of irreducible, pairwise non-isomorphic representations of $\Pi^{\blambda}(Q)$ and form the irreducible representations $Y=Y_1^{\mu_1}\otimes \cdots \otimes Y_r^{\mu_r}$ of $ \Pi^{\blambda}(Q)^{\otimes n}$. Then $Y\otimes W \! \!\uparrow \defn Ind_{S_{\mu}}^{S_n} (W\otimes Y)$ is an irreducible $\Pi^{\blambda}(Q)^{\otimes n}\rtimes S_n = \mathsf{A}_{n,\blambda,0}$-representation. All irreducible $\mathsf{A}_{n,\la,0}$-representations arise this way. Since Theorem \ref{defppthm}(2) classifies all finite dimensional irreducible representations of $\Pi^{\blambda}(Q)$ with $\blambda\cdot \delta \neq 0$, the above describes all finite dimensional representations for $H_{t,\underline{c}}(G)$ with $c_1 = 0$ and $t\neq 0$.

In \cite[Theorem 3.1]{em} and \cite[Theorem 1.3]{mon2} Etingof and Montarani use deformation theoretic arguments to find sufficient conditions on the parameters to ensure that a given $\mathsf{A}_{n,\blambda,0}$-representation extends to an $\mathsf{A}_{n,\blambda',\nu}$-representation. Later, \cite[Section 6]{gan} Gan introduced reflection functors for the algebras $\mathsf{A}_{n,\blambda,\nu}$, generalising those of Bernstein-Gelfand-Ponomarev. He used them to prove the necessity of Etingof and Montarani's conditions.
\begin{theorem}[\cite{em}, \cite{mon2}, \cite{gan}]
The $\mathsf{A}_{n,\blambda_0,0}$-representation $W\otimes Y\!\!\uparrow$ extends to a $\mathsf{A}_{n, \blambda', \nu}$-representation for some $(\blambda ' ,\nu)$ with  $\nu \neq 0$ if and only if 
\begin{enumerate}
\item for each $1\leq i\leq r$, the irreducible $S_{\mu_i}$-representation has rectangular Young diagram of size $a_i\times b_i$;
\item  $\Ext^1_{\Pi^{\blambda}(Q)}(Y_i, Y_j) = 0$ for any $i\neq j$;
\item $\blambda \cdot \alpha_i = (a_i  - b_i)\nu$ for each $1\leq i \leq r$  where $\alpha_i$ is the dimension vector of $Y_i$.\end{enumerate}
The deformation exists in a linear subspace of the parameter space of codimension $r$ (which is explicitly described) and it is unique.
\end{theorem}

When combined with Theorem \ref{GGmor}, this is more-or-less the extent of our knowledge of the finite dimensional representations for $H_{t,\underline{c}}(G)$ when $t=1$. 

\subsection{The $t=0$ case} In order to motivate the geometric constructions we will use later in this chapter we study the case $t=0$. There is a precise sense in which the representation theory of $H_{0,\underline{c}}(G)$ is now well-understood, but for a reason we explain below, explicit questions are still difficult to answer. 

\subsection{$\!\!\!\!\!$} There is a description of $Y_{\underline{c}} = \Spec Z(H_{0,\underline{c}})$ for any $c$ which generalises the two-dimensional Theorem \ref{defppthm}(2) to higher dimensions. Given a quiver $Q$ associated to $\Gamma$, set $Q_{\sf CM}$ to be the quiver obtained from $Q$ by introducing one new vertex, $\infty$, and one new arrow which begins at $\infty$ and ends at the extending vertex. Given $c\in \C[\mathcal{S}]^{\ad G}$ define $\blambda'(t,c) \in \C^{Q^0_{\sf CM}}$ by \begin{equation} \label{thelambdas}\blambda'(t,c)_k =  \begin{cases} \frac{1}{\delta_k^2} \blambda(t,c)_k \qquad & \text{if $k\in Q^0$ is not extending} \\
\blambda(t,c)_k -\nu & \text{if $k$ is the extending vertex} \\
-n \sum_{k\in Q_0} \blambda'(t,c)_k \delta_k & \text{if $k = \infty$.} 
\end{cases}\end{equation}
Set $\epsilon_{\infty}\in \C^{Q_{CM}^0}$ to be the unit vector concentrated at the $\infty$ vertex. 
\begin{theorem}[{\cite[Theorem 1.13]{EG}, \cite[Theorem 6.4]{mo2}}] \label{identcent} Set $\alpha = n\delta + \epsilon_{\infty}$, a dimension vector for $Q_{\sf CM}$. 
There is an isomorphism $$Y_{\underline{c}} \stackrel{\sim}\longrightarrow \Rep (\Pi^{\blambda'(0,c)}(Q_{\sf CM}), \alpha) ) // GL(\alpha).$$ 
Moreover, there is a natural Poisson structure on the representation variety and the isomorphism identifies symplectic leaves. 
\end{theorem}
To explain the proof, recall that $Y_{\underline{c}}$ is equipped with a tautological $G$-equivariant coherent sheaf $H_{0,\underline{c}}e$ which, by Theorem \ref{azumayathm}, is a vector bundle over the smooth locus with fibres carrying the regular representation of $G$.  In essence the isomorphism is induced on the smooth locus by sending the tautological sheaf to its $\Gamma^{n-1}\rtimes S_{n-1}$-invariant sub-bundle, where $S_{n-1}< S_n$ is the subgroup acting on indices $\{2, \ldots , n \}$: $G, x_1 $ and $ y_1$ still act on this sub-bundle and produce the tautological bundle on the smooth locus of $\Rep (\Pi^{\blambda'(0,c)}(Q_{\sf CM}), \alpha)  // GL(\alpha)$.

In \cite[Theorem 1.3]{mo2} Martino goes on to describe the symplectic leaves of the representation variety as the (connected components of the) the semisimple representation type strata. Since Crawley-Boevey, \cite{CBdim}, and Le Bruyn, \cite{leb}, have described the irreducible $\Pi^{\blambda'(0,c)}(Q_{\sf CM})$-representations, this then provides an effective combinatorial description of the symplectic leaves of $Y_{\underline{c}}$, and in particular of the singular locus of $Y_{\underline{c}}$, which by Theorem \ref{azumayathm} has ramifications for the representation theory of $H_{0,\underline{c}}$.

\subsection{$\!\!\!\!\!$} Taking $c = 0$ in the above theorem we see that \begin{equation} \label{c0IM}V/G \stackrel{\sim} \longrightarrow \Rep (\Pi^{0}(Q_{\sf CM}), \alpha) // GL(\alpha).\end{equation} Now Nakajima's quiver varieties provide many resolutions of singularities of $\Rep (\Pi^{0}(Q_{\sf CM}), \alpha)$ and hence of $V/G$. We briefly recall their construction here: it is often referred to as {\it hamiltonian reduction}.

Consider the moment map $$\mu_{Q_{\sf CM}} : \Rep ({{\overline{Q}_{\sf CM}}} , \alpha) \longrightarrow \gl (\alpha), \qquad \{ (X_{a}, X_{a^*}) : a\in Q_{\sf CM}^1\} \mapsto \big(\mathop{\sum_{a\in Q_{\sf CM}^1}}_{h(a) = k} X_a X_{a^*} - \mathop{\sum_{a\in Q_{\sf CM}^1}}_{t(a)=k} X_{a^*}X_a\big)_k.$$ Its zero fibre $\mu_{Q_{\sf CM}}^{-1}(0)$ is reduced, by \cite[Theorem 1.1.2]{ganginz2}, so there will be no geometric ambiguity in the definition to follow.\ff{Without the vertex at infinity, however, it is a famous open question whether $\mu^{-1}_Q(0)$ is reduced, particularly when $Q$ is $\tilde{A}_1$.}
Given $\theta \in \mathbb{Z}^{Q^0}$, set ${\theta}_{\sf CM} = (\theta_{\infty}, \theta)\in \mathbb{Z}^{Q_{\sf CM}^0}$ with $\theta_{\infty} = -n\delta \cdot \theta$. There is a corresponding character, $\det^{\theta_{\sf CM}}$, of $\GL(\alpha)$ defined by $\det^{{\theta}_{\sf CM}} (g_{\infty}, (g_k)_k) = \det(g_{\infty})^{\theta_{\infty}}\prod_{k\in Q^0} \det (g_k)^{{\theta}_k}$. We set $$\mathcal{M}_{\theta}(\Gamma, n) \defn \mu_{Q_{\sf CM}}^{-1}(0)//_{{\theta}_{\sf CM}} GL(\alpha) \defn \proj \bigoplus_{m\geq 0} \C[\mu_{Q_{\sf CM}}^{-1}(0)]^{\det^{m\hat{\theta}}}.$$ Here $\C[\mu_{Q_{\sf CM}}^{-1}(0)]^{\det^{m{\theta}_{\sf CM}}} \defn \{ f\in \C[\mu_{Q_{\sf CM}}^{-1}(0)] : (g_{\infty}, (g_k)_k) \cdot f = \det^{{\theta}_{\sf CM}}(g_{\infty}, (g_k)_k) f \}$, the space of $m{\theta}_{\sf CM}$-semi-invariants functions in $\mu^{-1}_{Q_{\sf CM}}(0)$, and $\proj$ stands for the projective variety associated to an $\mathbb{N}$-graded algebra. Projection onto to the degree zero component provides a canonical projective morphism $$\pi_{\theta} : \mathcal{M}_{\theta}(\Gamma, n) \longrightarrow \Rep (\Pi^{0}(Q_{\sf CM}), \alpha)//GL(\alpha)\cong V/G.$$ For generic choices of $\theta$ this is a symplectic resolution of $V/G$, \cite[Theorem 2.8]{nak}. Nakajima generally applies the condition $\theta \in \mathbb{N}^{Q^0}$ in his study of quiver varieties;  other choices of $\theta$ produce other interesting resolutions such as $\hil^n (\widetilde{\C^2/\Gamma})$, \cite[Theorem 4.9]{kuz}; the variation of the stability condition $\theta$ in the context of symplectic reflection algebras is studied in \cite{gordquiv} in general. 

Recall that $\mathcal{M}_{\theta}(\Gamma ,n)$ also has a representation theoretic description as the algebro-geometric quotient of the open subvariety of $\mu^{-1}_{Q_{\sf CM}}(0)$ consisting of those representations which have no proper subrepresentations $W$ whose dimension vector satisfies the inequality $\underline{\dim} W \cdot {\theta}_{\sf CM} < 0$, \cite[3.ii]{nak2} and \cite{king}. 

\subsection{$\!\!\!\!\!$} This description of resolutions of singularities of $V/G$, combined with Theorem \ref{identcent}, fleshes out some of the philosophy behind Theorem \ref{weresolve} explained in  \ref{res=def}: the resolutions have a hyper-K\"ahler structure and rotating the given complex structure produces the (smooth) deformations.

\subsection{$\!\!\!\!\!$} \label{gordsmiththm}
Now it is not hard to generalise the quiver variety picture above to show that $Y_{\underline{c}}$ for {\it any} $c\in \C[\mathcal{S}]^{ad G}$ -- which by Proposition \ref{sympsing} is a symplectic singularity -- has a symplectic resolution. We have seen in Theorem \ref{smoothpartcentre} that if $Y_{\underline{c}}$ is smooth then there is a Morita equivalence between $\Lmod{H_{0,\underline{c}}}$ and $\Coh (Y_{\underline{c}})$. In the general case we have the following theorem generalising Theorem \ref{smoothpartcentre}.

\begin{theorem}[{\cite[Theorem 1.2]{gordsmi}}] Let $G = \Gamma^n\rtimes S_n$ and $c \in \C[\mathcal{S}]^{\ad G}$. There is a symplectic resolution $\pi_{\underline{c}}: X_{\underline{c}} \longrightarrow Y_{\underline{c}}$ such that there is an equivalence of triangulated categories $$D^b (\Lmod{H_{0,\underline{c}}}) \stackrel{\sim}\longrightarrow D^b(\Coh (X_{\underline{c}}))$$
between the bounded derived category of finitely generated $H_{0,\underline{c}}$-modules and the 
bounded derived category of coherent sheaves on $X_{\underline{c}}.$
\end{theorem} 
This theorem is quite straightforward modulo one difficult point: the case $c =0$ which produces an equivalence between $D^b(\C[V]\rtimes G)$ and $D^b(\Coh (X))$ where $X$ is a symplectic resolution of $V/G$. A theorem of Bezrukavnikov and Kaledin, \cite[Theorem 1.1]{BK}, produces such as equivalence, but it is rather implicit and we require something more explicit: for $\Gamma = \{1\}$ this is essentially a combination of the celebrated $n!$ theorem of Haiman, \cite{haiman}, and the celebrated derived McKay correspondence of Bridgeland--King--Reid, \cite{BKR}; the case of general $\Gamma$ also follows from these works, since the $n!$ theorem is local and $\hil^n (\widetilde{\C^2/\Gamma})$ is modelled locally on $\hil^n ({\C^2})$. Given this, the proof proceeds by a deformation argument.

\subsection{$\!\!\!\!\!$} The bad news is that this proof of theorem above does not produce an explicit symplectic resolution $X_{\underline{c}}$, and as we have already seen there may be several non-isomorphic ones. In particular, it is not known at the moment whether $X_{\underline{c}}$ is a resolution arising from a quiver variety. However, inspired by Bondal and Orlov, \cite[Conjecture 7.1]{HvdB} conjectures that the derived categories of coherent sheaves on different symplectic (more generally, crepant) resolutions of $Y_{\underline{c}}$ should all be equivalent. 

\subsection{$\!\!\!\!\!$}
The theorem does have one consequence which {\it is} independent of the choice of symplectic resolution. Since $X_{\underline{c}}$ is a symplectic resolution of the normal variety $Y_{\underline{c}}$ we have $\mathcal{O}(X_{\underline{c}}) \cong \mathcal{O}(Y_{\underline{c}}) = Z(H_{0,\underline{c}})$, and under this isomorphism complexes of modules whose cohomology are annihilated by some power of a maximal ideal $\mathfrak{m}_y\in Z(H_{0,\underline{c}})$ are identified with complexes of sheaves whose cohomology are annihilated by some power of $\mathfrak{m}_y\in \mathcal{O}(X_{\underline{c}})$. On passing to Grothendieck groups this produces an isomorphism $K(H_{0,\underline{c}}(y))\otimes_{\mathbb{Z}}\C \cong K(\pi_{\underline{c}}^{-1}(y))\otimes_{\mathbb{Z}}\C$. However the complexified $K$-group of the fibre $\pi_{\underline{c}}^{-1}(y)$ is independent of the choice of resolution thanks to \cite[Proposition 6.3.2]{DenefLoser} and \cite[Theorem 2.12]{kaledinsurvey}, and is isomorphic to $H^*(\pi_{\underline{c}}^{-1}(y), \C)$. Thus the number of simple $H_{0,\underline{c}}(y)$-modules, which equals the rank of $K(H_{0,\underline{c}}(y))$, is the dimension of the cohomology $H^*(\pi_{\underline{c}}^{-1}(y), \C)$ for {\it any} symplectic resolution $\pi_{\underline{c}} : X_{\underline{c}} \longrightarrow Y_{\underline{c}}$. Now for a symplectic resolution of $Y_{\underline{c}}$ by quiver varieties, such cohomology spaces were studied by Nakajima in his geometric construction of representations of quantum enveloping algebras, \cite{nak2}, and this relates the representation theory of $H_{0,\underline{c}}(y)$ to Lie theoretic combinatorics, cf. Theorem \ref{defppthm}. 

\subsection{Back to $t=1$: differential operators on representations of quivers} Theorems \ref{centrethm} and \ref{identcent} show that $U_{0,\underline{c}}$ can be described in terms of functions on a representation variety of the double of the quiver $Q_{\sf CM}$. To generalise this quiver-theoretic description to $U_{1,\underline{c}}$ requires us to pass from commutative algebra to noncommutative algebra. Now $\Rep({{\overline{Q}}_{\sf CM}}, \alpha)$ can be interpreted as the cotangent bundle $T^*\Rep(Q_{\sf CM}, \alpha)$, so this suggests that we should look at the canonical noncommutative deformation of $\C[T^*\Rep(Q_{\sf CM}, \alpha)] $, namely the algebra of differential operators on $\Rep(Q_{\sf CM},\alpha)$, which we denote by $\D_{Q_{\sf CM}} \defn \D (\Rep(Q_{\sf CM}, \alpha))$. This $\C$-algebra has generators the pairwise commuting matrix coefficients $(X_a)_{ij}$ and their pairwise commuting partial derivatives $(\partial_a)_{ij} $ for each $a\in Q_{\sf CM}^1$ and $1\leq i \leq \alpha_{h(a)}$, $1\leq j\leq \alpha_{t(a)}$, and they satisfy the relations $[(\partial_a)_{ij}, (X_b)_{kl} ] = \delta_{ab} \delta_{ik}\delta_{jl}$ for all $a,b\in Q_{\sf CM}^1$ and appropriate $ i,j$. 

\subsection{$\!\!\!\!\!$}
The moment map $\mu_{Q_{\sf CM}} :  \Rep ({{\overline{Q}}} , \alpha) \longrightarrow \gl (\alpha)$ corresponds to an algebra mapping $\mu_{Q_{\sf CM}}^* : \sym(\gl (\alpha)) \longrightarrow \C [ T^*\Rep(Q_{\sf CM}, \alpha)]$ which has a natural noncommutative analogue which we now describe. The group $GL(\alpha)$ acts algebraically on $\Rep(Q_{\sf CM}, \alpha)$; differentiating this gives an action of $\gl(\alpha)$ on $\Rep(Q_{\sf CM}, \alpha)$ by vector fields which can be calculated to be $$(E_{st})_i \in \gl(\alpha) \mapsto \mathop{\sum_{a, t(a) = i}}_{1\leq u \leq \alpha_{h(a)}} (X_a)_{ut} (\partial_a)_{us} -  \!\!\!\!\mathop{\sum_{a, h(a) = i}}_{1\leq v \leq \alpha_{t(a)} } (X_a)_{sv}(\partial_a)_{tv} \in \vect(\Rep(Q_{\sf CM},\alpha)).$$ This extends to an algebra homomorphism \begin{equation} \label{NCmom} \tau_{Q_{\sf CM}}: U(\gl(\alpha)) \longrightarrow \D_{Q_{\sf CM}} ,\end{equation}
which is a noncommutative deformation of $\mu_{Q_{\sf CM}}^*$ in the sense that $\gr \tau_{Q_{\sf CM}} = \mu_{Q_{\sf CM}}^*$ if we equip left hand side of \eqref{NCmom} with the usual filtration induced from $\deg (\gl(\alpha)) = 1$ and the right hand side with the differential operator filtration.

\subsection{$\!\!\!\!\!$} Since $GL(\alpha)$ acts on $\Rep(Q_{\sf CM}, \alpha)$, it also acts on $\D_{Q_{\sf CM}}$ and the action is locally finite. We now have all but one of the ingredients required to define an analogue of the isomorphism of Theorem \ref{identcent}. The last is the construction of a character $\chi_{\underline{c}}$ of $\gl(\alpha)$ which will twist differential operators in order to produce the deformations we seek. For $(X_k)_k \in \gl(\alpha)$ set $$\chi_{\underline{c}} (X_k) = \begin{cases} \left(\blambda(|\Gamma|^{-1}, c)_k + n(\delta_i - \sum\delta_{t(a) = i}\delta_{h(a)})\right)\Tr(X_k) \qquad &\text{if $k\in Q_0$ is not extending} \\
\left(\blambda(|\Gamma|^{-1}, c)_k -\nu - n\right)\Tr(X_k) & \text{if $k$ is the extending vertex} \\
n\left(\nu - 1\right)\Tr(X_k)  & \text{if $k = \infty$}. \end{cases}$$
The following theorem has several authors: \cite[Theorem 1.17]{EG}, \cite[Theorem 2.5]{obl}, \cite[Theorem 1.4]{gordiff}, \cite[Theorem 1.4.4]{EGGO}
\begin{theorem} \label{radial}
Continue with the notation from above and let $\D_{Q_{\sf CM}} (\tau_{Q_{\sf CM}}-\chi_{\underline{c}})$ denote the left ideal of $\D_{Q_{\sf CM}}$ generated by the elements $ \tau_{Q_{\sf CM}}(X)-\chi_{\underline{c}}(X)$ for $X\in \gl(\alpha)$. There is a filtered algebra isomorphism $$\Phi_{\underline{c}}: \left(\frac{\D_{Q_{\sf CM}}}{\D_{Q_{\sf CM}} (\tau_{Q_{\sf CM}}-\chi_{\underline{c}})}\right)^{GL(\alpha)} \stackrel{\sim}\longrightarrow U_{1,\underline{c}}$$ whose associated graded mapping is the isomorphism \begin{align*}\gr \Phi_{\underline{c}} : \left(\frac{\C[T^*\Rep(Q_{\sf CM},\alpha))]}{\C[T^*\Rep(Q_{\sf CM},\alpha))]\mu_{Q_{\sf CM}^*}(\gl(\alpha))}\right)^{GL(\alpha)}&\cong \gr \left( \frac{\D_{Q_{\sf CM}}}{\D_{Q_{\sf CM}} (\tau_{Q_{\sf CM}}-\chi_{\underline{c}})} \right)^{GL(\alpha)} \\ &\stackrel{\sim}\longrightarrow \gr U_{1,\underline{c}} \cong \C[V]^G\end{align*} of \eqref{c0IM}.
\end{theorem}
It may seems surprising at first sight that $\left({\D_{Q_{\sf CM}}}/{\D_{Q_{\sf CM}} (\tau_{Q_{\sf CM}}-\chi_{\underline{c}})}\right)^{GL(\alpha)}$ is even an algebra: after all, we have factored out a left ideal from the simple ring $\D_{Q_{\sf CM}}$ and then taken $GL(\alpha)$-invariants. However, if we let $\gl(\alpha)$ act on $\D_{Q_{\sf CM}}$ via $\ad x:  u \mapsto [u, (\tau_{Q_{\sf CM}} - \chi_{\underline{c}})(X)]$ for $X\in \gl(\alpha)$, then we can identify \begin{align} \label{endoident} \left(\frac{\D_{Q_{\sf CM}}}{\D_{Q_{\sf CM}} (\tau_{Q_{\sf CM}}-\chi_{\underline{c}})}\right)^{GL(\alpha)} &\stackrel{\sim} \longrightarrow \left(\frac{\D_{Q_{\sf CM}}}{\D_{Q_{\sf CM}} (\tau_{Q_{\sf CM}}-\chi_{\underline{c}})}\right)^{\ad \gl(\alpha)} \nonumber \\ &\stackrel{\sim} \longrightarrow \End_{\D_{Q_{\sf CM}}}\left(\frac{\D_{Q_{\sf CM}}}{\D_{Q_{\sf CM}} (\tau_{Q_{\sf CM}}-\chi_{\underline{c}})}\right)^{\sf op},\end{align} where the first isomorphism holds by definition and the second is given by right multiplication. The algebra structure is now evident.  We call $\left({\D_{Q_{\sf CM}}}/{\D_{Q_{\sf CM}} (\tau_{Q_{\sf CM}}-\chi_{\underline{c}})}\right)^{GL(\alpha)}$ the {\it quantum hamiltonian reduction} of $Q_{\sf CM}$ with respect to $GL(\alpha)$.

\subsection{$\!\!\!\!\!$}
The proof of Theorem \ref{radial} for $\Gamma = \{ 1\}$ was given in \cite{ganginz2}, following a twist on the classical route of radial parts in \cite{EG}. Here the quiver $Q_{\sf CM}$ is simply 
\begin{equation} \label{SnQ} \xymatrix{ 
{\circ} \ar@(ul,dl)  
&& \ar[ll]\circ} 
\end{equation}The scalars $(\lambda \Id_n,\lambda)$ in the base change group $GL(\alpha)= GL(n, \C)\times \C^*$ act trivially on $\Rep({Q}_{\sf CM},\alpha) = \Mat_n(\C) \times \C^n$ and so, without loss of generality, we reduce the base change group to $GL(n,\C)$. Now consider the open set $U_{\sf CM}$ of $\Rep({Q}_{\sf CM},\alpha)$ consisting of pairs $(X,v)$ of a regular semisimple matrix $X$ and a vector $v$ such that $\C[X]v = \C^n$. Mapping the matrix $X$ to its characteristic polynomial produces a $GL(\alpha)$-equivariant morphism $$\rho: U_{\sf CM} \longrightarrow \hr/S_n$$ which is a principal $GL(n,\C)$-bundle. Thus any $S_n$-invariant function $f$ on $\hr$ can be lifted uniquely to a function $\hat{f}^{\chi}$ on $U_{\sf CM}$ which is $GL(n,\C)$-semi-invariant for some character $\chi$ of $\GL(n,\C)$. (And in fact, on a small enough neighbourhood of the slice $\hr\times {\bf 1} \defn \{ (X, (1,\ldots , 1)) \in \hr \times \C^n \} \subset U_{\sf CM}$, we can lift the function to one which is semi-invariant for any character of $\gl(n)$.) This allows us to construct a homomorphism $$\Phi_{\chi} : \D(\Rep(Q_{\sf CM}, \alpha))^{GL(\alpha)}\longrightarrow \D(\hr)^{S_n}$$ by defining $\Phi_{\chi}(D)(f) = D (\hat{f}^{\chi})\vert_{\hr\times {\bf 1}}$. Now recall that $U_{1,\underline{c}} \subset 
\D(\hr)^{S_n}$ thanks to Remark \ref{dunklexp}. Taking $\chi = \chi_{\underline{c}}$, one shows that the image of $\Phi_{\chi_{\underline{c}}}$ (after conjugating by an appropriate element in $\C[\hr]$) is $U_{1,\underline{c}}$ and that the kernel of $\Phi_{\chi_{\underline{c}}}$ contains $\left(\D_{Q_{\sf CM}} (\tau_{Q_{\sf CM}}-\chi_{\underline{c}})\right)^{GL(\alpha)}$. To check that this is actually all of the kernel one passes to $\gr \Phi_{\chi_{\underline{c}}}$ and uses the isomorphism \eqref{c0IM}.

For $\Gamma$ of type $A$ the theorem is proved by a generalisation of the argument above. However, the case of $\Gamma$ of type $D$ and $E$ is more complicated since there is no Dunkl embedding. It involves a careful study of several different mappings and differential operators on products of $\mathbb{P}^1$'s.  The details are in \cite{EGGO}. 

\subsection{The functor of hamiltonian reduction} 
Let $\Lmod{(\D_{Q_{\sf CM}},\gl(\alpha)_{\underline{c}})}$ denote the category of finitely generated $\D_{Q_{\sf CM}}$-modules which are locally finite as $(\tau_{Q_{\sf CM}}-\chi_{\underline{c}})(\gl(\alpha))$-modules. Given a module $M$ in this category, we can consider the set of fixed points $M^{\gl(\alpha)_{\underline{c}}} \defn \{ m \in M : \tau_{Q_{\sf CM}}(X) m = \chi_{\underline{c}}(X)m \text{ for all }X\in \gl(\alpha) \}.$ Clearly, this a $(\D_{Q_{\sf CM}}/\D_{Q_{\sf CM}} (\tau_{Q_{\sf CM}}-\chi_{\underline{c}}))^{GL(\alpha)}$-module and so by Theorem \ref{radial} we have constructed the {\it functor of hamiltonian reduction} 
$$\mathbb{H}_{\underline{c}} :\Lmod{(\D_{Q_{\sf CM}},\gl(\alpha)_{\underline{c}})} \longrightarrow \Lmod{U_{1,\underline{c}}}, \qquad M \mapsto M^{\gl(\alpha)_{\underline{c}}}.$$
Since $\gl(\alpha)$ is reductive and the objects of  $\Lmod{(\D_{Q_{\sf CM}},\gl(\alpha)_{\underline{c}})}$ are locally finite, the functor $\mathbb{H}_{\underline{c}}$ is exact. Tensoring over $U_{1,\underline{c}}$ by the $(\D_{Q_{\sf CM}}, (\D_{Q_{\sf CM}}/\D_{Q_{\sf CM}} (\tau_{Q_{\sf CM}}-\chi_{\underline{c}}))^{GL(\alpha)})$-bimodule $\D_{Q_{\sf CM}}/\D_{Q_{\sf CM}} (\tau_{Q_{\sf CM}}-\chi_{\underline{c}})$ is a left adjoint to $\mathbb{H}_{\underline{c}}$. The full subcategory $\ker \mathbb{H}_{\underline{c}} = \{ M \in (\Lmod{\D(\Rep(Q_{\sf CM},\alpha),\gl(\alpha)_{\underline{c}})}: \mathbb{H}_{\underline{c}}(M) = 0 \}$ is a Serre subcategory of $\Lmod{(\D_{Q_{\sf CM}},\gl(\alpha)_{\underline{c}})}$ and so we can form a quotient category. We have an equivalence $$\mathbb{H}_{\underline{c}} : \frac{\Lmod{(\D_{Q_{\sf CM}},\gl(\alpha)_{\underline{c}})}}{\ker \mathbb{H}} \stackrel{\sim} \longrightarrow \Lmod{U_{1,\underline{c}}}.$$ Thus, the representations of $U_{1,\underline{c}}$ can be studied in terms of a certain categories of $\D$-modules on $\Rep(Q_{\sf CM}, \alpha)$. This has been studied for $\Gamma = \{1\}$ first by Gan and Ginzburg, \cite{ganginz2}, and then by Finkelberg and Ginzburg, \cite{FGchar}, leading to a relationship with (a generalisation of) Lusztig's character sheaves.
\section{Very specific case: the symmetric group and Hilbert schemes}
The symplectic reflection groups $G = \mu_{\ell}^n\rtimes S_n  = G(\ell ,1 , n)$ acting on $(\C^n)^2= (\C^2)^n$ live in both families  of  \ref{groups}: they are complex reflection groups in their action on $\h = \C^n$; they are wreath products of kleinian subgroups of type $A$. Therefore we can apply the theory from Sections \ref{t1RCA} and \ref{t0RCA}, {\it and} from Section \ref{QUIVERS}. 

For $t = 0$ much of their representation theory is understood thanks to Theorem \ref{gordsmiththm}, and indeed in this case Conjecture \ref{gormartconj}(1) is confirmed.  

We will therefore concentrate on the case $t = 1$ for the rest of this section, and even specialise to the case $\ell = 1$. Thus $G = S_n$ and $V= \C^{2n}$.\ff{In fact we shall be slightly cavalier here: $\C^n$ is not the reflection representation of $S_n$ but rather $\h \oplus \triv$. The extra copy of $\triv$ results in a copy of $\D(\C)$ generated by $x_1+\cdots + x_n$ and $y_1 + \cdots + y_n$, and  an isomorphism $H_{1,c}(S_n, \C^{2n}) \cong \D(\C) \otimes H_{1,c}(S_n, \h\oplus \h^*)$. Since $\D(\C)$ is a particular {\it simple} algebra, we will elide the difference between $\C^n$ and $\h$ in what follows.} We will see the connection between finitely generated representations of $H_{1,c}(S_n)$ and coherent sheaves on the Hilbert scheme $\hil^n(\C^2)$. Many interesting and new phenomena arise for $\ell >1$, but the picture is not well understood at the moment, either algebraically or geometrically -- see \cite[Chapter 8]{haimansurvey}, \cite{gordquiv} and \cite{kuwabara}.

\subsection{Hilbert schemes}
We begin by motivating the appearance of Hilbert schemes. The celebrated theorem of Beilinson--Bernstein, \cite{BB}, shows that for a simple complex Lie algebra $\g$ and a regular integral dominant weight $\lambda+\rho$ there exists an equivalence of categories between $\Lmod{U(\g)_{\lambda}}$ and $\Lmod{\D^{\lambda}_{\mathcal{B}}}$, where $U(\g)_\lambda$ denotes the quotient of the enveloping algebra $U(\g)$ by the maximal ideal of its centre corresponding to $\lambda$, and $\D^{\lambda}_{\mathcal{B}}$ denotes the sheaf $\lambda$-twisted differential operators on the flag variety $\mathcal{B}$. There are filtrations on both $U(\g)_{\lambda}$ and $\D^{\lambda}_{\mathcal{B}}$ such that $\gr U(\g)_{\lambda} \cong \C[\mathcal{N}]$ where $\mathcal{N}$ is the set of nilpotent elements of $\g$ and $\gr \D^{\lambda}_{\mathcal{B}} \cong p_*\mathcal{O}_{T^*\mathcal{B}}$ where $p:T^*\mathcal{B} \longrightarrow \mathcal{B}$ is projection onto the base. Thus the Beilinson--Bernstein theorem can be considered as a quantisation of the Springer resolution $\pi: T^*\mathcal{B} \longrightarrow \mathcal{N}$, a symplectic resolution of singularities. Moreover the Beilinson--Bernstein theorem states that in the noncommutative setting the resolution $\pi$ has become an equivalence: this is not outrageous since the algebra $U(\g)_{\lambda}$ already has many very good homological properties, including finite global dimension, and so can already be considered to be a ``smooth" noncommutative algebra. 

\subsection{$\!\!\!\!\!$} In our setting the spherical sublagebras $U_{1,c}$ are analogues of the quotients $U(\g)_{\lambda}$, cf. \ref{firsttoyex}. Thanks to Proposition \ref{easyfilter} and \ref{simplicity1} the algebra $U_{1,c}$ has finite global dimension for most values of $c$. On the other hand $\gr U_{1,c} = \C [V]^{S_n}$. By definition, the right hand side is the coordinate ring of the $n$th symmetric product of the place $S^n\C^2$ consisting of $n$-unordered points on the plane. This is a symplectic singularity, but it has a resolution by the {\it Hilbert scheme of $n$ points of the plane}, $\Hilb$, which can be described as a quiver variety associated to the quiver $Q_{\sf CM}$ of \eqref{SnQ}, see for example \cite{nakbook}. Setting $$\mu_{Q_{\sf CM}}^{-1}(0)^{\sf ss} =  \{ (X,Y,i,j) \in \Mat_n(\C)^{\times 2}\times \C^n \times (\C^n)^* : [X,Y] + ij = 0, \C\langle X,Y\rangle i(1) =\C^n\}$$ and letting $GL(n,\C)$ act by $g\cdot (X,Y,i,j) = (gXg^{-1}, gYg^{-1}, gi, jg^{-1})$, the orbit mapping is a $GL(n, \C)$-principal bundle with image $\Hilb$ \begin{equation} \label{principal} p : \mu_{Q_{\sf CM}}^{-1}(0)^{\sf ss} \mapsto \mu_{Q_{\sf CM}}^{-1}(0)^{\sf ss}/GL(n, \C)\defn\Hilb.\end{equation}
 It is not hard to see that $j=0$ in $\mu^{-1}_{Q_{\sf CM}}(0)^{\sf ss}$ and so $X$ and $Y$ commute; thus the quiver data above encodes how $\C^n$ becomes a cyclic $\C[X,Y]$-module. Taking the annihilator in $\C[X,Y]$ of this representation produces the more familiar decsription of the Hilbert scheme $\Hilb = \{ I\triangleleft \C[X,Y] : \dim (\C [X,Y]/I) = n \}.$ In this description the symplectic resolution is easily described $$\pi : \Hilb \longrightarrow S^n\C^2, \qquad I \mapsto {\bsupp}(I),$$ where ${\bsupp}$ indicates support counted with multiplicity. 

\subsection{$\!\!\!\!\!$} \label{goodquestion}
By analogy with \cite{BB} it would make good sense to ask
\begin{center} {\it Is there a noncommutative deformation, $W$, of $\Hilb$ \\ such that $\Lmod{U_{1,c}}$ and $\Lmod{W}$ are Morita equivalent?}\end{center}

A positive answer to this would allow us to relate the representation theory of the symplectic reflection algebra and the algebraic geometry of the Hilbert scheme. So we would like to interpret $eL_c(\lambda)$ and $e\Delta_c(\lambda)$ geometrically, as well as other obvious $U_{1,c}$-modules such as $U_{1,c}$ and $eH_{1,c}$ as well as $eH_{1,c-1}\delta e$, the bimodule which induces the shift functor ${\sf S}_{c-1}$. On the other hand there is a rich seam of geometric combinatorics around the Hilbert scheme which we would like to understand representation theoretically. In particular there is the {\it Procesi bundle} $\mathcal{P}$, an $S_n$-equivariant vector bundle of rank $n!$ on $\Hilb$ whose fibres each carry the regular representation of $S_n$, \cite{haiman} (the existence of this bundle is equivalent to the $n!$ theorem describing $(q,t)$-Macdonald--Kostka polynomials), the tautological bundle $\mathcal{T}$ whose fibre above an ideal $I\in \Hilb$ is the $n$-dimensional space $\C[X,Y]/I$, as well as the distinguished ample line bundle $\mathcal{L} = \wedge^{\sf top} \mathcal{T}$. 

However, we would really also like that $\Lmod{W}$, if it exists, is a category that can be attacked with new tools. For instance, the Beilinson--Bernstein theorem allows us to translate problems about representations of $\g$ to problems on $\D$-modules on $\mathcal{B}$. But there are many tools to study $\D$-modules, particularly the Riemann--Hilbert correspondence. This transfers questions to the topology of sheaves on $\mathcal{B}$, and in particular leads to a description of the decomposition matrix for category $\mathcal{O}$ of $\g$, proving of the Kazhdan--Lusztig conjecture. 

\subsection{Noncommutative algebraic geometry} \label{NAG} The first answer to Question \ref{goodquestion} involves noncommutative algebraic geometry. One point of view on noncommutative algebraic geometry begins with the equation $$X = \Coh (X).$$ On the left hand side we have an algebraic variety $X$, on the right hand side its category of coherent sheaves. The equality means that either side always determines the other: in particular $X$ can be recovered from $\Coh (X)$. Thus the study of $X$ can be just as well thought of as the study of the category $\Coh (X)$. This point of view is very useful: we can say a noncommutative variety {\it is} simply a category $\mathbb{X}$ which has the same signal properties as $\Coh (X)$. Now such categories occur as various types of categories of modules for algebras and each such category can be studied with some geometric intuition. This philosophy has been pioneered very successfully; a recent survey is \cite{staffordvandenbergh}.

\subsection{$\!\!\!\!\!$} Following this we seek a category $\mathbb{X}_c$ which is Morita equivalent to $\Lmod{U_{1,c}}$ and which is also a ``deformation" of $\Coh (\Hilb)$. Here the deformation means that object $M$ of $\mathbb{X}_c$ can be given filtrations such that the associated graded object $\gr M$ is a coherent sheaf on $\Hilb$. Moreover, these filtrations should be compatible with filtrations on objects of $\Lmod{U_{1,c}}$ and the square $$\begin{CD} \mathbb{X}_c @> \gr >> \Coh (\Hilb) \\ @V \wr VV @VV \pi_{\ast} V \\ \Lmod{U_{1,c}} @> \gr >> \Coh (S^n\C^2) \end{CD}$$ should commute in the  appropriate sense.

\begin{theorem}[{\cite[Theorem 6.4]{gordst1}}]  \label{Zalg} Suppose that  $c\nless 0$. Then there exists a category $\mathbb{X}_c$ of coherent sheaves on a noncommutative variety that completes the above diagram.

\end{theorem}
The proof of this theorem is rather technical, but easy to understand in spirit. To deform $\Coh (\Hilb)$ we are going to deform a homogeneous coordinate ring of $\Hilb$. So take the ample line bundle $\mathcal{L}$ described in \ref{goodquestion} and the corresponding $\mathbb{N}$-graded algebra $R = \oplus_{i\geq 0} H^0(\Hilb , \mathcal{L}^{\otimes i})$. This algebra gives us yet another description of $\Hilb$, this time as $\proj R$. There are two explicit descriptions of the global sections of the tensor powers of $\mathcal{L}$: \begin{equation} \label{HaimansIMs} H^0(\Hilb, \mathcal{L}^{\otimes i}) \cong \C[\mu^{-1}(0]^{\det^i}  \cong  {\sf S}^i\end{equation} where ${\sf S} = \C[V]^{\sf sign}$ and ${\sf S}^i$ is the product of $i$ copies of ${\sf S}$ taken inside the polynomial algebra $\C[V]$. 

Following \ref{radial} it is clear how to deform $R$: replace $\C[\mu^{-1}(0)]^{\det^i}$ by $\left(\D_{Q_{\sf CM}}/(\D_{Q_{\sf CM}}(\tau_{Q_{\sf CM}} - \chi_c) \right)^{\det^i}$. However, if we take the direct sum of all of these to produce an analogue of $R$ we run into a technical problem -- the sum does not have a well-defined multiplication since $\left(\D_{Q_{\sf CM}}/(\D_{Q_{\sf CM}}(\tau_{Q_{\sf CM}} - \chi_c) \right)^{\det^i}$ is a $(U_{1,c+i}, U_{1,c})$-bimodule. This is solved by the formalism of {\it $\mathbb{Z}$-algebras} which replaces $R$ by the deformation $$\hat{R}_c = \bigoplus_{i\geq j \geq 0} \left(\frac{\D_{Q_{\sf CM}}}{(\D_{Q_{\sf CM}}(\tau_{Q_{\sf CM}} - \chi_{c+j})} \right)^{\det^{i-j}}.$$ With this in hand, we recall a fundamental theorem of Serre which states that the coherent sheaves on $\proj R$ can be described algebraically as the category of finitely generated graded $R$-modules, $\Lgmod{R}$, factored by the full subcategory of all torsion graded $R$-modules, $\tors{R}$. Mimicking this, we take $\mathbb{X}_c$ to be $\Lgmod{\hat{R}_c} / \tors{\hat{R}_c}$.

To see why $\mathbb{X}_c$ is equivalent to $\Lmod{U_{1,c}}$ we need a noncommutative version of the second isomorphism of \eqref{HaimansIMs}. Since $\delta\in \C[\h]$ is a $\sign$ semi-invariant element, an obvious noncommutative analogue of ${\sf S}$ is the $(U_{1,c+1},U_{1,c})$-bimodule $eH_{1,c+1}\delta e$ from \eqref{sphershift}; it follows that the product $(eH_{1,c+i}\delta e)\cdots (eH_{1,c+2}\delta e)(eH_{1,c+1}\delta e)$ should play the role of ${\sf S}^i$. It is shown in \cite[Main Theorem]{GGS} by a variation on Theorem \ref{radial} that this product is isomorphic to $\left(\D_{Q_{\sf CM}}/(\D_{Q_{\sf CM}}(\tau_{Q_{\sf CM}} - \chi_c)\right)^{\det^i}$. Thus $\mathbb{X}_c$ can be built out of the bimodules that induce the composition of shift functors ${\sf S}^{c+i}_{c+j}\defn S_{c+i-1}\circ \cdots \circ S_{c+j} : \Lmod{U_{1,{c+j}}} \longrightarrow \Lmod{U_{1,c+i}}$. By  \ref{shiftequiv} each ${\sf S}^{c+i}_{c+j}$ is a Morita equivalence if $c\nless 0$; it then follows that the functor $$\Phi_c : \Lmod{U_{1,c}} \longrightarrow \mathbb{X}_c, \qquad M \mapsto \bigoplus_{i\geq 0} S_c^{c+i} (M)$$ induces the required equivalence.

\subsection{$\!\!\!\!\!$} The explicit construction of $\Phi_c$ in \ref{Zalg} shows that a filtered module in $\Lmod{U_{1,c}}$ produces a filtered object of $\mathbb{X}_c$; taking the associated graded object of this then produces a coherent sheaf on $\Hilb$. Not every finitely generated $U_{1,c}$-module has a natural filtration on it, but for some that do the corresponding coherent sheaf on $\Hilb$ can be calculated, \cite{gordst2}:
\begin{itemize}
\item If $M = U_{1,c}$ with its differential operators filtration then the corresponding sheaf is $\mathcal{O}_{\Hilb}$.
\item If $M = eL_{m+1/h}(\triv)$ with the filtration of \ref{coinvariants} then the corresponding sheaf is the restriction $\mathcal{L}^{\otimes m}\vert_Z$, the restriction of $\mathcal{L}^{\otimes m}$ to the punctual Hilbert scheme $Z = \pi^{-1}(0)\subset \Hilb$.
\item If $M = e\Delta_c(\lambda)$ for $\lambda \in \irr(S_n)$ with the filtration induced from the filtration by ${\bf h}$-degree, then the corresponding coherent sheaf is $(\mathcal{P}/\mathfrak{h}\mathcal{P})^{\lambda}$, the $\lambda$ isotypic component of the quotient of $\mathcal{P}$ by $\C[\h^*]_+\mathcal{P}$.
\item If $M = eH_{1,c}$ with its filtration inherited from the differential operators filtration on $H_{1,c}$, then the corresponding coherent sheaf is $\mathcal{P}$, the Procesi bundle on $\Hilb$.
\item if $M = eH_{1,c+1}\delta e$ with the filtration induced from the differential operator filtration on $\D(\hr)\rtimes G$, then the corresponding sheaf is $\mathcal{L}$.
\end{itemize}
The first and last items are tautologies. The second is related to Haiman's work on $\co{\h\oplus \h^*}{S_n}$ and higher Catalan numbers and provides the correct context for \ref{coinvariants}. The third item follows from the fourth, and it is the fourth that is the most ``expensive" since it uses some key ingredients in the proof of the $n!$ theorem, \cite{haiman}. Informally, it is very satisfying since by \ref{simplicity1} the  $U_{1,c}$-module $eH_{1,c}$ is projective of rank $n!$ for almost all values of $c$; the identification of the associated graded of $\Phi_c(eH_{1,c})$ with $\mathcal{P}$ is, however, quite challenging and using $eH_{1,c}$ to give a new proof of the $n!$ theorem would seem to need more understanding than we have at the moment. 

\subsection{$\!\!\!\!\!$} Choosing different filtrations on the {\it same} object $M\in \Lmod{U_{1,c}}$ can produce {\it different} coherent sheaves on $\Hilb$. However, taking the support counted with multiplicity of these sheaves is independent of the choice of filtration and so to each $M\in \Lmod{U_{1,c}}$ we can associate $\mathbf{Ch}(M)$, a cycle in $\Hilb$. This is particularly interesting for objects $M = e\hat{M}$ where $\hat{M} \in \mathcal{O}_c$ since it induces an isomorphism between $K(\mathcal{O}_c)\otimes_{\Z} \C$ and $H_{\sf top}^{BM}(\pi^{-1}(\C^n\times \{0\}), \C)$ where $H_{\sf top}^{BM}$ is the top Borel--Moore homology group, \cite[Corollary 6.10]{gordst2}. Nakajima, \cite[Chapter 9]{nakbook}, and Grojnowski, \cite{groj}, identified the direct sum over $n\geq 0$ of these Borel--Moore spaces with the ring of symmetric functions. Under these identifications the $[L_c(\lambda)]$'s are identified with the specialisation at $v=1$ of the $U_v(\widehat{sl}_{d})$ lower canonical basis on symmetric functions introduced by Leclerc and Thibon, \cite{LT}, where $d$ is the order of $c$ in the abelian group $\C/\mathbb{Z}$, \cite[Proposition 6.11]{gordst2}.

\subsection{Microlocalisation} In recent work Kashiwara and Rouquier  have given another answer to Question \ref{goodquestion} by combining the radial parts approach of Gan and Ginzburg and the Morita equivalences of \eqref{shiftequiv} with the theory of microlocalisation to quantise the sheaf of regular functions on $\Hilb$ and find an analogue of the Beilinson--Bernstein theorem for symplectic reflection algebras associated to $S_n$. 

\subsection{$\!\!\!\!\!$} It is rather easy to explain the technical difficulty faced in quantising $\Hilb$ and how it is overcome in \cite{KR}.  The sheaf of functions on the cotangent bundle $T^*\C^n$ has a deformation $\mathcal{W}_{T^*\C^n}$ over the ring of formal Laurent series $\C[\hbar^{-1}, \hbar]]$ provided by the Moyal product. For $u,v \in \mathcal{O}_{T^*\C^n}$ this means that we have a deformed multiplication given by
$$u \star v = \sum_{\underline{\alpha}\in \mathbb{Z}_{\geq 0}^n} \hbar^{|\underline{\alpha}|} \frac{1}{\underline{\alpha}!} \partial^{\alpha_1}_{y_1} \cdots \partial^{\alpha_n}_{y_n} (u) \partial^{\alpha_1}_{x_1}\cdots \partial^{\alpha_n}_{x_n} (v).$$ It is easy to check that the subalgebra generated by the $x_i$ and $\hbar^{-1}y_i$ for $1\leq i \leq n$ is isomorphic to the ring of differential operators $\D(\C^n)$. A crucial difference, however, between $\mathcal{W}_{T^*\C^n}$ and $\D (\C^n)$ is that the former is a sheaf on (i.e. localises over) $T^*\C^n$ since the terms defining each degree in the Moyal product are bidifferential operators in $u$ and $v$ and hence localisable, whilst the latter is a sheaf only on the base $\C^n$.

Now if $Z$ is a smooth variety, then the cotangent bundle $T^*Z$ has a deformation $\D_Z$, and this can be thought of as glueing together the $\D(\C^n)$'s for different patches of $Z$. If $Z = \mathcal{B}$, the flag manifold of a simple complex Lie group, this deformation is the sheaf of algebras which appears in the Beilinson-Bernstein theorem. Note that it does not depend on the formal parameter $\hbar$. 

If we want to deal with a general complex algebraic symplectic variety $X$ then Darboux's Theorem tells us that $X$ is (holomorphically) locally of the form $T^*\C^n$, but not necessarily globally. Thus to deform we must glue on $X$ and hence we must use $\mathcal{W}_{T^*\C^n}$ as our local model. In this way we will form the {\it $W$-algebra} $\mathcal{W}_X$. This, however, introduces two problems: the various $\mathcal{W}_{T^*\C^n}$ do not necessarily glue together to form a sheaf of algebras, but rather a sheaf of ``algebroid stacks" due to the data of automorphisms, see \cite{schapira}; a formal Laurent parameter $\hbar$ has appeared which will have no analogue in the world of $U_{1,c}$, and hence threaten our goal of a Beilinson--Bernstein type theorem. The first problem is not too bad since there is still a respectable category of modules for algebroid stacks; and, in any case, for varieties with enough cohomology vanishing one might even expect the stackiness to disappear. The second problem is serious and we have to make a restriction on the the types of symplectic varieties we consider in order to overcome it: we insist that $X$ has a {\it contracting $\C^*$-action}. This means that $\C^*$ acts on $X$ such that the induced action on forms we have $\lambda \cdot \omega_X = \lambda^m \omega_X$ for some $m\in \mathbb{Z}_{>0}$ and for all $\lambda \in \C^*$. If we declare $\hbar$ to be a $\C^*$-eigenvector with weight $m$, then $\mathcal{W}_X$ inherits a $\C^*$-action and in a procedure analagous to \ref{gradeddef} we can remove the formal parameter $\hbar$ by studying appropriate categories of $\C^*$-equivariant $\mathcal{W}_X$-modules. In the case $X = T^*Z$ this category will be equivalent to $\D_Z$; it will, however, exist for algebraic symplectic varieties of more general type. 

\subsection{$\!\!\!\!\!$} The discussion above is relevant to the study of symplectic reflection algebras because $\Hilb$ is {\it not} the cotangent bundle of any variety (for $n>2$; for $n=2$ it is $T^(\mathbb{P}^1\times \C)$) -- if it were it would basically be  the cotangent bundle of the punctual Hilbert scheme $\pi^{-1}(0)$, but this is singular for $n>2$ -- and so to quantise we would like to use $\mathcal{W}_{\Hilb}$. The scalar $\C^*$-action on $\C^2$ induces a $\C^*$-action on $\Hilb$ such that $\lambda \cdot \omega_{\Hilb} = \lambda^2 \omega_{\Hilb}$ (in the quiver-theoretic description this action is $\lambda \cdot (X, Y, i, j) = (\lambda X, \lambda Y, \lambda i, \lambda j)$). 

\begin{theorem}[{\cite[Theorem 4.9]{KR}}] 
There is a $\mathcal{W}$-algebra $\mathcal{W}_c$ on $\Hilb$ with contracting $\C^*$-action such that the following functors define quasi-inverse equivalences between $\Coh^{gd}_{\C^*} (\mathcal{W}_c)$ and $\Lmod{U_{1,c}(S_n)}$: \begin{align*} {\sf{global \,  sections}} : & \,\, \mathcal{M} \mapsto \Hom_{\Coh^{gd}_{\C^*} (\mathcal{W}_c)}( \mathcal{W}_c, \mathcal{M}) \\ {\sf localisation} : & \,\, \mathcal{W}_c\otimes_{U_{1,c}} M \mapsfrom M. \end{align*}
\end{theorem}

Here $\Coh^{gd}_{\C^*} (\mathcal{W}_c)$ denotes a category of $\C^*$-equivariant $\mathcal{W}_c$-modules with a standard {\it good} generation property, explained in \cite[Section 2.3]{KR}. (There is a restriction on $c$: see \cite[Theorem 4.9]{KR} for details.)

The algebra $\mathcal{W}_c$ appearing in the theorem is a sheaf of algebras on $\Hilb$, described as follows. Recall that $\D_{Q_{\sf CM}}$ is the ring of differential operators on $\Rep(Q_{\sf CM}, \alpha)$. Let $\mathcal{W}_{Q_{\sf CM}}$ denote the $\mathcal{W}$-algebra on $T^* \Rep(Q_{\sf CM}, \alpha)$, as defined above. Now, following the Gan and Ginzburg recipe of \eqref{endoident}, we set : $$\mathcal{W}_c \defn \C[\hbar^{1/2}] \otimes_{\C[\hbar]} p_{\ast}\left( \End_{\mathcal{W}_{Q_{\sf CM}}} \left(\frac{\mathcal{W}_{Q_{\sf CM}}}{\mathcal{W}_{Q_{\sf CM}}(\tau_{Q_{\sf CM}} - \chi_c)}\right)\right)^{\GL(\alpha), \sf op}$$ where $p$ is the mapping \eqref{principal} which sends elements of $\mu_{Q_{\sf CM}}^{-1}(0)$ to their $GL(n,\C)$ orbits in $\mu_{Q_{\sf CM}}^{-1}(0)/GL(n,\C)  = \Hilb$. (The change of base ring from $\C[\hbar]$ to $\C[\hbar^{1/2}]$ is a technicality to allow the $\C^*$-action on $\omega_{\Hilb}$ have weight $1$.) The essential points of the proof of the theorem then consist in identifying the noncommutative analogue of ample bundle $\mathcal{L}$ in terms of the $\mathcal{W}$-algebra, relating it to the Heckman-Opdam shift functors ${\sf S}_c$, and using \ref{shiftequiv} to prove that $\Hilb$ is $\mathcal{W}$-affine, which means that the global sections and localisation functors are quasi-inverse to one another. 

\subsection{\!\!\!\!\!} What is required now is a good set of tools to probe $\Lmod{\mathcal{W}_c}$ and with it a geometric or topological description of important objects in $\Lmod{U_{1,c}}$. In particular it would be very interesting to understand the sheaf that corresponds to $eH_{1,c}$: it should be a quantisation of the Procesi bundle $\mathcal{P}$.

\section{Problems}
Here are problems for next year. (They overlap with the problems in \cite{rou}, \cite{openproblems} and a few of the original set of problems in \cite[Section 17]{EG} -- the other questions in \cite{EG} have been solved.)

\smallskip
\noindent 
{\bf 1.} (Chapter 2) Study symplectic reflection algebras $H_{t,c}(G)$ when $G$ does not belong to Family (1) or (2) of \ref{groups}.


\smallskip
\noindent
{\bf 2.} (Chapters 4 and 5) Determine when $U_{1,c}$ and $H_{1,c}$ are Morita equivalent; determine general conditions which ensure that ${\sf S}_c: \Lmod{H_{1,c}} \longrightarrow \Lmod{H_{1,c+1}}$ (or its generalisations of \ref{generalshift}) is an equivalence.

\smallskip
\noindent
{\bf 3.} (Chapter 4)
Show that for any nonzero finitely generated $H_{1,c}$-module $M$ the generalised Bernstein inequality holds, $
{\sf GKdim}(M ) \geq \frac{1}{2} {\sf GKdim}(H_{1,c}/{\sf Ann}(M)).$ 

\smallskip
\noindent
{\bf 4.} (Chapter 4) Study the function $d: \C[\mathcal{S}]^{\ad G} \longrightarrow \N$ defined by $d(c) = \min  \{ \text{\sf GKdim} (I): I \text{ irreducible $H_{1,c}$-representation}\}$.

\smallskip
\noindent
{\bf 5.}
(Chapter 5) Determine $[\Delta_{c}(\lambda): L_c(\mu)]$ for $\lambda ,\mu \in \irr(G)$. (This has been carried out for $G=S_n$ by \cite{RouquierqSch} and \cite{suz}; see \cite[Section 2]{yvonne} for a conjecture in the case $G = G(\ell ,1 ,n)$.)

\smallskip
\noindent
{\bf 6.} (Chapter 5)
Suppose that $c\in \C[\mathcal{S}]^{\ad G}$ and let $d\in \Z[\mathcal{S}]^{\ad G}$. Is there an equivalence between $D^b(\mathcal{O}_c)$ and $D^b(\mathcal{O}_{c+d})$?

\smallskip
\noindent
{\bf 7.}
(Chapter 5) Give an explicit algebraic construction of $P_{\KZ_c}$ for any $c\in \C[\mathcal{S}]^{\ad G}$. 

\smallskip
\noindent
{\bf 8.} (Chapter 5)
Is $\mathcal{O}_c$ Koszul? If so, what is its Koszul dual?

\smallskip
\noindent
{\bf 9.}
(Chapter 6) Determine $[\overline{\Delta}_c(\lambda): \overline{L}_c(\mu)]$ for $\lambda,\mu \in \irr(G)$.

\smallskip
\noindent
{\bf 10.}
(Chapter 6) Give a representation theoretic construction of a symplectic resolution to $\h\times \h^*/G_4$. Does it have a hyper-K\"ahler structure?

\smallskip
\noindent
{\bf 11.}
(Chapter 8) Find natural filtrations on objects in $\Lmod{U_{1,c}}$ or $\Lmod{H_{1,c}}$.
 
\smallskip
\noindent
{\bf 12.}
(Chapter 7 and 8) Use $eH_{1,c}(\Gamma\wr S_n)$ to confirm the generalised $n!$ conjecture \cite[Conjecture 7.2.19]{haimansurvey} or give a new proof of the $n!$ theorem.

\bibliographystyle{plain}

\end{document}